\documentclass{article}

\usepackage{amsmath}
\usepackage{amssymb}
\usepackage{stmaryrd}
\usepackage{epsfig}
\usepackage{psfrag}
\usepackage[matrix, arrow, curve]{xy}
\usepackage{myarrows}
\usepackage{theorem}
\usepackage{graphicx}
\usepackage{mathrsfs}
\usepackage{comment}
\theorembodyfont{\rmfamily}

\newtheorem{definition}{Definition}[section]

\theorembodyfont{\it}

\newtheorem{theorem}[definition]{Theorem}%\includegraphics[]{orb-14-4.pdf}
\newtheorem{proposition}[definition]{Proposition}
\newtheorem{lemma}[definition]{Lemma}

\newtheorem{corollary}[definition]{Corollary}

\theorembodyfont{\rmfamily}
\newtheorem{example}[definition]{Example}
\newtheorem{remark}[definition]{Remark}

\DeclareMathOperator{\aut}{aut}
\DeclareMathOperator{\coeq}{coeq}
\DeclareMathOperator{\cof}{cof}
\DeclareMathOperator{\colim}{colim}
\DeclareMathOperator{\cone}{cone}
\DeclareMathOperator{\cov}{Cov}
\DeclareMathOperator{\hocolim}{hocolim}
\DeclareMathOperator{\holim}{holim}
\DeclareMathOperator{\fib}{fib}
\DeclareMathOperator{\gauge}{gauge}
\DeclareMathOperator{\map}{map}
\DeclareMathOperator{\Map}{\mathbb{M}ap}
\DeclareMathOperator{\stackMap}{{\mathscr M\!\mathit a\mathit p}}
\DeclareMathOperator{\Hom}{\mathbb{H}om}

\DeclareMathOperator{\tmap}{map}
\DeclareMathOperator{\pair}{pair}

\DeclareMathOperator{\Orb}{Orb}
\DeclareMathOperator{\sk}{sk}

\def\C{\mathscr C}
\def\D{\mathscr D}

\def\F{\mathscr F}

\def\M{\mathscr M}

\def\X{\mathscr X}
\def\Y{\mathscr Y}
\def\Z{\mathscr Z}

\def\FF{\mathit F}

\def\GG{\mathbb G}
\def\HH{\mathbb H}

\def\LL{\mathbb L}

\def\RR{\mathbb R}

\def\ZZ{\mathbb Z}

\def\cB{\mathbb B}
\def\cE{\mathbb E}
\def\cG{\mathbb G}
\def\cH{\mathbb H}
\def\cK{\mathbb K}
\def\cO{\mathbb O}
\def\cP{\mathbb P}
\def\cU{\mathbb U}

\def\ccF{\mathcal F}
\def\ccC{\mathcal C}

\def\caseone{{\rm$\put(2.75,-.3){\small 1}\bigcirc$}}
\def\casetwo{{\rm$\put(2.9,-.2){\small 2}\bigcirc$}}

\def\gpd{\mathrm{gpd}}

\def\Homone{\Hom^{\put(2.2,.3){\rm\tiny 1}\bigcirc}}
\def\Homtwo{\Hom^{\put(2.1,.3){\rm\tiny 2}\bigcirc}}
\def\gpdone{\Map^{\put(2.2,.3){\rm\tiny 1}\bigcirc}}
\def\gpdtwo{\Map^{\put(2.1,.3){\rm\tiny 2}\bigcirc}}
\def\mapone{\map^{\put(2.2,.3){\rm\tiny 1}\bigcirc}}

\def\half{\textstyle{\frac{1}{2}}}

\def\Id{\mathrm{Id}}
\def\torb{{\Orb'}}

\def\pr{\mathrm{pr}}
\def\proof{\noindent {\em Proof.}\,\,}
\def\qed{\hfill $\square$}
\def\sk{\mathrm{sk}}
\def\top{\text{\rm top}}

\def\dontshow#1{}

\begin{document}

\title{Homotopy Theory of Orbispaces}
\author{David Gepner and Andr\'e Henriques}
%\date{}
\maketitle

\tableofcontents

\section{Introduction}

\subsection{Background \rm\dontshow{sec:b}}\label{sec:b}

An orbispace is supposed to be something that is locally the quotient of a space by a group.
Different meanings of the words {\em space} and {\em group} thus lead to different notions of orbispace,
which have been introduced at different times in the various fields of mathematics.
Table 1 summarizes the most important ones.

\vspace{.3cm}
\centerline{
\begin{minipage}{15.2cm}
\begin{center}
\begin{tabular}{|c|c|c|}
\hline
\phantom{$\Big($}``space'' & ``group'' & ``orbispace''\\
\hline
\hline
\phantom{$\Big($}manifold & finite group & orbifold\\
\hline
\phantom{$\Big($}polyhedron & discrete group & complex of groups\\
\hline
\phantom{$\Big($}algebraic variety & finite group & Deligne-Mumford stack\\
\hline
\phantom{$\Big($}algebraic variety & algebraic group & Artin stack\\
\hline
\phantom{$\Big($}topological space & topological group & topological stack\\
\hline
\end{tabular}
\end{center}
\centerline{Table 1.}
\end{minipage}
}
\vspace{.3cm}

\noindent Several formalisms have been developed to describe these various notions.
Some are specific to the situation for which they have been introduced, while others generalize without difficulty.

Orbifolds were first introduced by Satake \cite{Sat56} and later popularized by Thurston \cite{Thu79}.
Their approach is very geometric, but unfortunately rather tricky to adapt outside of the context of smooth manifolds.
Complexes of groups were introduced by Haefliger \cite{Hae91}; a good reference is \cite[Chapter III.$\mathcal{C}$]{BH99}.
His definiton is purely combinatorial and thus inherently specific to the world of polyhedra and simplicial complexes.

A geometric approach to orbispaces that works well independently of the ambient category is via {\em groupoid objects}.
It has been used both in algebraic geometry and differential geometry; for good surveys, we recommend \cite{Moe02} and \cite{MM05}.
However, the most powerful formalism is undoubtedly that of {\em stacks}, largely due to ideas of Grothendieck.
Stacks were introduced by Deligne and Mumford in \cite{DM69} and later generalized by Artin \cite{Art74}.
A good introduction to stacks in the topological context is provided by Noohi \cite{Noo05}.
The theory of stacks admits in turn two main variants, the first of which uses {\em fibered categories} and is developed in detail in \cite{LMB00}, while the second uses {\em sheaves of groupoids}.
We give a brief and very informal introduction to these two formalisms.

\vspace{.3cm}

\noindent{$\bullet$ \it Groupoid objects:}
An orbispace $\X$ always admits a surjective map $U\to \X$ from an ordinary space $U$.
We may then form the fibered product $V:=U\times_\X U$.
It comes with two projections maps $V\rrarrow U$ which we shall call ``source'' and ``target'', and
we call the map $V\times_U V\to V$
given by projection onto the first and last factors $V\times_U V\simeq U\times_\X U\times_\X U\to U\times_\X U$ by the name ``multiplication''.
Then $V\rrarrow U$ becomes a {\em groupoid object} with $U$ as space of objects and $V$ as space of arrows.
The idea of this approach is that $\X$ is completely determined by this groupoid object, which acts as an atlas for $\X$.
One then defines an orbispace simply as a topological groupoid\footnote{
Some people define orbispaces as {\em equivalence classes} of groupoid objects,
where two groupoid objects are identified if they are Morita equivalent.
We warn the reader that this approach makes it impossible to correctly set up the notion of map between orbispaces.
}.

\vspace{.3cm}
\noindent{$\bullet$ \it Stacks:}
For each space $U$ we can consider the set of maps from $U$ to the orbispace $\X$.
The idea is that the collection of all these maps is enough to encode $\X$.
This may be done in two different ways:
\begin{itemize}
\item Define a category $\mathscr C$ whose objects are maps $U\to\X$ and whose arrows are commutative (in the $2$-categorical sense) triangles
\[
\xymatrix@R=.6cm{
V\ar[rr]\ar[dr]&&U\ar[dl]\\ &\X.
}
\]
We then have a forgetful functor $F$ from $\mathscr C$ to spaces given by sending $U\to \X$ to $U$.
A {\em fibered category} is such a pair $(\mathscr C,F)$, and this data fully encodes the orbispace $\X$.

\item The set of maps from $U$ to $\X$ has a natural groupoid structure.
So the assignment
\[
U\mapsto \hom(U,\X)
\]
is a contravariant functor from spaces to groupoids.
Moreover, this functor is actually a {\em sheaf of groupoids} (in the $2$-categorical sense),
and this sheaf determines the orbispace $\X$.
\end{itemize}

\noindent In this paper we use both groupoid objects and sheaves of groupoids as models for orbispaces, and as we make no further mention of fibered categories, we shall refer to sheaves of groupoids simply as stacks.
Since we will be working entirely in the topological context,
we will only be concerned with {\em topological groupoids} and {\em topological stacks}.
We also introduce a new formalism, convenient for homotopy theory, which we call $\Orb$-spaces.

\subsection{Elmendorf's theorem \rm\dontshow{sec:ET}}\label{sec:ET}

By way of motivation, let us briefly recall the basics of the homotopy theory of $G$-spaces.
We first fix a family\footnote{We do not require that $\ccF$ be closed under taking subgroups.}
$\ccF$ of closed subgroups of $G$, closed under conjugation.
Although everything we do is relative to $\ccF$, we generally omit it from the notation.

A $G$-space is said to be cellular
if it can be constructed by successively attaching cells of the form $D^n\times G/H$ for $H\in\ccF$.
And a map $f:M\to N$ is called a weak equivalence
is for every $H\in\ccF$, the induced map between fixed point sets
\(
M^H=\map_G(G/H,M)\to N^H=\map_G(G/H,N)
\)
is a weak equivalence of spaces.

\begin{definition}
The {\em orbit category} $\Orb_G$ of the topological group $G$ (we suppress $\ccF$ from the notation) is the full subcategory of $G$-spaces on the $G$-orbits $O=G/H$ for $H\in\ccF$.
An $\Orb_G$-space is a continuous contravariant functor from $\Orb_G$ to spaces.
\end{definition}

Cellular objects and weak equivalences also exist in the category of $\Orb_G$-spaces.
The latter are the maps that induce weak equivalences upon evaluating at any orbit $O\in\Orb_G$.

\begin{theorem}[Elmendorf's theorem]\label{et's}\dontshow{et's}
The functor
\[
\Phi:\big\{G\text{\rm-spaces}\big\}\,\longrightarrow\, \big\{\Orb_G\text{\rm-spaces}\big\}
\]
sending a $G$-space $M$ to the $\Orb_G$-space $\Phi(M):O\mapsto \map_G(O,M)$ admits a ``homotopy inverse''
\[
\Psi: \big\{\Orb_G\text{\rm-spaces}\big\}\,\longrightarrow\, \big\{G\text{\rm-spaces}\big\}
\]
and natural transformations from $\Psi\Phi$ and $\Phi\Psi$ to the respective identity functors that provide weak equivalences \dontshow{frv}
\begin{equation}\label{frv}
\Psi\Phi(M)\stackrel{\sim}{\longrightarrow} M,\hspace{2cm} \Phi\Psi(X)\stackrel{\sim}{\longrightarrow} X
\end{equation}
for all $G$-spaces $M$ and $\Orb_G$-spaces $X$.
Moreover, if $M$ and $X$ are cellular, then the maps (\ref{frv}) are homotopy equivalences.
\end{theorem}
The above theorem can be interpreted as giving an explicit equivalence of homotopy theories
between the categories of $G$-spaces and $\Orb_G$-spaces.
The original proof \cite{Elm83} restricts to the case in which $G$ is a compact Lie group and $\ccF$ is the family of all closed subgroups of $G$.
A proof that works in full generality can be found in \cite[Section VI.6]{May96}.

\begin{remark}
A modern way of rephrasing Elmendorf's theorem is in terms of a Quillen equivalence
\[
L:\:\big\{\Orb_G\text{-spaces}\big\}\put(7,3){$\longrightarrow$}\put(7,-2){$\longleftarrow$}\hspace{1cm}\big\{G\text{-spaces}\big\}\::R.
\]
The translation from the previous formulation is given by $\Phi=R$ and $\Psi=L\circ\cof$, where $\cof$ is the cofibrant replacement functor in $\Orb_G$-spaces.
The existence of a natural transformation $\Phi\Psi\to 1$ is then due to the exceptional fact that the unit map
$X\to RL(X)$ is an isomorphism whenever $X$ is cofibrant.
\end{remark}

The main goal of this paper is to prove an analog of Elmendorf's theorem for orbispaces.
As in the equivariant case, we first fix a family $\ccF$ of allowed isotropy groups.
The difference is that now, $\ccF$ is a class of arbitrary topological groups, as opposed to being a set of subgroups of a fixed group $G$.
This leads to the notion of a cellular orbispace, by which we mean an orbispace which can be constructed inductively from the empty orbispace by successively attaching
cells\footnote{There's a technical point to be made here:
the operation of attaching cells needs to be made in the category of topological groupoids, and not in that of topological stacks.
This is because the Yoneda functor does not commute with attaching cells.}
 of the form $D^n\times \M_H$ for $H \in \ccF$.
Here $\M_H$ denotes the classifying stack for principal $H$-bundles\footnote{The stack $\M_H$ is denoted $\M_{\cB H}$ in the main body of this paper.};
it's an orbit stack in the sense that its coarse moduli space, or ``underlying space'', consists of a single point.

The hope, then, is that the same arguments which work for $G$-spaces would apply verbatim in this more general setting of orbispaces.
Namely, we could try to define $\Orb$ as the full subcategory of orbispaces on the orbit stacks $\M_G$ with $G$ in $\ccF$ an allowed isotropy group,
$\Orb$-spaces as continuous contravariant functors from $\Orb$ to spaces, and $\Phi$ as the functor which associates to an orbispace $\X$ the $\Orb$-space
\begin{equation}\label{xsd}
\Phi\X(\M_G):=\map(\M_G,\X).
\end{equation}
We could then attempt to show that $\Phi$ induces an equivalence of homotopy theories.

However, making sense of the functor $\Phi$ presents certain difficulties.
The main problem is that in general there is only a {\em stack}, as opposed to a {\em space}, of maps between stacks.
%the mapping space $\map(\Y,\X)$ between two stacks $\Y$ and $\X$ is typically not a space but rather a stack,
In particular, $\Orb$ is not naturally a topological category, and consequently the notion of an $\Orb$-space is unclear.
Moreover, even after replacing $\Orb$ by an appropriately equivalent topological category, the functor $\Phi$ defined in (\ref{xsd}) does not take values in spaces, and therefore cannot possibly be an $\Orb$-space.
We shall address each of these difficulties in turn.

\subsection{Statement of results \rm\dontshow{s:Res}}\label{s:Res}

Fix a family $\ccF$ of allowed isotropy groups, by which we mean an essentially small class of topological groups closed under isomorphisms,
and subject to a mild paracompactness condition (see Section \ref{se:ctg}).
Once again, although everything is done relative to $\ccF$, we shall generally omit it from the notation.
Given a group $G$, let
$\cB G:=(G\rrarrow\ast)$ denote the topological groupoid with trivial object space $\ast$, arrow space $G$, and composition given by multiplication in $G$.

We define $\Orb$ to be the topologically enriched category whose objects are the groups in $\ccF$ and whose space of morphisms from the group $H$ to the group $G$ is the fat geometric realization\dontshow{pdoo}
\begin{equation}\label{pdoo}
\Orb(H,G):=\|\Map(\cB H,\cB G)\|
\end{equation}
of the topological groupoid of maps from $\cB H$ to $\cB G$.
This category is the topological version of the full subcategory of stacks on the orbit stacks $\M_G$,
and the map $\Orb(H,G)\to\map(\M_H,\M_G)$ is indeed a weak equivalence by Lemma \ref{l:vd} and Theorem \ref{wty}.

\begin{definition}
An {\em $\Orb$-space} is a continuous contravariant functor from $\Orb$ to spaces.
\end{definition}

Our main theorems assert that the homotopy theory of cellular stacks (with isotropy in $\ccF$) is equivalent to the homotopy theory of $\Orb$-spaces.
In order to make sense of the above statement, we introduce the intermediate category of topological groupoids an construct functors \dontshow{3Ct}
\begin{equation}\label{3Ct}
\big\{\text{\rm Stacks}\big\}\,\,\stackrel{\,\textstyle\M}{\longleftarrow\!\!\!-}\,\,\big\{\text{\rm Topological groupoids}\big\}\,\,\stackrel{\textstyle\RR}{-\!\!\!\longrightarrow} \big\{\Orb\text{\rm-spaces}\big\}
\end{equation}
to stacks and to $\Orb$-spaces, respectively.
We then argue that both $\M$ and $\RR$ induce an equivalence of homotopy theories, which is to say that both $\M$ and $\RR$ are
fully faithful and essentially surjective in the appropriate homotopical sense.

In order to make sense of these conditions, we require a considerable amount of extra structure on our categories.
Namely, we need notions of {\em weak equivalence}, {\em fibrant object} and {\em cofibrant object}.
The weak equivalences should satisfy the two out of three property, and a map between fibrant and cofibrant objects should be a weak equivalence if and only if it's a homotopy equivalence.
Moreover, the weak homotopy type of the space of maps from a cofibrant source to a fibrant target should be invariant under weak equivalence, in both source and target, provided the source remains cofibrant and the target remains fibrant.
%
%The weak homotopy type of the space of maps out of a fixed cofibrant source should be invariant under weak equivalence of fibrant
%
%Maps out of cofibrant objects should send weak equivalences between fibrant objects to weak equivalences of spaces, and
%maps into fibrant objects should send weak equivalences between cofibrant objects to weak equivalences of spaces.
%And finally, a map between fibrant-cofibrant objects should be a weak equivalence iff it is a homotopy equivalence.
In the presence of appropriate fibrant and cofibrant replacement functors (an additional piece of data), we can then define the {\em derived mapping space}
%$\map'(Y,X)$ to be the space of maps between the cofibrant replacement of $Y$, and the fibrant replacement of $X$
\[
\map'(Y,X):=\map(\cof Y,\fib X)
\]
as the space of maps between the cofibrant replacement of $Y$ and the fibrant replacement of $X$.
Note however, that in the absence of the appropriate replacement functor,
we only have derived mapping spaces when the source is already cofibrant or the target is already fibrant.

Now let $f$ be a functor between such categories, which respects the extra structure
in the sense that it sends weak equivalences to weak equivalences, fibrant objects to fibrant objects, and cofibrant objects to cofibrant objects.
We say that $f$ is an {\em equivalence of homotopy theories} if it is:

\vspace{.3cm}
\noindent $\bullet$
{\em Homotopically fully faithful:} For all pairs of fibrant and cofibrant objects $X$ and $Y$ in the source category,
the induced map on (derived) mapping spaces
$
\map'(Y,X)\to\map'\big(f(Y),f(X)\big)
$
is a weak homotopy equivalence.

\vspace{.3cm}
\noindent $\bullet$
{\em Homotopically essentially surjective:} For any fibrant and cofibrant object $X$ of the target category there is a fibrant and cofibrant object $Y$ of the source category and a weak equivalence between $f(Y)$ and $X$.
\vspace{.3cm}

%{\em homotopically fully faithful} if the induced map on derived mapping spaces
%$
%\map'(Y,X)\to\map'\big(f(Y),f(X)\big)
%$
%is a weak homotopy equivalence for all pairs of fibrant-cofibrant objects $X$ and $Y$ in the source category,
%that $f$ is {\em homotopically essentially surjective} if for any fibrant-cofibrant object $X$ of the target
%category there is a fibrant-cofibrant object $Y$ of the source category and a weak equivalence $f(Y)\simeq X$,
%and, finally, that $f$ is an {\em equivalence of homotopy theories} if it is both homotopically fully faithful
%and homotopically essentially surjective.

\noindent
If one defines the {\em homotopy category} of a category with said structure as the quotient of the full topological subcategory on the
fibrant-cofibrant objects obtained by identifying homotopic arrows, then an equivalence of homotopy theories induces an (ordinary) equivalence of homotopy categories.

The category of $\Orb$-spaces is a topological model category and therefore enjoys all the extra structure we need.
This is unfortunately not the case for the other two categories, at least to our knowledge\footnote{While there are simplicial and/or topological model category structures on topological groupoids as well as topological stacks, they do not model the homotopy theory of interest to us.}.
For example, we don't introduce any cofibrant replacement functors\footnote{
We're not saying that they don't exist. We just don't introduce any.
} in topological stacks or topological groupoids, but only the notion of a cofibrant object.
%As a consequence, the derived mapping space $\map'(Y,X)$ can only be defined if $Y$ is already cofibrant.
Moreover, the category of stacks causes additional problems because is not even topologically enriched.
This is dealt with in Section \ref{sec:hts}, where we explain how to enrich $\text{\rm Stacks}$ over the homotopy category of spaces.

In order to help the reader understand the subtleties involved we have compiled in Table 2 a list of the structure that is and is not
available in each of the three categories we use to model orbispaces.
Note that the existence of fibrant replacement functors in topological stacks and $\Orb$-spaces is trivial as all objects are already fibrant.
Due to the varying amounts of extra structure, we shall have to interpret slightly differently the conditions of homotopy essential surjectivity and homotopy fully faithfulness for the two functors $\M$ and $\RR$.

The homotopical essential surjectivity of $\M$ is more-or-less axiomatic, for we simply declare that the subcategory of relevant stacks is the essential image under $\M$ of the category of cellular groupoids; the fact that this is also the essential image of the {\em fibrant} cellular groupoids is the content of Proposition \ref{ijd} and Corollary \ref{fcic}.
The homotopical full faithfulness of $\M$ is our first main result.

\begin{theorem}\label{T15}\dontshow{T15}
Let $\cH$ be a cellular topological groupoid and $\cG$ be a fibrant topological groupoid.
Then the functor $\M$ induces a weak equivalence
\(
\map(\cH,\cG)\approx \map(\M_\cH,\M_\cG).
\)
\end{theorem}
\vspace{.3cm}

\centerline{
\begin{minipage}{15.2cm}
\begin{center}
\begin{tabular}{| p{1.5cm}|| p{3.8cm}| p{3.75cm}| p{3.75cm}|}
\hline
& \hspace{.2cm} $\phantom{\Big(}$ Topological stacks & \hspace{.07cm} Topological groupoids & \hspace{.9cm} $\Orb$-spaces\\
\hline\hline
$\overset{\phantom{.}}{\text{W}}$eak \hfill \newline
equiva- lences &
Homotopy equivalences; only defined  between cellular stacks (section \ref{sec:hts})&
Isomorphisms on all \hfill \newline  homotopy groups \hfill \newline (section \ref{sec:HG}) &
Weak equivalences upon evaluation on objects of $\Orb$ (section \ref{catorb})\\
\hline
$\overset{\phantom{.}}{\text{C}}$ofibrant objects &
Cellular stacks \hfill \newline i.e. stacks associated to \hspace{.2cm} cellular groupoids
&
Cellular groupoids \hfill \newline (section \ref{se:ctg}) &
Retracts of cellular $\Orb$-spaces (section \ref{catorb}) \\
\hline
$\overset{\phantom{.}}{\text{F}}$ibrant \hspace{.2cm} objects &
All &
Fibrant groupoids \hfill \newline (section \ref{se:Fr}) &
All\\
\hline
$\overset{\phantom{.}}{\text{R}}$eplace- ment functors &
No cofibrant replacement &
Fibrant replacement; \hfill \newline
no cofibrant replacement (section \ref{se:Fr})&
Cofibrant \hspace{.25cm} replacement (section \ref{catorb})\\
\hline
$\overset{\phantom{.}}{\text{D}}$erived mapping spaces &
Defined up to homotopy; only defined if source is cellular (section \ref{sec:hts})&
$\map(\HH,\fib\GG)$; \hspace{1cm} only defined if source is cellular &
$\map(\cof Y,X)$\\
\hline
\end{tabular}
\end{center}
\centerline{Table 2.}
\end{minipage}
}
\vspace{.3cm}

To show the homotopical essential surjectivity and homotopical full faithfulness of $\RR$ we exhibit it as the right derived functor $\RR=R\circ \fib$ of a functor $R$.
The latter has a left adjoint $L$, which in turn has a left derived functor $\LL=L\circ\cof$.
We then show at once that $\RR$ is homotopically essentially surjective (\ref{30}) and homotopically fully faithful (\ref{32}),
that $\LL$ is homotopically essentially surjective (\ref{29}) and homotopically fully faithful (\ref{31}),
and that $\LL$ is both a homotopy left inverse (\ref{29}) as well as a homotopy right inverse (\ref{30}) of $\RR$.

\begin{theorem}\label{T16}\dontshow{T16}
Let $X$, $Y$ be $\Orb$-spaces and $\cG$, $\cH$ topological groupoids with $\cH$ cellular.
Then the derived unit and counit natural transformations
\begin{eqnarray}
\cof&\longrightarrow&\RR\LL\label{30}\\
\LL\RR&\longrightarrow&\fib\label{29}
\end{eqnarray}
of the derived adjunction
$
\LL:\:\big\{\Orb\text{-spaces}\big\}\put(7,3){$\longrightarrow$}\put(7,-2){$\longleftarrow$}\hspace{1cm}\big\{\text{Topological groupoids}\big\}\::\RR
$
are objectwise weak equivalences, and they induce natural weak equivalences
\begin{eqnarray}
\map(\cof Y,X)&\approx&\map(\LL Y,\fib\LL X)\label{31}\\
\map(\HH,\fib\GG)&\approx&\map(\cof\RR\cH,\RR\cG)\label{32}
\end{eqnarray}
on derived mapping spaces.
\end{theorem}

\begin{remark}
Strictly speaking, as $\RR$ need not preserve cofibrant objects, it doesn't make sense to say that it's an equivalence of homotopy theories
(at least in the above sense).
This problem is easily fixed by replacing $\RR$ by $\cof\circ\hspace{.05cm} \RR$,
the latter being a functor which preserves weak equivalences, fibrant objects, cofibrant objects, and which is both
homotopically fully faithful and homotopically essentially surjective.
Another solution is to consider the functor $\LL$. It also has all the listed properties,
and thus is an equivalence of homotopy theories in the above sense.
We choose to ignore that issue, as it does not further interfere with our results.
\end{remark}

We conclude with a couple of applications to equivariant homotopy theory.
In order to obtain a meaningful comparison, we must first restrict to the subcategory of orbispaces and
{\em representable} maps, which is to say those maps which induce closed monomorphisms on all isotropy groups.
This restriction comes from the fact that $G$-equivariant maps are necessarily injective on isotropy.

Given a topological group $G$ and a family $\ccF$ of closed subgroups, closed under conjugation and taking subgroups, one can form \cite{Lue05}
a $G$-space $E_\ccF G$ characterized by the property that $(E_\ccF)^H$ is contractible for $H\in\ccF$ and empty otherwise.
Let $B_\ccF G:=[E_\ccF G/G]$ denote the quotient stack.
One then gets an equivalence of homotopy theories
\[
\Big\{\parbox{5.7cm}{Topological stacks $\X$ equipped with a representable map $\X\to B_\ccF G$}\Big\}
\quad\approx\quad
\Big\{\parbox{4.6cm}{$G$-spaces with stabilizers in $\ccF$}\Big\}.
\]
This is the precise sense in which one recovers the homotopy theory of $G$-spaces from that of topological stacks.
Much more surprising is that one can also go the other way round and recover the homotopy theory of topological stacks (and representable maps thereof) from that of $G$-spaces.

As before, we fix a class $\ccF$ of allowed isotropy groups and consider the homotopy theory of topological stacks with isotropy in $\ccF$
and representable maps thereof.
Theorems \ref{T15} and \ref{T16} still apply in this context, provided we modify the right hand side of (\ref{pdoo}) by restricting to the full subgroupoid $\Map^\mathrm{rep}(\cB H,\cB G)$ of $\Map(\cB H,\cB G)$ on those maps which induce closed inclusions of $H$ into $G$.
We then have the following result.

\begin{theorem}
Let $\ccF$ be as above, and $\Orb^\mathrm{rep}$ be the topological category with objects $\ccF$ and morphisms
\[
\Orb^\mathrm{rep}(H,G):=\big\|\,\mathbb{M}\mathrm{ap}^{\mathrm{rep}}(\cB H,\cB G)\,\big\|.
%\big\|\,\mathbb{R}\mathrm{epresentable}\mathbb{M}\mathrm{aps}\,(\cB H,\cB G)\,\big\|.
\]
Then there is a topological group $G$, depending on the chosen family $\ccF$,
such that the topological category $\Orb_{G}$ of $G$-orbits with isotropy in $\ccF$ is weakly equivalent to $\Orb^\mathrm{rep}$.
\end{theorem}
Combining Elmendorf's Theorem \ref{et's} and our Theorems \ref{T15} and \ref{T16} (adapted to the case of representable maps),
we thus get an equivalence of homotopy theories
between $G$-spaces with isotropy in $\ccF$, and topological stacks with isotropy in $\ccF$ and representable maps thereof.

\subsection{A convention concerning representable maps}

Recall that a map of orbispaces is said to be {\em representable} or {\em faithful} if it induces closed inclusions on all isotropy groups.
When constructing orbispaces by attaching cells, it is natural to require that the attaching maps be representable; this is the condition which ensures that the resulting orbispace is locally isomorphic to a quotient stack.

On the other hand, when considering maps between orbispaces, there's no reason a priori to restrict to representable maps.
It is therefore tempting to allow arbitrary maps between orbispaces but force attaching maps to be faithful.
Unfortunately, this does not work well from a homotopy theoretic point of view, for forming homotopy colimits of diagrams of cellular objects would tend to take us outside of the realm of cellular objects.
%From the point of view of homotopy theory, especially that of cellular approximation, we should
%it rather unnatural to allow all maps between orbispaces but to restrict those which can be used as attaching maps for cells.

Thus, we are confronted with the following dilemma:
\begin{enumerate}
\item[\caseone]
We allow all maps between orbispaces, at the expense of having to accept as cellular certain orbispaces which are not local quotients.
%But then, we have to accept as cellular many orbispaces which are not local quotients.
\item[\casetwo]
We require our attaching maps be representable, at the expense of losing the ability to consider non-representable maps.
%But then, we lose the possibility to speak about non-representable maps.
\end{enumerate}
As both setups have their advantages and disadvantages, we decided to treat the two cases in parallel.
Thus, we develop cases \caseone \- and \casetwo \- simultaneously throughout the paper
%, leaving the reader free to decide which notion he or she prefers,
%While we could (and perhaps should, for the sake of simplicity) have written the paper entirely from one or the other perspective,
%we see no reason to limit the scope of our theory,
and explicitly point out the few instances in which the distinction is important.

\section{Topological Groupoids \rm\dontshow{s:TG}}\label{s:TG}

\subsection{Conventions and examples \rm\dontshow{s:VX}}\label{s:VX}

Define the category of {\em topological groupoids} in one of the following two ways:
\begin{enumerate}
\item[\caseone] The category of topological groupoids and continuous functors between them.
\item[\casetwo] The category of topological groupoids and continuous {\em faithful}
functors, which is to say those that induce closed inclusions on automorphism groups of objects.
\end{enumerate}

The reader is free to decide which one of the above interpretations he or she prefers so long as it used consistently throughout the course of this paper.

Most of the categories we use are naturally enriched over the category of topological spaces\footnote{
As usual, we work with compactly generated Hausdorff spaces, in order to have the convenient adjunction between products and mapping spaces \cite{Ste67}.},
and we shall need to distinguish between enriched mapping spaces and ordinary morphism sets.
If $X$ and $Y$ are objects of such a category, we shall write $\map(Y,X)$ for the {\em space} of morphisms from $Y$ to $X$, and $\hom(Y,X)$ for the underlying {\em set} of morphisms from $Y$ to $X$.

In addition, some of our categories are naturally $2$-categories, or even topological $2$-categories.
They will always have invertible $2$-morphisms, meaning that they are enriched over (possibly topological) groupoids.
Given objects $X$ and $Y$ of a topological $2$-category, we shall write $\Map(Y,X)$ and $\Hom(Y,X)$ for the topological groupoid and underlying ordinary groupoid, respectively, of morphisms from $Y$ to $X$.
A summary of these notational conventions is provided in Table 3.
As an illustration, if $\cH$ and $\cG$ are topological groupoids, the symbol $\map(\cH,\cG)$ will be used to denote the space of continuous (possibly faithful) functors $\cH\to\cG$; it is the object space of the topological groupoid $\Map(\cH,\cG)$.

\vspace{.3cm}
\centerline{
\begin{minipage}{15.2cm}
\begin{center}
\begin{tabular}{| l || c | c |}
\hline
& \phantom{\Big(}Non-topologized\phantom{\Big(} &\quad\phantom{\Big(} Topologized\phantom{\Big(\quad}\\ \hline\hline
Set of maps\phantom{\Big(} & $\hom\,(\,\text{-}\,,\,\text{-}\,)$ & $\map\,(\,\text{-}\,,\,\text{-}\,)$ \\ \hline
Groupoid of maps\phantom{\Big(} & $\Hom\,(\,\text{-}\,,\,\text{-}\,)$ & $\Map\,(\,\text{-}\,,\,\text{-}\,)$ \\ \hline
\end{tabular}\\
\end{center}
\centerline{Table 3.}
\end{minipage}
}
\vspace{.3cm}

\noindent If we need to specify that we are in one of the two aforementioned cases \caseone, \casetwo, we shall write a superscript $\gpdone$, $\gpdtwo$, $\mapone$, etc.
For example, the adjunction map
\[
\Map(\cK\times\cH,\cG)\to \Map\big(\cK,\Map(\cH,\cG)\big)
\]
is only an isomorphism in case \caseone;
in the other case, we just get inclusions
\[
\gpdtwo(\cK\times\cH,\cG)\hookrightarrow \gpdtwo\big(\cK,\gpdtwo(\cH,\cG)\big)\hookrightarrow\gpdone(\cK\times\cH,\cG).
\]

Given a topological group $G$ we shall write $\cB G$ for the topological groupoid with one object and morphism space $G$.
Note that
\(
\map(\cB H,\cB G)
\)
is none other than the space of continuous homomorphisms from $H$ to $G$ in case \caseone, or the subspace of closed monomorphisms in case \casetwo.

The following examples introduce certain well known classes of groupoids which play an important role in the theory of orbispaces.
See \cite[Section 2.1]{MM05} for the corresponding examples in the context of Lie groupoids.

\begin{example}\label{extp}\dontshow{extp}
The {\em unit groupoid} of a topological space $T$ has both object and arrow spaces equal to $T$, and all groupoid structure maps are the identity of $T$
(provided we identify the fibered product $T\times_T T$ with $T$).
In other words, given $x,y\in T$ there is an arrow $x\to y$ if and only if $x=y$, and in that case the arrow is unique.
In practice it will be clear from the context whether we mean the space $T$ or the unit groupoid associated to the space $T$ and we will not employ different notation for the latter.
Note that the unit groupoid functor is left adjoint to the forgetful functor which associates to a topological groupoid its space of objects.
\end{example}

\begin{example}
The {\em pair groupoid} $\pair(U)$ of a topological space $U$ has object space $U$ and arrow space $U\times U$ with source and target maps the two projections $U\times U\rrarrow U$.
It has the property that between any two objects $x,y\in U$ there is exactly one arrow $x\to y$.
The pair groupoid functor is right adjoint (in case \caseone) to the forgetful functor which associates to a topological groupoid its space of objects.
\end{example}

\begin{example}
The pair groupoid of the set $[n]=\{0,1,\ldots,n\}$ will be denoted $\Delta^n_\gpd$.
We think of $\Delta^n_\gpd$ as a groupoid version of the $n$-simplex with vertices the objects and edges the arrows.
A sequence $x_0\leftarrow\cdots\leftarrow x_n$ of $n$ composable arrows of a groupoid $\cG$ extends uniquely to a functor $\Delta^n_\gpd\to\cG$.
%This means that $\Delta^n_\gpd$ classifies $n$-tuples of composable arrows (in either case, since pair groupoids only have identity automorphisms).
We shall write %\dontshow{Gn}
\begin{equation*}%\label{Gn}
\cG_n:=\map\big(\Delta^n_\gpd,\cG\big)\cong\cG_1\times_{\cG_0}\cdots\times_{\cG_0}\cG_1
\end{equation*}
for the space of such sequences of composable arrows.
\end{example}

\begin{example}
We also have relative versions of pair groupoids.
Applying the same construction in the category of spaces over a base space $T$, the {\em relative pair groupoid} $\pair_T(U)$ associated to a space $U\to T$ over $T$ has object space $U$ and arrow space $U\times_T U$.
If $U$ is an open cover of $T$ admitting a partition of unity then the geometric realization of the nerve of $\pair_T(U)$ is homotopy equivalent to $T$ itself.
Since we will restrict to paracompact spaces, this gives one sense in which topological groupoid $\pair_T(U)$ is equivalent to (the unit groupoid of) the topological space $T$.
\end{example}

\begin{example}\label{txh}\dontshow{txh}
The following is a generalization of the relative pair groupoid construction.
Given a topological groupoid $\GG$ and a space $U\to\GG_0$ over $\cG_0$, the {\em restriction} $\cG_U$ of $\cG$ to $U$
has object space $U$ and arrow space $U\times_{\cG_0}\cG_1\times_{\cG_0}U$.
If $\cG\cong T$ is a unit groupoid then $\cG_U$ is just $\pair_T(U)$.
Again, if $U\to\cG_0$ is a cover then $\cG_U$ should be regarded as being equivalent to $\cG$.
We will later make this precise using the language of topological stacks.
\end{example}

The next few examples involve groupoids constructed from the action of a topological group on a topological space.

\begin{example}
Let $G$ be a topological group and let $X$ be a $G$-space.
The {\em action groupoid} $G\ltimes X$ has $X$ as space of objects and $G\times X$ as space of arrows.
An arrow from the object $x$ to the object $y$ is just a group element $g\in G$ such that $g\cdot x = y$.
We can also describe $G\ltimes X$ as the quotient \dontshow{roq}
\begin{equation}\label{roq}
G\ltimes X\cong\big(\pair(G)\times X\big)\big/G.
\end{equation}
This is a groupoid version of the homotopy quotient.
%, which is recovered taking the quotient of the geometric realization of $\pair(G)\times X$.
Note that the one-object groupoid $\cB G=G\ltimes\ast$ is a special case of an action groupoid.
\end{example}

Given two actions $H\acts Y$ and $G\acts X$, a homomorphism $\varphi:H\to G$, and an $H$-equivariant map $f:Y\to X$ where $H$ acts on $X$ via $\varphi$, there is an obvious map of action groupoids $\varphi\ltimes f:H\ltimes Y\to G\ltimes X$.
However, not all groupoid morphisms $H\ltimes Y\to G\ltimes X$ arise this way:
identifying $G\ltimes G$ with $\pair(G)$ in the obvious way, we see that any map $H\to G$ (not necessarily a group homomorphism) induces a morphism of groupoids
$H\ltimes H\to G\ltimes G$.

\begin{example}\label{Ggd}\dontshow{Ggd}
The {\em gauge groupoid} $\gauge_G(P)$ of a principal $G$-bundle $P\to X$ has object space $X$ and arrow space $(P\times P)/G$.
An arrow in $\gauge_G(P)$ from the object $x\in X$ to the object $y\in X$ can also be described as a $G$-equivariant map $P_x\to P_y$ from the fiber of $P$ over $x$ to the fiber of $P$ over $y$.
An equivalent description of the gauge groupoid of $P\to X$ is as the quotient of $\pair(P)$
\[
\gauge_G(P)\cong\pair(P)/G.
\]
by the free action of $G$.
We shall often drop $G$ from the notation and simply write $\gauge(P)$ for $\gauge_G(P)$.
\end{example}

Gauge groupoids can also be characterized abstractly as those groupoids which are connected in a certain categorical sense.

\begin{lemma}
A groupoid $\cG$ is isomorphic to $\gauge_G(P)$ for some group $G$ and principal $G$-bundle $p:P\to X$ if and only if the product of the source and target maps
$s\times t:\cG_1\to \cG_0\times \cG_0$ admits local sections.
\end{lemma}

\proof
If $\cG=\gauge(P)$, then local sections $U\to P$, $V\to P$ of $p$ defined around points $x,y$ of $X$ provide
a local section $U\times V\to P\times P\twoheadrightarrow (P\times P)/G$ of $s\times t$ around the point $(x,y)$ of $X\times X=\cG_0\times \cG_0$.
Conversely, if $s\times t:\cG_1\to \cG_0\times \cG_0$ admits local sections, then the choice of any object
$x\in\cG_0$ gives a principal $\aut(x)$-bundle $P:=s^{-1}(x)\to\cG_0$ such that $\cG\cong\gauge_{\aut(x)}(P)$.
\qed

\begin{lemma}\label{jspq}\dontshow{jspq}
Let $G$ be a topological group and $H$ a closed subgroup with the property that $G\to G/H$ is a locally trivial principal $H$-bundle.
Then there is a natural isomorphism of groupoids $G\ltimes(G/H)\cong \gauge_H(G)$.
\end{lemma}

\proof
Both $G\ltimes(G/H)$ and $\gauge_H(G)$ have $G/H$ as object space.
An arrow $g_1H\to g_2H$ in the former groupoid is an element $g\in G$ satisfying $gg_1H=g_2H$, while an arrow in the latter groupoid is an $H$-equivariant map $\varphi:g_1H\to g_2H$.
These two data are equivalent to each other via the formulae $\varphi(x)=gx$ and $g=\varphi(g_1)g_1^{-1}$.
\qed

\subsection{Cellular topological groupoids \rm\dontshow{se:ctg}}\label{se:ctg}

In order to avoid set-theoretic issues without limiting the flexibility of our theory, let us choose, once and for all,
an arbitrary but fixed family of topological groups $\ccF$, which we refer to as our family of {\em allowed isotropy groups}.
We require that $\ccF$ is closed under isomorphism and that its isomorphism classes of objects be indexable by a set.
Recall that a space is called paracompact if every open cover can be refined to a locally finite open cover, and that
this condition is equivalent to the existence of enough partitions of unity \cite{Ern53} (see \cite[Section 6.4]{Mun75} for an introduction to paracompactness).
We also add the technical condition that the groups $G$ in $\ccF$, as well as their finite powers $G^n$, be paracompact spaces\footnote{
We don't know if this condition is implied by the simple fact that $G$ is paracompact (we believe it's not). \label{foot}}.
Some examples we have in mind are the families of finite groups, countable discrete groups, and compact Lie groups.

\begin{remark}
Our paracompactness condition, though annoying to state, is quite weak in practice.
For example, it is satisfied if the groups $G$ are metrizable or if they are CW-complexes.
It is a special case of the following more general notion of paracompactness for topological groupoids.
\end{remark}

\begin{definition} \label{d:pc} \dontshow{d:pc}
We shall say that a topological groupoid $\cG$ is {\em paracompact} if all the spaces $\cG_n:=\map(\Delta^n_\gpd,\GG)$ are paracompact.
\end{definition}

Recall that a $G$-cell complex is a $G$-space built from attaching $G$-cells, that is, $G$-spaces of the form $D^n\times O$, where $O=G/H$ is a $G$-orbit with allowed isotropy group.
We shall define cell groupoids analogously,
as the ones that can be build by attaching cells of the form $D^n\times \cO$, where $\cO$ is an {\em orbit groupoid}.
Examples of orbit groupoids include the groupoids $\cB H$ for $H \in \ccF$,
as well as the action groupoids $G\ltimes(G/H)\cong\gauge_H(G)$ associated to $G$-orbits.

Anticipating Proposition \ref{BGstr}, we may note that $\cB H$ and $G\ltimes(G/H)$ become equivalent upon passage to associated stacks, and also that there are many other groupoids sharing this same property.
As will be shown in Proposition \ref{BGstr}, these are exactly the groupoids of the form $\gauge_H(P)$.
We'll take the point of view that all of them deserve to be called orbits (at least if paracompact).

\begin{definition}\label{celgr}\dontshow{celgr}
An orbit groupoid is a paracompact topological groupoid of the form $\gauge_G(P)$ for $G$ an allowed isotropy group and $P$ a principal $G$-bundle.
A cell groupoid is a topological groupoid of the form $D^n\times\cO$, where $\cO$ is an orbit groupoid.
\end{definition}

The cell groupoids are the basic building blocks for cellular groupoids.

\begin{definition}
A cellular groupoid is a topological groupoid built from attaching cell groupoids.
That is, it's a topological groupoid of the form $\colim \cG_\alpha$,
where $\alpha$ runs through the set of ordinals smaller than a given ordinal,
such that, if $\alpha-1$ exists, then $\GG_\alpha$ is obtained from $\GG_{\alpha-1}$ as the pushout
\[
\begin{matrix}\xymatrix{
S^{n-1}\times
\cO
\ar[r]\ar[d]^f&D^{n}\times
\cO
\ar[d]\\
\cG_{\alpha-1}\ar[r]& \cG_{\alpha}}\end{matrix}
\]
of a cell groupoid $D^{n}\times
\cO
$ along an attaching map $f\in\hom(S^{n-1}\times
\cO
,\GG_{\alpha-1})$, and if $\alpha-1$ does not exist, then $\cG_\alpha$ is obtained as the colimit $\underrightarrow{\lim}_{\beta<\alpha}\cG_\beta$.
\end{definition}

\begin{remark}
We emphasize that, whereas in case \caseone, the attaching map $f$ is an arbitrary morphism of groupoids, in case \casetwo, we require that $f$ be faithful.
Thus the classes of cellular groupoids differ considerably between cases \caseone\phantom{,\!\!} and \casetwo.
\end{remark}

\begin{remark}
If one goes a step further and defines a cofibrant groupoid as a retract of a cellular groupoid, then most (if not all) of our assertions about cellular groupoids remain true for this more general class.
%We avoid this terminology and instead reserve the term for the cofibrant objects in the model category of $\Orb$-spaces.
\end{remark}

Note that if $G$ is a topological group (which needn't be cellular as a space) and $X$ is a cellular $G$-space with stabilizers in $\ccF$, then the action groupoid $G\ltimes X$ is a cellular groupoid, at least provided the paracompactness condition is satisfied (which it almost always is in practice).
This allows us to view as cellular certain groupoids which arise from actions of infinite-dimensional groups.

\begin{example}
Let $G$ be a compact Lie group, let $P\to M$ be a principal $G$-bundle over a compact manifold $M$ with universal cover $N\to M$ (a principal $\pi_1(M)$-bundle), and let $C$ be the space of flat connections on $P$.
%Fix a principal $G$-bundle $P$ over a compact manifold $M$, and let $C$ be the space of flat connections on $P$.
The gauge group $C^\infty(M,G)$ acts on $P$ and thus also on $C$.
Fixing a point $p\in P$, the monodromy of a connection provides a homomorphism $\pi_1(M)\to G$; conversely, given a homomorphism $\pi_1(M)\to G$, we get a flat $G$-bundle $N\times_{\pi_1(M)} G$.
%, where $\widetilde M$ denotes the universal cover of $M$.
Thus, letting $X\subset \mapone(\pi_1(M),G)/G$ be the subspace (a union of connected components)
consisting of those homomorphisms such that $N\times_{\pi_1(M)} G\simeq P$,
%From the above discussion, we deduce that
the stack theoretic quotients $[C/C^\infty(M,G)]$ and $[X/G]$ are equivalent.
But $X$ is a $G$-CW-complex, so $C$ is a $C^\infty(M,G)$-CW-complex,
%and the paracompactness condition is satisfied;
and consequently $C^\infty(M,G)\ltimes C$ is cellular (as a groupoid) with respect to the family of compact Lie groups.
\end{example}

In many respects orbit groupoids behave like points.
For instance, mapping spaces from orbit groupoids commute with certain kinds of colimits.
We introduce the following substitute for the notion of a cofibration of topological groupoids to make this precise.

Given a topological groupoid $\cG$,
a subspace $Z\subset \cG_0$ will be called {\em saturated} if there are no arrows from $Z$ to its complement,
which is to say that for each arrow $g:x\to y$ we have $x\in Z\Leftrightarrow y\in Z$.
A {\em saturated inclusion} of topological groupoids is a map of the form $\cG_Z\to\cG$ such that $Z\subset\cG_0$
is a saturated subspace and $\cG_Z$ is the restriction of $\cG$ to $Z$ (see Example \ref{txh}).
If $Z\subset\GG_0$ is closed, then we shall say that $\GG_Z\to\GG$ is a closed saturated inclusion.
For example, if $Z\hookrightarrow Y$ is a (closed) inclusion of spaces and $\cG$ is a groupoid
then $Z\times\cG\hookrightarrow Y\times\cG$ is a (closed) saturated inclusion.
Attaching a cell $D^n\times\cO$ to a groupoid $\cG$ via an attaching map $S^{n-1}\times\cO\to\GG$
thus involves pushing out along the closed saturated inclusion $S^{n-1}\times\cO\to D^n\times\cO$.

\begin{lemma}
The pushout of a (closed) saturated inclusion $\cH_Z\hookrightarrow\cH$ along a map $\cH_Z\to\cG$ is again a (closed) saturated inclusion.
\end{lemma}

\proof
Let $\cK$ be the pushout and $Y$ the image of $\cG_0$ in $\cK_0$.
Then $\cG\cong\cK_Y$, and moreover if $Z$ is closed in $\cH_0$ then $Y\cong \cG_0$ is closed in $\cK_0$ since the latter is the pushout of $\cG_0\leftarrow Z \hookrightarrow \cH_0$.
\qed
\vspace{.3cm}

Let $\pi$ be the functor that assigns to a topological groupoid $\cG$ its {\em space} of isomorphism classes of objects $\pi(\cG):=\mbox{coeq}(\cG_1\rrarrow\cG_0)$.
Then given a saturated inclusion $\cH\hookrightarrow\cG$, the diagram
\[
\xymatrix{
\cH\ar[r]\ar[d]&\cG\ar[d]\\
\pi\cH\ar[r]&\pi\cG
}
\]
is always a pullback square (here we identify $\pi\cH$ and $\pi\cG$ with their corresponding unit groupoids).

\begin{lemma}\label{orfilt}\dontshow{orfilt}
Let $\GG$ be the colimit of a (possibly transfinite) sequence of closed saturated inclusions $\cG^1\hookrightarrow\cG^2\hookrightarrow\cdots$.
Then for any orbit groupoid $\cO$, the natural map
$$
\underrightarrow{\lim}\,\map(\cO,\cG^n)\longrightarrow\map(\cO,\cG)
$$
is an isomorphism of spaces.
\end{lemma}

\proof
Since $\cG^n\hookrightarrow\cG^{n+1}$ are saturated inclusions, we have a sequence of pullback squares
\[
\xymatrix{
\cG^1\ar[r]\ar[d]&\cG^2\ar[r]\ar[d]&\ldots\ar@{}[d]|{\displaystyle\cdots}\ar[r]&\cG\ar[d]\\
\pi\cG^1\ar[r]&\pi\cG^2\ar[r]&\ldots\ar[r]&\pi\cG
}
\]
Thus, for every groupoid $\cH$, the square
\dontshow{amw}
\begin{equation}\label{amw}
\begin{matrix}
\xymatrix{
\underrightarrow{\lim}\,
\map(\cH,\cG^n)\ar[r]\ar[d]&\map(\cH,\cG)\ar[d]\\
\underrightarrow{\lim}\,
\mapone(\cH,\pi\cG^n)\ar[r]&\mapone(\cH,\pi\cG).
}
\end{matrix}
\end{equation}
is also a pullback.
If $\cH=\cO$ is an orbit groupoid, the bottom arrow of (\ref{amw}) is just
\[
\underrightarrow{\lim}\,\mapone(\cO,\pi\cG^n)\cong \underrightarrow{\lim}\,\map(\pi\cO,\pi\cG^n)\longrightarrow \map(\pi\cO,\pi\cG)\cong \mapone(\cO,\pi\cG),
\]
which is an isomorphism since $\pi\cO=*$.
Therefore the top arrow of (\ref{amw}) is also an isomorphism.
\qed\vspace{.3cm}

Given a closed saturated inclusion $\cG\hookrightarrow \cH\sqcup_{\cH_Z}\cG$, obtained from pushing out $\cH_Z\hookrightarrow \cH$ along a map $\cH_Z\to \cG$,
and given an orbit groupoid $\cO$, we had hoped that the continuous bijection
\[
\map(\cO,\cH)\underset{\map(\cO,\cH_Z)}{\sqcup}\map(\cO,\cG)\to\map\Big(\cO,\cH\!\underset{\hspace{.2cm}\cH_Z}{\sqcup}\!\cG\Big)
\]
would turn out to be an isomorphism.
Unfortunately, the topologies really are different; nevertheless, we have the following lemma.

\begin{lemma}\label{orpush}\dontshow{orpush}
Let $\cG$, $\cH$ be two topological groupoids and $S^{n-1}\times \cH\to \cG$ be a map.
Then for any orbit groupoid $\cO$, the map \dontshow{sya}
\begin{equation}\label{sya}
\map(\cO,\cG)\underset{S^{n-1}\times\map(\cO,\cH)}{\sqcup}D^n\times\map(\cO,\cH)\,\longrightarrow\,
\map\Big(\cO,\cG\underset{S^{n-1}\times\cH}{\sqcup}D^n\times\cH\Big)
\end{equation}
is a homotopy equivalence.
Moreover, the inclusion
\(
\map(\cO,\cG)\hookrightarrow\map\big(\cO,\cG
\sqcup_{S^{n-1}\times\cH}
D^n\times\cH\big)
\)
is a Hurewicz cofibration (i.e. a neighborhood deformation retraction).
\end{lemma}

\proof
The neighborhood $\map\big(\cO,\cG
\sqcup_{S^{n-1}\times\cH}
(D^n-\frac{1}{2}D^n)\times\cH\big)$
of $\map(\cO,\cG)$ in
$\map\big(\cO,\cG
\sqcup_{S^{n-1}\times\cH}
D^n\times\cH\big)$
admits a deformation retraction coming from the deformation retraction $(D^n-\frac{1}{2}D^n)\searrow S^{n-1}$.
Extending that retraction via the isomorphism $\map(\cO,\frac{1}{2}D^n\times \cH)\cong\map(\cO,D^n\times \cH)$ provides the desired homotopy inverse of (\ref{sya}).
\qed

\subsection{Geometric realization of a topological groupoid \rm\dontshow{s:gr}}\label{s:gr}

The fat geometric realization \cite{Seg74} of a simplicial space $X$ is the quotient
\[
\|X\|=\coprod_{n}\Delta^n_\top\times X_n\big/\sim
\]
which results from identifying $(f_* p,x)\in\Delta^n_\top\times X_n$ with $(p,f^*x)\in\Delta^m_\top\times X_m$ for any {\em strictly increasing} map $f:[m]\to [n]$.
Here $\Delta^n_\top$ denotes the standard topological $n$-simplex.
The fat realization has a filtration
$
\|X\|\cong\colim_d\sk_d\|X\|
$
by skeleta such that $\sk_d\|X\|$ is obtained from $\sk_{d-1}\|X\|$ as the pushout \dontshow{skF}
\begin{equation}\label{skF}
\begin{matrix}
\xymatrix{
X_d\times\partial\Delta^d_\top\ar[r]\ar[d] & \sk_{d-1}\|X\|\ar[d]\\
X_d\times        \Delta^d_\top\ar[r]       & \sk_d    \|X\|.}
\end{matrix}
\end{equation}

If the simplicial space $X$ is the nerve of a groupoid, then its fat realization admits the following more explicit description.
A point of $\Delta^n_\top\times \cG_n$ is precisely a composable sequence of morphisms of {\em weighted} objects $r\cdot x$ of $\cG$ \dontshow{obF}
\begin{equation}\label{obF}
r_0\cdot x_0\stackrel{\alpha_1}{\leftarrow}r_1\cdot x_1\stackrel{\alpha_2}{\leftarrow}\cdots\stackrel{\alpha_{n-1}}{\leftarrow} r_{n-1}\cdot x_{n-1}\stackrel{\alpha_n}{\leftarrow}r_n\cdot x_n,
\end{equation}
with $0\leq r_i\leq 1$, $\sum r_i =1$, and
$x_0\stackrel{\alpha_1}{\leftarrow}x_1\leftarrow\ldots\leftarrow x_{n-1}\stackrel{\alpha_n}{\leftarrow}x_n$
a sequence of composable morphisms in $\cG$.
The equivalence relation is generated by  reducing sequences, that is, identifying a sequence in which some $r_i=0$ with the sequence obtained by omitting $r_i\cdot x_i$ and composing $\alpha_i$ with $\alpha_{i+1}$.
If $i=0$ we instead omit $\alpha_i$, and if $i=n$ we instead omit $\alpha_n$.

\begin{remark}
If we further identify a sequence in which some $x_{i-1}\stackrel{\alpha_i}{\leftarrow}x_i$ is an identity arrow with the sequence obtained by omitting $\alpha_i$ and combining $r_{i-1}\cdot x_{i-1}$ with $r_i\cdot x_i$ into $(r_{i-1}+r_i)\cdot x_i$, then we obtain the usual geometric realization $|\cG|$.
The projection $\|\cG\|\to|\cG|$ is a homotopy equivalence provided the identity map $\cG_0\to\cG_1$ is a neighborhood deformation retract.
\end{remark}

\begin{example}
The fat geometric realization of the pair groupoid of a nonempty space is always contractible.
If $G$ is a topological group, the translation action of $G$ on itself induces a proper free action on $\|\pair(G)\|$.
It follows that  $\|\pair(G)\|$ is a model for $EG$.
\end{example}

\begin{example}
The fat geometric realization of an action groupoid $G\ltimes X$ is a model for the homotopy quotient \dontshow{x:hq}
\begin{equation}\label{x:hq}
\|G\ltimes X\|\cong X\times_G \|\pair(G)\|\cong X\times_G EG,
\end{equation}
as can easily be seen from (\ref{roq}).
\end{example}

\begin{remark}
Unlike the usual geometric realization, which commutes with products on the nose, the fat geometric realization only commutes with products
(and other finite limits) up to homotopy.
More precisely, given any two simplicial spaces $X$ and $Y$, the fat realization $\|X\times Y\|$ is a retract of the product $\|X\|\times \|Y\|$ via maps $\iota:\|X\times Y\|\to\|X\|\times\|Y\|$ and $r:\|X\|\times\|Y\|\to\|X\times Y\|$, in such a way that given any three simplicial spaces $X$, $Y$, and $Z$,
the relations\dontshow{iri}
\begin{equation}\label{iri}
\begin{split}
\iota\circ(\iota\times 1)&= \iota\circ(1\times \iota):\|X\times Y\times Z\|\to\|X\|\times\|Y\|\times\|Z\|,\\
r\circ(r\times 1)&= r\circ(1\times r):\|X\|\times\|Y\|\times\|Z\|\to\|X\times Y\times Z\|
\end{split}
\end{equation}
are satisfied.
Note that the fat geometric realization commutes with fibered products of simplicial spaces;
in particular it commutes with all finite limits fiberwise over $\|\!*\!\|$.
In other words, the fat realization functor commutes with finite limits, provided we regard it as taking values in the category of spaces over $\|\!*\!\|$.
\end{remark}

Recall that a topological groupoid $\cG$ is paracompact if all the spaces $\cG_n$ are paracompact.
Since the product of a paracompact space and a compact space is paracompact, the mapping cylinder of a map between paracompact spaces is paracompact,
and the colimit of a (possibly transfinite) sequence of NR inclusions (that is, inclusions of closed subspaces $Z\subset X$ such that $Z$ has a neighborhood $U\subset X$ that retracts onto $Z$) of paracompact spaces is paracompact, we see that:

\begin{lemma}\label{pcrap}\dontshow{pcrap}
Cellular groupoids are paracompact.\qed
\end{lemma}

\begin{lemma}\label{pcannoy}\dontshow{pcannoy}
If $\cG$ is paracompact, then $\|\cG\|$ is paracompact.\qed
\end{lemma}

Recall that a map $p:U\to T$ is said to be {\em \'etale} if every point $u\in U$ has a neighborhood $V\subset U$ such that $V$ is homeomorphic via $p$ to its image in $T$.

\begin{definition} \label{d:cov} \dontshow{d:cov}
A {\em cover} of a topological space $T$ is a surjective \'etale map $U \to T$.
We write $\cov(T)$ for the full subcategory of spaces over $T$ consisting of the covers.
\end{definition}

Typical examples of covers include the projection map $U=\coprod U_\alpha\to T$, where $\{U_\alpha\}$ is an open cover of $T$,
as well as covering spaces (universal covers, for instance).

Unfortunately, the category $\cov(T)$ is not filtered, and thus not well suited for taking colimits of (contravariant) functors.
For that reason, we introduce two new categories $\cov^0(T)$, $\cov^2(T)$ whose objects are the covers of $T$.
Letting $\pi_{-1}$ be given by $\pi_{-1}(X)=*$ if $X\not =\emptyset$, and $\pi_{-1}(\emptyset)=\emptyset$,
we define \dontshow{d2cv}
\begin{equation}\label{d2cv}
\begin{split}
\hom_{\cov^0(T)}(V,U)&:=\pi_{-1}\hom_{\cov(T)}(V,U),\\
\Hom_{\cov^2(T)}(V,U)&:=\pair\,\hom_{\cov(T)}(V,U).
\end{split}
\end{equation}
The category $\cov^0(T)$ has a morphism $V\to U$ if and only if U refines $V$, and in that case the morphism is unique;
$\cov^2(T)$ is then a 2-category which is equivalent to $\cov^0(T)$.
Note that, as opposed to $\cov(T)$, the categories $\cov^0(T)$ and $\cov^2(T)$ are filtered.

\begin{definition} \label{d:pou} \dontshow{d:pou}
A {\em partition of unity} on a topological space $T$ is a triple $(U,\varphi,<)$, where $U=(U,p:U\to T)$ is a cover of $T$ with the property that the inverse image $p^{-1}(t)$ of any point $t\in T$ is finite,
$\varphi:U\to \RR_{>0}$ is a continuous function satisfying $$\sum_{u\in p^{-1}(t)}\varphi(u)=1$$ for all $t\in T$,
and $<$ is a total ordering of each fiber $p^{-1}(t)$ which varies continuously with $t$.
The latter condition means that $\{(u,v)\in U\times_T U\,|\,u<v\}$ is both closed and open in $U\times_T U$.
\end{definition}

%Let $T$ be a topological space and $p:U\to T$ a cover with the property that the inverse image $p^{-1}(t)$ of any point $t\in T$ is a finite set.
%Then a {\em partition of unity} on $T$, subordinate to the cover $p:U\to T$, is a pair $(\varphi,<)$,
%where $\varphi:U\to \RR_{>0}$ is a continuous function satisfying $$\sum_{u\in p^{-1}(t)}\varphi(u)=1$$ for all $t\in T$,
%and $<$ is a total order on each set $p^{-1}(t)$.
%Moreover, we require that $<$ varies continuously with $t$,
%so that $\{(u,v)\in U\times_T U\,|\,u<v\}$ is both closed and open in $U\times_T U$.
%\end{definition}

We shall see that this nonstandard notion of partition of unity is especially well suited to working with fat geometric realizations.
Unlike $|\cG|$, the space $\|\cG\|$ comes with a preferred
%cover $p_{\|\cG\|}:U_{\|\cG\|}\to\|\cG\|$ and
partition of unity $(U_{\|\cG\|},\varphi_{\|\cG\|},<)$.
A point in $U_{\|\cG\|}$ consists of an element
\(
(r_0\cdot x_0\stackrel{\alpha_1}{\leftarrow}\cdots\stackrel{\alpha_n}{\leftarrow}r_n\cdot x_n)
\)
of $\|\cG\|$, together with a specified weighted object $r_i\cdot x_i$ such that $r_i\not = 0$, and will be denoted
\(
(r_0\cdot x_0\stackrel{\alpha_1}{\leftarrow}\cdots\stackrel{\alpha_i}{\leftarrow}\underline{r_i\cdot x_i}\stackrel{\alpha_{i+1}}{\leftarrow}\cdots\stackrel{\alpha_n}{\leftarrow}r_n\cdot x_n)
\).
The map $U_{\|\cG\|}\to\|\cG\|$ forgets the distinguished object, $\varphi_{\|\cG\|}$ is given by
\[
\varphi_{\|\cG\|}:(r_0\cdot x_0\stackrel{\alpha_1}{\leftarrow}\cdots\stackrel{\alpha_i}{\leftarrow}\underline{r_i\cdot x_i}\stackrel{\alpha_{i+1}}{\leftarrow}\cdots\stackrel{\alpha_n}{\leftarrow}r_n\cdot x_n)
\mapsto r_i,
\]
%The fiber of $p_{\|\cG\|}:U_{\|\cG\|}\to\|\cG\|$ over the point
%\[
%(r_0\cdot x_0\stackrel{\alpha_1}{\leftarrow}\cdots\stackrel{\alpha_i}{\leftarrow}r_i\cdot x_i\stackrel{\alpha_{i+1}}{\leftarrow}\cdots\stackrel{\alpha_n}{\leftarrow}r_n\cdot x_n)\in\|\cG\|
%\]
%is the set of $i\in\{0,\ldots,n\}$ such that $r_i\neq 0$, so that a point of $U_{\|\cG\|}$ consists of a point of $\|\cG\|$ and a choice of distinguished {\em nonzero} weighted object, written
%\[
%(r_0\cdot x_0\stackrel{\alpha_1}{\leftarrow}\cdots\stackrel{\alpha_i}{\leftarrow}\underline{r_i\cdot x_i}\stackrel{\alpha_{i+1}}{\leftarrow}\cdots\stackrel{\alpha_n}{\leftarrow}r_n\cdot x_n).
%\]
%The fact that the (nonzero) weights sum to unity determines the map $\varphi_{\|\cG\|}:U_{\|\cG\|}\to\RR_{>0}$,
%\[
%\varphi_{\|\cG\|}:(r_0\cdot x_0\stackrel{\alpha_1}{\leftarrow}\cdots\stackrel{\alpha_i}{\leftarrow}\underline{r_i\cdot x_i}\stackrel{\alpha_{i+1}}{\leftarrow}\cdots\stackrel{\alpha_n}{\leftarrow}r_n\cdot x_n)
%\mapsto r_i,
%\]
and the total order $<$ on the fibers is just the rule that
\[
(r_0\cdot x_0\stackrel{\alpha_1}{\leftarrow}\cdots\stackrel{\alpha_i}{\leftarrow}\underline{r_i\cdot x_i}\stackrel{\alpha_{i+1}}{\leftarrow}\cdots\stackrel{\alpha_n}{\leftarrow}r_n\cdot x_n)<
(r_0\cdot x_0\stackrel{\alpha_1}{\leftarrow}\cdots\stackrel{\alpha_j}{\leftarrow}\underline{r_j\cdot x_j}\stackrel{\alpha_{j+1}}{\leftarrow}\cdots\stackrel{\alpha_n}{\leftarrow}r_n\cdot x_n),
\]
whenever $i<j$.

Let $p_{\|*\|}:U_{\|*\|}\to \|\!*\!\|$ and $(\varphi_{\|*\|},<_{\|*\|})$ be the canonical partition of unity on the fat realization
$\|\!*\!\|$ of the terminal groupoid $*$.
Then any space over $\|\!*\!\|$ also acquires a partition of unity by pullback from $\|\!*\!\|$;
%This construction holds for the fat geometric realization of any simplicial space,
%for it may be obtained simply by pulling back the canonical cover $\dot p:\dot U\to \|\!*\!\|$ and partition of unity $(\dot\varphi,\dot <)$ on the fat realization $\|\!*\!\|$ of the terminal simplicial space.
in fact, the space $\|\!*\!\|$ classifies partitions of unity up to isomorphism.

\begin{lemma}
The assignment which sends the map $f:T\to\|\!*\!\|$ to the
%cover $f^*p_{\|*\|}:f^*U_{\|*\|}\to T$ and
partition of unity
$$
f^*(U_{\|*\|},\varphi_{\|*\|},<_{\|*\|}):=(f^*U_{\|*\|},f^*\varphi_{\|*\|},f^*<_{\|*\|})
$$
induces a bijection between the set of maps from $T$ to $\|\!*\!\|$ and the set of isomorphism classes of partitions of unity $(U,\varphi,<)$ on $T$.
%with $U\to T$ a cover and $(\varphi,<)$ a partition of unity.
\end{lemma}

\proof
A point in $\|\!*\!\|$ is the same thing as a collection $(r_0,\ldots, r_n)$ of positive real numbers satisfying $\sum r_i=1$.
Given a partition of unity $(U,\varphi,<)$ on $T$, the corresponding map $f:T\to\|\!*\!\|$ sends $t\in T$ to the collection of values $\varphi(u)$
as $u$ runs through the ordered set $p^{-1}(t)$.
\qed

\subsection{Groupoid actions and principal groupoid bundles \rm\dontshow{secPGB}}\label{secPGB}

%The following definition is well known, see for example \cite[section 2.4]{MM05}.
\begin{definition}
Let $\cG$ be a topological groupoid and $X\to\GG_0$ a space over $\cG_0$.
A {\em left action of $\cG$ on $X\to\GG_0$} consists of a map $a:\cG_1\times_{\cG_0} X\to X$ over $\cG_0$ which is associative and unital in the sense that the diagrams
$$
\xymatrix{
\cG_1\times_{\cG_0}\cG_1\times_{\cG_0} X\ar[r]^-{m\times 1}\ar[d]_-{1\times a} & \cG_1\times_{\cG_0} X\ar[d]^-a\\
\cG_1\times_{\cG_0} X\ar[r]^-a & X}
\qquad\qquad
\xymatrix{
\cG_0\times_{\cG_0} X\ar[r]^-{i\times 1}\ar[d]_-\cong & \cG_1\times_{\cG_0} X\ar[d]^-a\\
X\ar[r]^-= & X}
$$
commute.
A {\em $\GG$-space} is a space $X\to\GG_0$ over $\GG_0$ equipped with a (left) action of $\GG$.
\end{definition}

\begin{remark}
It is useful to think of a space over $\cG_0$ as a continuous $\cG_0$-indexed family of spaces
%, or the object function of a continuous functor from $\GG$ to spaces,
and a left action of $\cG$ on a space over $\cG_0$ as a continuous contravariant functor from $\cG$ to spaces.
\end{remark}

\begin{example}
Given a topological space $T$, the associated {\em trivial action} of $\cG$ is the action on the space $\cG_0\times T$ given by $a:\cG_1\times_{\cG_0} (\cG_0\times T)\cong\cG_1\times T\put(4.5,5){$\scriptstyle t\times 1$}-\!\!\!\!\longrightarrow\cG_0\times T$.
It corresponds to the notion of the functor from $\cG$ to spaces with constant value $T$.
\end{example}

\begin{example}
The {\em regular action} of $\cG$ on $\cG_1$ is given by $a=m:\cG_1\times_{\cG_0}\cG_1\to\cG_1$.
\end{example}

Given a space $X$ over $\cG_0$ equipped with an action of $\cG$,
the {\em quotient} $X/\cG$ is the coequalizer
$$
\cG_1\times_{\cG_0} X\rrarrow X\rightarrow X/\cG
$$
of the action map and the projection map.
Thinking of an action as a functor from $\cG$ to spaces, we can also describe $X/\cG$ as the colimit of that functor.
\medskip

The following definition is well-known, see for example \cite[Section 2.5]{MM05}.
\begin{definition}
Let $P\to\cG_0$ be a space equipped with an action of $\cG$ and let $T$ be the quotient.
We say that $P\to T$ is a {\em principal $\cG$-bundle} if
$P\to T$ admits local sections and if the map
$$\cG_1\times_{\cG_0} P\to P\times_T P$$
induced by the action and projection maps is a homeomorphism.
\end{definition}

\begin{remark}
If $P$ has a $\cG$-action with quotient $T$ such that the map $P\to T$ has local sections and if for each object
$s\in\cG_0$ the group $\aut(s)$ acts freely on the fiber $P_s$ over $s$,
then the map $\cG_1\times_{\cG_0}P\to P\times_T P$ is automatically a bijection.
The requirement that it be a homeomorphism should be thought of as a rather small additional assumption.
\end{remark}

Note that the total space $P$ of a principal $\cG$-bundle admits two reference maps, one to $\cG_0$ and one to $T$.
The collection of all these maps can be encoded in the following little diagram:
\[
\begin{matrix}
\cG_1\acts P\\
\downarrow\downarrow\hspace{.1cm}\swarrow\hspace{.15cm}\downarrow\hspace{.05cm}\\
\hspace{0cm}\cG_0\hspace{.45cm}T
\end{matrix}
\]

\begin{example}
Given a space $T$ over $\cG_0$, the {\em trivial} principal $\cG$-bundle on $T$ is given by \dontshow{swr}
\begin{equation}\label{swr}
\quad\begin{matrix}
\cG_1\acts \cG_1\times_{\cG_0} T\\
\downarrow\downarrow\hspace{.2cm}\swarrow\hspace{.5cm}\downarrow\hspace{.8cm}\\
\cG_0\hspace{.92cm}T\hspace{.75cm}
\end{matrix}
\end{equation}
\end{example}

Trivial principal $\cG$-bundles $P\to T$ are characterized by the existence of a section $T\to P$.
More precisely, such a section is equivalent to an isomorphism with a bundle of the form (\ref{swr}).
To see that, suppose that $P\to T$ is a principal $\cG$-bundle equipped with a section $\iota:T\to P$.
Let $f:P\to\cG_0$ and $a:\cG_1\times_{\cG_0} P\to P$ be the structure maps of the $\cG$-space $P$.
Then $T$ sits over $\cG_0$ via the composite $f\circ\iota:T\to P\to\cG_0$, and the composite
$$
\cG_1\times_{\cG_0} T\put(4.5,5.5){$\scriptstyle 1\times\iota$} -\!\!\!\longrightarrow \cG_1\times_{\cG_0} P\stackrel{a}{\longrightarrow} P
$$
is an isomorphism.
So our condition that $P\to T$ admits local sections can be reinterpreted as saying that $P$ is locally of the form (\ref{swr}).
Any principal $\cG$-bundle is therefore {\em locally} trivial, though not necessarily {\em globally} trivial.

Given a topological groupoid $\GG$ and a $\GG$-space $X$ (we often omit the map to $\GG_0$ from the notation), we obtain an {\em action groupoid}
\[
\cG\ltimes X:=\big(\cG_1\times_{\cG_0} X\rrarrow X\big),
\]
generalizing the earlier notion for group actions, whose structure maps are formally the same as before.
%Here, the source and target maps are the action and the projection, while the inverse, identity and multiplication maps are given by the same %formulae as before.
We also have the {\em gauge groupoid}
%of a principal $\GG$-bundle $P\to T$, introduced in Example \ref{Ggd} can be generalized to principal groupoid bundles.
%Given a principal $\cG$-bundle $P\to Y$, we let
\[
\gauge(P):=\big((P\times_{\cG_0}P)/\cG\rrarrow Y\big)
\]
of a principal $\GG$-bundle $P\to T$; that is, the relative pair groupoid $\pair_{\GG_0}(P)$ of the space $P$ over $\GG_0$ carries a (diagonal) action by $\GG$ (individually on object and arrow spaces), the quotient of which is the gauge groupoid $\gauge(P)=\pair_{\cG_0}(P)/\cG$.

\begin{example}\label{d:EG} \dontshow{d:EG}
The {\em translation groupoid} $\cE\cG$ associated to a topological groupoid $\cG$ is the action groupoid
$\cG\ltimes\cG_1$ for the regular action of $\cG$ on $\cG_1$.
Its objects are the arrows of $\cG$ and its arrows are the pairs of composable arrows of $\cG$; more generally, we have $(\cE\cG)_n=\cG_{n+1}$.
In the special case $\cG=\cB G$, we shall write $\cE G$ instead of $\cE(\cB G)$. Note the isomorphism $\cE G\cong\pair(G)$; it follows that
\(
\|\cE G\|\cong\|\pair(G)\|\cong EG,
\)
justifying our terminology.
\end{example}

Given a principal $\cG$-bundle $P\to T$, let $U\to T$ be a cover for which the projection $P\times_T U\to U$ admits a section $U\to P\times_T U$.
This induces a groupoid homomorphism $\pair_T U\to\cG$ via the maps $U\to P\to \cG_0$ and $U\times_T U\to P\times_T P\cong\cG_1\times_{\cG_0} P\to \cG_1$.
Conversely, given a groupoid homomorphism $\pair_T U\to\cG$ for some cover $U\to T$, the trivial principal $\cG$-bundle $\cG_1\times_{\cG_0}U\to U$
is such that its two pullbacks to $U\times_T U$ come with an isomorphism (induced from multiplication in $\cG$) which satisfies the obvious cocycle condition on $U\times_T U\times_T U$.
Hence $\cG_1\times_{\cG_0}U$ descends to a principal bundle on $T$.
These two constructions are inverse to one another, yielding the bijection
\dontshow{ugb}
\begin{equation}\label{ugb}
\left\{
\parbox{4cm}{Isomorphism classes of principal $\cG$-bundles on $X$}\hspace{.05cm}
\right\}\,\,\cong\,\, \underset{U\in\cov(T)}{\colim}\, \pi_0\Hom(\pair_T U,\cG).
\end{equation}
Here $\pi_0\Hom(\pair_T U,\cG)$ denotes the set of isomorphism classes of objects in $\Hom(\pair_T U,\cG)$, i.e.
the set of groupoid homomorphisms $\pair_T U\to\cG$ modulo natural transformation.

\begin{remark}
Given $f:\pair_T U\to \cG$ as above and two maps $q_1,q_2:V\to U$ over $T$, the resulting objects
$q_1^*f,q_2^*f\in\Hom(\cH_V,\cG)$ are isomorphic via $\tau:q_1^*f\Rightarrow q_2^*f$, $\tau(v):=f(q_1(v),1,q_2(v))$.
It follows that the above diagram $\pi_0\Hom(\pair_T (-),\cG):\cov(T)\to\{\text{Groupoids}\}$ factors through $\cov^0(X)$.
Thus we could have replaced $\cov(T)$ by $\cov^0(T)$ in the right hand side of (\ref{ugb}),
the latter being better behaved since $\cov^0(T)$ is filtered.
\end{remark}

\begin{definition}
A {\em principal $\cG$-bundle on a groupoid} $\cH$ consists of a principal $\cG$-bundle $P$ on $\cH_0$ together with an isomorphism
$\alpha:d_1^*P\to d_0^*P$ of principal $\cG$-bundles on $\cH_1$ which satisfies the cocycle condition on
$\cH_2$ in the sense that the diagram \dontshow{cddP}
\begin{equation}\label{cddP}
\xymatrix@R.1cm{ & d_2^* d_0^* P \ar@{=}[r] & d_0^* d_1^* P \ar[rd]^{d_0^*\alpha} & \\
d_2^* d_1^* P \ar[ur]^{d_2^*\alpha}\ar@{=}[dr] & & & d_0^* d_0^* P \\
& d_1^* d_1^* P \ar[r]^{d_1^*\alpha} & d_1^* d_0^* P\ar@{=}[ur] & }
\end{equation}
commutes.
An {\em isomorphism of principal $\GG$-bundles} $\eta:(Q,\beta)\to(P,\alpha)$ on $\cH$ is an isomorphism $\eta:Q\to P$ of principal $\cG$-bundles on $\cH_0$ such that $\alpha\circ d_1^*\eta = d_0^*\eta\circ\beta:d_1^*Q\to d_0^*P$.
\end{definition}

\begin{remark}
A principal $\cG$-bundle on $\cH$ is equivalent to a pair of principal $\cG$-bundles $P_0\to\cH_0$ and $P_1\to\cH_1$
together with isomorphisms $\alpha_i:P_1\to d_i^*P_0$ such that the composite $\alpha:=\alpha_0\alpha_1^{-1}$ satisfies the obvious cocycle condition.
Together these two spaces comprise the objects and arrows of a groupoid $\cP\to\cH$ over $\cH$.
\end{remark}

A principal $\cG$-bundle on $\cH$ gives rise to a right action of $\cH$ on $P$ via the map \dontshow{aktm}
\begin{equation}\label{aktm}
P\times_{\cH_0}\cH_1=d_1^*P\stackrel{\alpha}{\longrightarrow}d_0^*P \longrightarrow P,
\end{equation}
and that action commutes with the left action of $\cG$ on $P$.
The data of a space $P$ over $\GG_0$ and $\HH_0$, with commuting actions of $\cG$ and $\cH$, such that the map $P\to \cH_0$ is a principal $\cG$-bundle, is what's known as a {\em Hilsum-Skandalis morphism}
from $\cH$ to $\cG$.
These were introduced in \cite{HS87} as refinements of groupoid homomorphisms.
We shall depict a Hilsum-Skandalis morphism by the diagram \dontshow{hsp}
\begin{equation}\label{hsp}
\begin{matrix}
\cG_1\acts P\acted \cH_1\\
\downarrow\downarrow\hspace{.1cm}\swarrow\hspace{.3cm}\searrow\hspace{.1cm}\downarrow\downarrow\hspace{.1cm}\\
\phantom{.}\cG_0\hspace{1.17cm}\cH_0.
\end{matrix}
\end{equation}
Note that a Hilsum-Skandalis morphism is just an equivalent reformulation of a principal $\cG$-bundle on $\cH$;
indeed, the isomorphism $\alpha$ is recovered by the formula $\alpha(p,h)=(ph,h)$.

\begin{remark}\label{ryb}\dontshow{ryb}
If we put on (\ref{hsp}) the additional assumption that $P\to\cG_0$ is a principal $\cH$-bundle,
then we get the notion of a {\em Morita equivalence} between $\cH$ and $\cG$.
The existence of such a Morita equivalence is exactly what is needed for the stacks associated to $\cH$ and $\cG$ to be equivalent.
\end{remark}

The following proposition explains why it is natural to regard the data of (\ref{hsp}) as a specifying a morphism from $\cH$ to $\cG$.

\begin{proposition}\label{zuo}\dontshow{zuo}
Given topological groupoids $\cH$ and $\cG$, there is a natural equivalence of groupoids
\dontshow{u2}
\begin{equation}\label{u2}
\big\{
\text{Principal $\cG$-bundles on $\cH$}
\big\}\,\,\simeq\,\, \underset{U\in\cov^2(\cH_0)}{\hocolim}\,\, \Homone(\cH_U,\cG),
\end{equation}
where $\hocolim$ denotes the lax colimit in the category of groupoids (see Appendix \ref{A:2}).
\end{proposition}

\proof
Here we explain how to get a groupoid map from a principal bundle and vice versa; see Appendix \ref{A:2} for a complete proof (which also makes clear why $\cov^2$ is needed in place of $\cov$).
%is provided in Appendix \ref{A:2}.

Let $P\to\cH_0$ be a principal $\cG$-bundle and $\alpha:d_1^*P\to d_0^*P$ an isomorphism satisfying the cocycle condition (\ref{cddP}).
Pick a cover $f:U\to\cH_0$ over which $P$ trivializes and a section $U\to P$, thereby identifying $f^*P$ with the fibered product $\cG_1\times_{\cG_0} U$.
The cover $U\times_{\cH_0}\cH_1\times_{\cH_0} U$ of $\cH_1$ refines both $d_0^*U=U\times_{\cH_0}\cH_1$ and $d_1^*U=\cH_1\times_{\cH_0}U$,
and the trivialization of $f^*P$ induces trivializations of the restrictions of $d_0^*P$ and $d_1^*P$ to $U\times_{\cH_0}\cH_1\times_{\cH_0} U$.
The map $\alpha$ induces an isomorphism between these trivialized bundles,
which in turn corresponds to a map $U\times_{\cH_0}\cH_1\times_{\cH_0}U\to\cG_1$,
intertwining the two projections $U\times_{\cH_0}\cH_1\times_{\cH_0} U\to U$ with the source and target maps $\cG_1\to\cG_0$, respectively.
The cocycle condition amounts to asserting that the square
$$
\xymatrix{
(U\times_{\cH_0}\cH_1\times_{\cH_0}U)\times_U(U\times_{\cH_0}\cH_1\times_{\cH_0}U)\ar[r]\ar[d] & \cG_1\times_{\cG_0}\cG_1\ar[d]\\
 U\times_{\cH_0}\cH_1\times_{\cH_0}U                                              \ar[r]       & \cG_1}
$$
commutes, the vertical maps being those induced by composition in $\cG$ and $\cH$.
Hence we have constructed a functor $\cH_U\to\cG$.

Now suppose given a cover $U\to\cH_0$ and a groupoid homomorphism $\cH_U\to\cG$.
Precomposing with the inclusion $\pair_{\cH_0}(U)\to\cH_U$ gives a map $\pair_{\cH_0}(U)\to\cG$ classifying a principal $\cG$-bundle $P$ on $\cH_0$,
and the bundles $d_0^*P$ and $d_1^*P$ are classified by the two diagonal arrows in the diagram
\[
\xymatrix{
\pair_{\cH_1}(U\times_{\cH_0}\cH_1)\ar[r]^(.59){d_0}\ar[dr]&\pair_{\cH_0}(U)\ar[d]&
\pair_{\cH_1}(\cH_1\times_{\cH_0}U).\ar[l]_(.58){d_1}\ar[dl]\\&\cG
}
\]
An isomorphism $d_0^*P\cong d_1^*P$ corresponds to a natural transformation between the composites
\(
\pair_{\cH_1}(U\times_{\cH_0}\cH_1\times_{\cH_0}U)\!\to\!\pair_{\cH_1}(U\times_{\cH_0}\cH_1)\!\to\!\cG
\) and
\(
\pair_{\cH_1}(U\times_{\cH_0}\cH_1\times_{\cH_0}U)\!\to\!\pair_{\cH_1}(\cH_1\times_{\cH_0}U)\!\to\!\cG,
\)
which is to say a map
\[
\pair_{\cH_1}(U\times_{\cH_0}\cH_1\times_{\cH_0}U)\to\Map(\Delta^1_\gpd,\cG).
\]
Restricting the natural map $\cH_1\to\Map(\Delta^1_\gpd,\cH)$ to the cover $U\times_{\cH_0}\cH_1\times_{\cH_0}U\to\cH_1$,
\[
\pair_{\cH_1}(U\times_{\cH_0}\cH_1\times_{\cH_0}U)\to \Map(\Delta^1_\gpd,\cH)_{U\times_{\cH_0}\cH_1\times_{\cH_0}U}\cong \Map(\Delta^1_\gpd,\cH_U),
\]
and composing with the map $\Map(\Delta^1_\gpd,\cH_U)\to\Map(\Delta^1_\gpd,\cG)$ coming from $\cH_U\to\cG$ gives such a map.
The fact that the resulting isomorphism $\alpha:d_1^*P\to d_0^*P$ satisfies the cocycle condition is a consequence of the fact that
$\cH_U\to\cG$ preserves composition in the respective groupoids.
\qed

\begin{remark}\label{fb}\dontshow{fb}
The morphism groupoids $\Homone$ and $\Homtwo$ agree whenever the source has only trivial automorphism groups of objects, as is the case of the relative pair groupoid $\pair_T U$ of (\ref{ugb}).
%The source groupoid $\pair_T U$ of (\ref{ugb}) has only trivial automorphism groups of objects,
%so there is no difference between $\Homone$ and $\Homtwo$.
In (\ref{u2}), however, this is no longer the case,
so it is natural to ask which principal $\cG$-bundles on $\cH$ come from faithful functors $\cH_U\to\cG$, i.e. to elements of \dontshow{u4}
\begin{equation}\label{u4}
\underset{U\in\cov^2(\cH_0)}{\hocolim}\,\,\Homtwo(\cH_U,\cG).
\end{equation}
These are the bundles $P$ such that, for any object $x$ of $\cH_0$,
the topological group $\aut(x)$ acts {\em faithfully} on the fiber $P_x$ of $P$ over $x$ (and the image of $\aut(x)$ in $\aut(P_x)$ is closed).
\end{remark}

\begin{definition}
A principal $\cG$-bundle $E\to B$ is {\em universal}
if every principal $\cG$-bundle $P\to X$
over a paracompact base admits a $\cG$-bundle map
\[
\xymatrix@R=.5cm{
P\ar[r]\ar[d]&E\ar@<-.5ex>[d]\\
X\ar[r]&B,
}
\]
unique up to homotopy.
\end{definition}

The fat geometric realization $\|\cE\cG\|$ of the translation groupoid comes equipped with maps
\[
\begin{matrix}
\cG_1\acts \|\cE\cG\|\\
\downarrow\downarrow\hspace{.2cm}\swarrow\hspace{.3cm}\downarrow\hspace{.4cm}\\
\cG_0\hspace{.55cm}\|\cG\|\hspace{.15cm}
\end{matrix}
\]
that give it the structure of a principal $\cG$-bundle over $\|\cG\|$.
Indeed, a point in $\|\cE\cG\|$ can be represented as a sequence of composable morphisms
\[
(x\stackrel{\alpha_0}{\leftarrow}r_0\cdot x_0\stackrel{\alpha_1}{\leftarrow}r_1\cdot x_1\stackrel{\alpha_2}{\leftarrow}\cdots\stackrel{\alpha_n}{\leftarrow}r_n\cdot x_n),\qquad  r_i\in[0,1],\,\,\,\,\,\sum_{i=0}^nr_i=1,
\]
where all but the first object are given a weight.
The projections to $\|\cG\|$ and $\cG_0$ are then given by
\[
\begin{split}
&(x\stackrel{\alpha_0}{\leftarrow}r_0\cdot x_0\stackrel{\alpha_1}{\leftarrow}\cdots\stackrel{\alpha_n}{\leftarrow}r_n\cdot x_n)
\mapsto (r_0\cdot x_0\stackrel{\alpha_1}{\leftarrow}\cdots\stackrel{\alpha_n}{\leftarrow}r_n\cdot x_n),\\%\quad\text{and}\\
&(x\stackrel{\alpha_0}{\leftarrow}r_0\cdot x_0\stackrel{\alpha_1}{\leftarrow}\cdots\stackrel{\alpha_n}{\leftarrow}r_n\cdot x_n)\mapsto x,
\end{split}
\]
respectively, and $g:y\leftarrow x$ acts by
\[
g\cdot(x\stackrel{\alpha_0}{\leftarrow}r_0\cdot x_0\stackrel{\alpha_1}{\leftarrow}\cdots\stackrel{\alpha_n}{\leftarrow}r_n\cdot x_n)
=(y\stackrel{g\alpha_0}{\leftarrow}r_0\cdot x_0\stackrel{\alpha_1}{\leftarrow}\cdots\stackrel{\alpha_n}{\leftarrow}r_n\cdot x_n).
\]
It is then easy to verify the all the conditions for being a principal $\cG$-bundle are satisfied.
Moreover, there is a canonical section $s_{\|\cG\|}:U_{\|\cG\|}\to\|\cE\cG\|$ from the preferred cover $U_{\|\cG\|}$ of $\|\cG\|$ given by
\[
(r_0\cdot x_0\stackrel{\alpha_1}{\leftarrow}\cdots\stackrel{\alpha_i}{\leftarrow}\underline{r_i\cdot x_i}\stackrel{\alpha_{i+1}}{\leftarrow}\cdots\stackrel{\alpha_n}{\leftarrow}r_n\cdot x_n)
\mapsto
(x_i
\put(2,6){$\scriptstyle(\alpha_1\alpha_2\cdots\alpha_i)^{-1}$}
\longleftarrow\!\!\!-\!\!\!-\!\!\!-\!\!\!-\!\!\!-\:
r_0\cdot x_0\stackrel{\alpha_1}{\leftarrow}\cdots\stackrel{\alpha_n}{\leftarrow}r_n\cdot x_n).
\]

\begin{proposition}\label{itc}\dontshow{itc}
The principal $\cG$-bundle $\|\cE\cG\|\to\|\cG\|$ is a universal $\cG$-bundle.
\end{proposition}

\proof
Let $P\to T$ be a principal $\cG$-bundle over a paracompact space $T$ and let $U\to T$ be a cover equipped with a section $s:U\to P$.
Since $T$ is paracompact, we may refine $U$ to a cover whose fibers are finite and pick a partition of unity $(U,\varphi,<)$.

Given a point $t$ in $T$ with preimages $t_0,\ldots, t_n$ in $U$, write $p_i:=s(t_i)$ and $x_i$ for their images in $\cG_0$.
Since $P$ is principal, there exist unique morphisms $g_i:x_{i-1}\leftarrow x_i$ such that $g_i\cdot p_i=p_{i-1}$.
Given a point $p\in P$ over $t\in T$, with image $x\in\cG_0$, we let $g$ be the unique morphisms such that $g\cdot p_0=p$.
We can then define our desired bundle map $f:(P\to X)\to(\|\cE\cG\|\to\|\cG\|)$ by \dontshow{hrp}
\begin{equation}\label{hrp}
\begin{split}
f(t)&=\big(\varphi(t_0)\cdot x_0\stackrel{g_1}{\leftarrow}\varphi(t_1)\cdot x_1\leftarrow\cdots\stackrel{g_n}{\leftarrow}\varphi(t_n)\cdot x_n\big)\\
f(p)&=\big(x\stackrel{g}{\leftarrow}\varphi(t_0)\cdot x_0\stackrel{g_1}{\leftarrow}\varphi(t_1)\cdot x_1\leftarrow\cdots\stackrel{g_n}{\leftarrow}\varphi(t_n)\cdot x_n\big).
\end{split}
\end{equation}

Recall that the bundle $\|\cE\cG\|\to\|\cG\|$ comes equipped with a preferred partition of unity $(U_{\|\cG\|},\varphi_{\|\cG\|},<)$ and a section $s_{\|\cG\|}$.
Given a $\cG$-bundle map $f:(P\to T)\to(\|\cE\cG\|\to\|\cG\|)$, one can pullback the above data to a get a corresponding partition of unity and section over $T$.
% which pulls backs to a partition of unity and section for bundle $P\to T$ along a bundle map $f:(P\to T)\to(\|\cE\cG\|\to\|\cG\|)$.
It is not hard to see that the procedure (\ref{hrp}), applied to $\|\cE\cG\|\to\|\cG\|$ itself using the partition of unity $(U_{\|\cG\|},\varphi_{\|\cG\|},<)$ and section $s_{\|\cG\|}$,
recovers the identity map of the bundle $\|\cE\cG\|\to\|\cG\|$.
Thus, for any principal $\cG$-bundle $P\to T$, the formulae (\ref{hrp}) define a bijection between pairs $((U,\varphi,<),s)$ and bundle maps $f:(P\to T)\to(\|\cE\cG\|\to\|\cG\|)$.

We now show that any two such bundle maps $f$ and $f'$ are homotopic.
Let $(U,\varphi,<)$, $(U',\varphi',<')$ and $s$, $s'$ denote the pulled back partitions of unity and sections, respectively.
Then the cover $U'':=(U\times [0,1))\sqcup (U'\times(0,1])\to T\times [0,1]$ admits a section
$s''=(s\times \Id)\sqcup(s'\times \Id):U''\to P\times [0,1]$ and a partition of unity
given by
\[
\varphi''(u,\epsilon)=\begin{cases}
(1-\epsilon)\varphi(u)&\text{if}\quad u\in U\\
\epsilon\,\varphi'(u)&\text{if}\quad u\in U'\\
\end{cases}
\]
and
\[
(u,\epsilon)<''(v,\epsilon)\,\stackrel{\mathrm {def}}{\Longleftrightarrow}\,
\begin{cases}
u<v&\text{if}\,\, u,v\in U\\
u<'v&\text{if}\,\, u,v\in U'\\
u\in U,\, v\in U'&\text{otherwise.}
\end{cases}
\]
The map defined by the equations (\ref{hrp}) using $(U'',\varphi'',<'')$ is then a homotopy between $f$ and $f'$.
\qed

\subsection{Fibrant topological groupoids \rm\dontshow{se:Fr}}\label{se:Fr}

Recall that two topological groupoids $\cG$ and $\cH$ are said to be {\em categorically equivalent} if there are continuous functors $\sigma:\cH\to\cG$, $\tau:\cG\to\cH$ and continuous natural transformations $\sigma\circ\tau\Rightarrow 1_\cG$, $\tau\circ\sigma\Rightarrow 1_\cH$.
If instead of natural transformations we have
homotopies $\cH\times[0,1]\to \cG$, and $\cG\times[0,1]\to\cH$ from $g\circ f$ to $1_\cH$ and from $f\circ g$ to $1_\cG$,
%homotopies $\sigma\circ\tau\Rightarrow 1_\cG$, $\tau\circ\sigma\Rightarrow 1_\cH$,
then we say that $\cG$ and $\cH$ are {\em homotopy equivalent}.
Neither of these notions implies the other.
Also of interest to us is the notion of {\em Morita equivalence}, briefly alluded to in Remark \ref{ryb}; it is implied by categorical equivalence and unrelated to homotopy equivalence.

In this section, we introduce the notion of a fibrant topological groupoid and construct a fibrant replacement functor $\cG\mapsto \fib(\cG)$.
The groupoid $\fib(\cG)$ is Morita equivalent to $\cG$ but better suited for receiving maps from other groupoids.
Recall that a principal $\cG$-bundle is called {\em universal} if any other $\cG$-bundle over a paracompact base
admits a map to it, unique up to homotopy.

\begin{definition}\label{df:fib}\dontshow{df:fib}
A topological groupoid $\cG$ is {\em fibrant} if $\cG_1\to\cG_0$ is a universal $\cG$-bundle.
\end{definition}

\begin{definition}
The {\em fibrant replacement} $\fib(\cG)$ of a topological groupoid $\cG$ is the gauge groupoid of the universal principal
$\cG$-bundle $\|\cE\cG\|\to\|\cG\|$; that is,
\[
\fib(\cG):=\gauge\big(\|\cE\cG\|\big)=\big(\pair_{\cG_0}\|\cE\cG\|\big)/\cG=\big((\|\cE\cG\|\times_{\cG_0}\|\cE\cG\|)/\cG \rrarrow\|\cG\|\big).
\]
\end{definition}
A priori, it is not clear that $\fib(\cG)$ is a fibrant groupoid; this is the content of Proposition \ref{half}.
We begin with a more explicit description of the objects and arrows of $\fib(\cG)$.
The object space is just the realization of $\cG$ and its points have been described in (\ref{obF}).
The arrow space is the space of composable sequences of morphisms of weighted objects \dontshow{rxa}
\begin{equation}\label{rxa}
(r_0\cdot x_0\stackrel{\alpha_1}{\leftarrow}\cdots\stackrel{\alpha_n}{\leftarrow}r_n\cdot x_n\stackrel{\beta}{\leftarrow}
s_0\cdot y_0\stackrel{\gamma_1}{\leftarrow}\cdots\stackrel{\gamma_m}{\leftarrow}s_m\cdot y_m)
\end{equation}
where
$\sum r_i=\sum s_i=1$,
modulo the same equivalence relation as the one used in defining the fat geometric realization.
We shall call $\beta$ the {\em middle element}\hspace{.05cm}\footnote
{
Strictly speaking, the middle element only makes sense given a representative (\ref{rxa}) of a morphism.
} of the morphism (\ref{rxa}).
Composition is given by
\[
\begin{split}
&(r_0\cdot x_0\stackrel{\alpha_1}{\leftarrow}\cdots\stackrel{\alpha_n}{\leftarrow}r_n\cdot x_n\stackrel{\beta}{\leftarrow}
s_0\cdot y_0\stackrel{\gamma_1}{\leftarrow}\cdots\stackrel{\gamma_m}{\leftarrow}s_m\cdot y_m)
(s_0\cdot y_0\stackrel{\gamma_1}{\leftarrow}\cdots\stackrel{\gamma_m}{\leftarrow}s_m\cdot y_m\stackrel{\delta}{\leftarrow}
t_0\cdot z_0\stackrel{\varepsilon_1}{\leftarrow}\cdots\stackrel{\varepsilon_l}{\leftarrow}t_l\cdot z_l)\\
&\hspace{2cm}:=(r_0\cdot x_0\stackrel{\alpha_1}{\leftarrow}\cdots\stackrel{\alpha_n}{\leftarrow}r_n\cdot x_n
\put(2,6){$\scriptstyle
\beta\gamma_1\cdots\gamma_m\delta
$}
\longleftarrow\!\!\!-\!\!\!-\!\!\!-\:
t_0\cdot z_0\stackrel{\varepsilon_1}{\leftarrow}\cdots\stackrel{\varepsilon_l}{\leftarrow}t_l\cdot z_l),
\end{split}
\]
inversion by
\[
\begin{split}
&(r_0\cdot x_0\stackrel{\alpha_1}{\leftarrow}\cdots\stackrel{\alpha_n}{\leftarrow}r_n\cdot x_n\stackrel{\beta}{\leftarrow}
s_0\cdot y_0\stackrel{\gamma_1}{\leftarrow}\cdots\stackrel{\gamma_m}{\leftarrow}s_m\cdot y_m)^{-1}\\
&\hspace{2cm}:=(s_0\cdot y_0\stackrel{\gamma_1}{\leftarrow}\cdots\stackrel{\gamma_m}{\leftarrow}s_m\cdot y_m
\put(2,6){$\scriptstyle
(\alpha_1\cdots\alpha_n\beta\gamma_1\cdots\gamma_m)^{-1}
$}
\longleftarrow\!\!\!-\!\!\!-\!\!\!-\!\!\!-\!\!\!-\!\!\!-\!\!\!-\!\!\!\!-\!\!\!-\!\!\!-\:
r_0\cdot x_0\stackrel{\alpha_1}{\leftarrow}\cdots\stackrel{\alpha_n}{\leftarrow}r_n\cdot x_n),
\end{split}
\]
and the identities by
\[
1_{(r_0\cdot x_0\stackrel{\alpha_1}{\leftarrow}\cdots\stackrel{\alpha_n}{\leftarrow}r_n\cdot x_n)}
:=(r_0\cdot x_0\stackrel{\alpha_1}{\leftarrow}\cdots\stackrel{\alpha_n}{\leftarrow}r_n\cdot x_n
\put(2,6){$\scriptstyle
(\alpha_1\cdots\alpha_n)^{-1}
$}
\longleftarrow\!\!\!-\!\!\!-\!\!\!\!-\!\!\!\!-\:
r_0\cdot x_0\stackrel{\alpha_1}{\leftarrow}\cdots\stackrel{\alpha_n}{\leftarrow}r_n\cdot x_n).
\]
The inclusion of principal $\cG$-bundles $(\cG_1\to\cG_0)\hookrightarrow(\|\cE\cG\|\to\|\cG\|)$ induces a natural transformation $\eta:1\to \fib$ by taking gauge groupoids:
\[
\eta:\cG=\gauge(\cG_1)\,\rightarrow\,\gauge\big(\|\cE\cG\|\big)=\fib(\cG).
\]
More explicitly, it sends an object $x$ of $\cG$ to the object $(1\cdot x)$ of $\fib(\cG)$ and a morphism $x_0\stackrel{\alpha}{\leftarrow}x_1$ of $\cG$ to the morphism
$(1\cdot x_0\stackrel{\alpha}{\leftarrow}1\cdot x_1)$ of $\fib\cG$.

\begin{example}\label{dugg}\dontshow{dugg}
Given a topological group $G$, let $\cU G:=\fib \cB G$ be the gauge groupoid of the universal principal $G$-bundle $EG=\|\cE G\|\to\|\cB G\|$.
It can be written explicitly as
\[
\cU G\cong (\pair EG)/G\cong \big((EG\times EG)/G\rrarrow EG/G\big).
\]
We shall call it the {\em universal gauge groupoid} of the group $G$.
\end{example}

We now define a natural transformation $\mu:\fib^2\to\fib$.
Observe that an object $x$ of $\fib^2(\cG)$ is represented by a sequence of the form
\[
x=\big(
s_0\cdot(t_{00}\cdot x_{00}\stackrel{\alpha_{01}}{\leftarrow}\cdots\stackrel{\alpha_{0n_0}}{\leftarrow}t_{0n_0}\cdot x_{0n_0})
\stackrel{\alpha_1}{\leftarrow}\cdots\stackrel{\alpha_m}{\leftarrow}
s_m\cdot(t_{m0}\cdot x_{m0}\stackrel{\alpha_{m1}}{\leftarrow}\cdots\stackrel{\alpha_{mn_m}}{\leftarrow}t_{mn_m}\cdot x_{mn_m})
\big)
\]
in which the $\alpha_{ij}$ are morphisms in $\cG$, the $\alpha_i$ are morphisms in $\fib(\cG)$, and the $s_i$ and $t_{ij}$ are
such that
$$
1 = \sum_i s_i = \sum_i\Big(s_i\sum_j t_{ij}\Big) = \sum_{i,j} s_i t_{ij}.
$$
Identifying the $\alpha_i$ with their middle elements $x_{i-1, n_{i-1}}\!\leftarrow x_{i 0}$, we obtain an object
\[
\mu(x):=(
s_0t_{00}\cdot x_{00}\stackrel{\alpha_{01}}{\leftarrow}\cdots\stackrel{\alpha_{0n_0}}{\leftarrow}s_0t_{0n_0}\cdot x_{0n_0}
\stackrel{\alpha_1}{\leftarrow}\cdots\stackrel{\alpha_m}{\leftarrow}
s_mt_{m0}\cdot x_{m0}\stackrel{\alpha_{m1}}{\leftarrow}\cdots\stackrel{\alpha_{mn_m}}{\leftarrow}s_mt_{mn_m}\cdot x_{mn_m})
\]
in $\fib(\cG)$.
The transformation $\mu$ on morphisms is defined similarly.

\begin{remark}
The transformation $\mu$ would fail to be well defined if we had used $|\cG|$ instead of $\|\cG\|$ to define $\fib(\cG)$.
Indeed, the two objects
$$
\big(
\half\cdot(\half\cdot x\stackrel{\alpha}{\leftarrow}\half\cdot y)\stackrel{\alpha^{-1}}{\leftarrow}\half\cdot(\half\cdot x\stackrel{\alpha}{\leftarrow}\half\cdot y)
\big)
\qquad\text{and}\qquad
\big(
1\cdot(\half\cdot x\stackrel{\alpha}{\leftarrow}\half\cdot y)
\big)
$$
would be identified, but not their images.
\end{remark}

\begin{remark}\label{monad}\dontshow{monad}
The triple $(\fib,\mu,\eta)$ defines a monad on the category of topological groupoids.
\end{remark}

\begin{proposition}\label{him}\dontshow{him}
For each topological groupoid $\cG$, the natural map $\mu(\cG):\fib^2(\cG)\to\fib(\cG)$ is a homotopy equivalence.
\end{proposition}

\proof
The natural transformation $\eta\circ\fib:\fib\to\fib^2$ gives a section of $\mu:\fib^2\to\fib$.
So it suffices to show that the natural transformation
$$
\fib^2\stackrel{\mu}{\longrightarrow}\fib\stackrel{\eta\circ\fib}{\longrightarrow}\fib^2
$$
is homotopic to the identity.
We define such a homotopy by sending the point
$$
\big(
s_0\cdot(t_{00}\cdot x_{00}\stackrel{\alpha_{01}}{\leftarrow}\cdots\stackrel{\alpha_{0n_0}}{\leftarrow}t_{0n_0}\cdot x_{0n_0})
\stackrel{\alpha_1}{\leftarrow}\cdots\stackrel{\alpha_m}{\leftarrow}
s_m\cdot(t_{m0}\cdot x_{m0}\stackrel{\alpha_{m1}}{\leftarrow}\cdots\stackrel{\alpha_{mn_m}}{\leftarrow}t_{mn_m}\cdot x_{mn_m});u
\big)
$$
in $\fib^2(\cG)\times[0,1]$ to the point
\begin{eqnarray*}
us_0\cdot(t_{00}\cdot x_{00}\stackrel{\alpha_{01}}{\leftarrow}\cdots\stackrel{\alpha_{0n_0}}{\leftarrow}t_{0n_0}\cdot x_{0n_0})
\stackrel{\alpha_1}{\leftarrow}\cdots\stackrel{\alpha_m}{\leftarrow}
us_m\cdot(t_{m0}\cdot x_{m0}\stackrel{\alpha_{m1}}{\leftarrow}\cdots\stackrel{\alpha_{mn_m}}{\leftarrow}t_{mn_m}\cdot x_{mn_m})
\stackrel{\beta^{-1}}{\leftarrow}
\\
\phantom{\stackrel{\beta}{\leftarrow}}
\bar{u}\cdot(s_0t_{00}\cdot x_{00}\stackrel{\alpha_{01}}{\leftarrow}\cdots\stackrel{\alpha_{0n_0}}{\leftarrow}s_0t_{0n_0}\cdot x_{0n_0}
\stackrel{\alpha_1}{\leftarrow}\cdots\stackrel{\alpha_m}{\leftarrow}
s_mt_{m0}\cdot x_{m0}\stackrel{\alpha_{m1}}{\leftarrow}\cdots\stackrel{\alpha_{mn_m}}{\leftarrow}s_mt_{mn_m}\cdot x_{mn_m}),
\end{eqnarray*}
where $\beta$ is the composite
$x_{00}\stackrel{\alpha_{01}}{\leftarrow}\cdots
\stackrel{\alpha_{mn_m}}{\leftarrow}x_{mn_m}$,
and $\bar{u}=1-u$.
The homotopy on morphisms is defined analogously.
\qed
\vspace{.3cm}

\begin{corollary}\label{fibsq}\dontshow{fibsq}
The natural transformation $\mu:\fib^2\to\fib$ is a homotopy inverse for each of the natural transformations
$\eta\circ\fib,\fib\circ\eta:\fib\to\fib^2$.\qed
\end{corollary}

\begin{proposition}\label{fibfactor}\dontshow{fibfactor}
The map
$$
\psi:\map\big(\fib\cH,\fib\cG\big)\to\map\big(\cH,\fib\cG\big)
$$
obtained by restricting along $\eta:\cH\to\fib\cH$ is a homotopy equivalence.
\end{proposition}
\proof
Let $\phi$ be the composite
$$
\phi:\map(\cH,\fib\cG)\stackrel{\fib}{\longrightarrow}\map(\fib\cH,\fib^2\cG)\stackrel{\mu} {\longrightarrow}\map(\fib\cH,\fib\cG),
$$
and note that
$
\psi\circ\phi:\map(\cH,\fib\cG))\to\map(\fib\cH,\fib\cG)\to\map(\cH,\fib\cG)
$
is the identity.
So it suffices to find a homotopy from the identity of $\map(\fib\cH,\fib\cG)$ to the composite
$$
\phi\circ\psi:\map(\fib\cH,\fib\cG)\to\map(\cH,\fib\cG)\to\map(\fib\cH,\fib\cG).
$$

To any arrow $g=(r_0\cdot x_0\stackrel{\alpha_1}{\leftarrow}\cdots r_n\cdot x_n\stackrel{\beta}{\leftarrow}s_0\cdot y_0\stackrel{\gamma_1}
{\leftarrow}\cdots s_m\cdot y_m)\in\fib(\cG)_1$ corresponds a path $\delta_g:[0,1]\to\fib(\cG)_0$ between $t(g)$ and $s(g)$, given by \dontshow{aqy}
\begin{equation}\label{aqy}
\delta_g(u):=((1-u)r_0\cdot x_0\stackrel{\alpha_1}{\leftarrow}\cdots\stackrel{\alpha_n}{\leftarrow}(1-u)r_n\cdot x_n\stackrel{\beta}
{\leftarrow}us_0\cdot y_0\stackrel{\gamma_1}{\leftarrow}\cdots\stackrel{\gamma_m}{\leftarrow}us_m\cdot y_m).
\end{equation}
Recall that $\beta$ is called the middle element of $g$.
Given a morphism $f:\fib(\cH)\to\fib(\cG)$ and an object
$y=(s_0\cdot y_0\stackrel{\gamma_1}{\leftarrow}\cdots\stackrel{\gamma_m}{\leftarrow}s_m\cdot y_m)\in\fib(\cH)$,
we need to produce a path between $f(y)$ and $\phi(\psi(f))(y)$.
Let $g\in\fib(\cG)_1$ be the morphism given by $s(g)=y$, $t(g)=\phi(\psi(f))(y)$, and whose middle element is the middle element of
\[
f(s_0\cdot y_0\stackrel{\gamma_1}{\leftarrow}\cdots\stackrel{\gamma_m}{\leftarrow}s_m\cdot y_m
\put(2,6){$\scriptstyle
(\gamma_1\gamma_2\ldots\gamma_m)^{-1}
$}
\longleftarrow\!\!\!-\!\!\!-\!\!\!-\!\!\!-\!\!\!-\!\!\!-\:
1\cdot y_0)
\]
Then $\delta_g$ is the desired path.

The construction on arrows is similar.
It uses the fact that to any commutative square in $\fib(\cG)$, one can assign a path $[0,1]\to\fib(\cG)_1$ between the two vertical arrows of the square.
\qed
\vspace{.3cm}

\begin{lemma}\label{Kri}\dontshow{Kri}
Given two topological groupoids $\cG$, $\cH$, then $\fib(\cG\times \cH)$ is a retract of $\fib(\cH)\times \fib(\cG)$ via maps
$\iota:\fib(\cG\times \cH)\put(3,6){$\curvearrowleft$}\hookrightarrow\fib(\cH)\times \fib(\cG):r$.
Moreover, given three groupoids $\cG$, $\cH$, $\cK$ we have equalities \dontshow{irk}
\begin{equation}\label{irk}
\begin{split}
\iota\circ(\iota\times 1)&= \iota\circ(1\times \iota):\fib(\cH\times \cG\times \cK)\to\fib(\cH)\times \fib(\cG)\times \fib(\cK),\\
r\circ(r\times 1)&= r\circ(1\times r):\fib(\cH)\times \fib(\cG)\times \fib(\cK)\to\fib(\cH\times \cG\times \cK).
\end{split}
\end{equation}
\end{lemma}

\proof
At the level of objects, this is the statement (\ref{iri}) about fat geometric realizations, which was left as an exercise to the reader.
We now define the maps $\iota$ and $r$ on arrows, leaving the rest of the proof to the reader.
Given an arrow $g\in \fib(\cG\times \cH)_1$ with middle element $(\beta_1,\beta_2)$,
let $\iota(g)\in \fib(\cH)\times \fib(\cG)$ be the pair $(g_1,g_2)\in\fib(\cH)_1\times \fib(\cG)_1$, with sources and targets given by
$\iota(s(g))=(s(g_1),s(g_2))$, $\iota(t(g))=(t(g_1),t(g_2))$ and middle elements $\beta_1$, $\beta_2$.
In the other direction, given arrows $g_1\in \fib(\cH)$, $g_2\in \fib(\cG)$ with middle elements $\beta_1$, $\beta_2$,
the image $r(g_1,g_2)$ is the arrow $r(s(g_1),s(g_2))\to r(t(g_1),t(g_2))$ with middle element $(\beta_1,\beta_2)$.
\qed
\vspace{.3cm}

We now justify our use of the terminology ``fibrant replacement'' by showing that $\fib(\cG)$ is indeed a fibrant groupoid.

\begin{proposition}\label{half}\dontshow{half}
Let $\cG$ be any topological groupoid.
Then $\fib(\cG)$ is fibrant.
\end{proposition}

\proof
The categories of $\cG$-principal bundles and $\fib(\cG)$-principal bundles are equivalent via functors
\begin{eqnarray*}
\big(P\to X\big)\,&\longmapsto&\, \Big((P\times_{\cG_0}\|\cE\cG\|)/\cG_1\to X\Big),\\
\Big((Q\times_{\fib(\cG)_0}\|\cE\cG\|)/\fib(\cG)_1\to X\Big)
\,\,&\raisebox{1.7mm}{\rotatebox{180}{$\longmapsto$}}&\,\,
\big(Q\to X\big).
\end{eqnarray*}
The image of the universal principal $\cG$-bundle $\|\cE\cG\|\to\|\cG\|$ is the principal
$\fib(\cG)$-bundle
\[
(\|\cE\cG\|\times_{\cG_0}\|\cE\cG\|)/\cG_1=\fib(\cG)_1\to\|\cG\|=\fib(\cG)_0.
\]
The result then follows from Proposition \ref{itc}.
\qed

\begin{lemma}\label{stur}\dontshow{stur}
Let $\cG$ be a fibrant groupoid and let $\iota:\cH\put(3,6){$\curvearrowleft$}\hookrightarrow\cG:r$ be a retract of $\cG$.
Then $\cH$ is fibrant.
\end{lemma}

\proof
According to the dictionary (\ref{ugb}), a principal $\cH$-bundle $P\to T$ corresponds to an element $\alpha\in\hocolim_U\Hom(\pair_X U,\cH)$.
Constructing a bundle map $(P\to T)\to(\cH_1\to\cH_0)$ is then equivalent to factoring $\alpha$ through the inclusion $\cH_0\hookrightarrow\cH$.
Since $\cG$ is fibrant, the composite $\iota\circ\alpha$ factors through $\cG_0\hookrightarrow\cG$.
Composing with $r_0$ then produces the desired lift $T\to\cH_0$.

Given bundle maps $f,f':(P\to T)\to(\cH_1\to\cH_0)$, we may compose $f_0,f_0':T\to\cH_0$ with the inclusion $\cH_0\hookrightarrow\cH$ in order to produce objects
$\alpha, \beta \in \Hom(T,\cH)$.
Since both $\alpha$ and $\beta$ classify $P$ we also get a morphism $\tau:\alpha\to\beta$,
given explicitly by the formula $\tau(t)=f'(p)^{-1}f(p)$ for every $p\in P$ with image $t\in T$.
This morphism induces an element
$\theta\in\Hom(\pair_{T\times[0,1]}(T\times[0,1)\sqcup T\times (0,1]),\cH)$
whose restrictions to the ends are $\alpha$ and $\beta$.
Finding a homotopy between $f$ and $f'$ is then equivalent to finding
an extension of $\theta$ to $\pair_{T\times[0,1]}(T\times[0,1)\sqcup T\times (0,1]\sqcup T\times [0,1])$.
As before, we first find an extension of $\iota\circ\theta$ and then compose it with $r$.
\qed

\begin{proposition}\label{ccFR}\dontshow{ccFR}
A paracompact topological groupoid $\cG$ is fibrant if and only if $\cG$ is a retract of $\fib(\cG)$.
\end{proposition}

\proof
If $\cG$ is a retract of $\fib(\cG)$ then $\cG$ is fibrant by Proposition \ref{half} and Lemma \ref{stur}.
For the other direction, let $\cG$ be a fibrant groupoid.
Since $\cG_1\to \cG_0$ is universal, it comes with a map of $\cG$-principal bundles $r:(\|\cE\cG\|\to\|\cG\|)\to(\cG_1\to \cG_0)$ which is homotopy left inverse to the inclusion $i:(\cG_1\to \cG_0)\hookrightarrow(\|\cE\cG\|\to\|\cG\|)$, and since $i$ is a neighborhood deformation retract in the category of $\cG$-principal bundles, we can modify $r$ so that it becomes a strict left inverse of $i$.
It follows that $\cG=\gauge(\cG_1)$ is a retract of $\fib(\cG)=\gauge(\|\cE\cG\|)$.
\qed
\vspace{.3cm}

In the unenriched setting, a fully faithful and essentially surjective morphism of groupoids is automatically a categorical equivalence.
In the topological setting, however, this is not enough
(see \cite[section 2.4]{Moe02} for precise definitions of these notions in the enriched context); we must also demand that the map in question admits a one-sided inverse.

\begin{lemma}\label{hece}\dontshow{hece}
If $\cG$ is paracompact and fibrant then $\eta:\cG\to\fib(\cG)$ is both a homotopy and a categorical equivalence.
\end{lemma}

\proof
By Proposition \ref{ccFR}, we may write $\cG$ as retract of $\fib(\cG)$.
It follows that $\eta:\cG\to\fib(\cG)$ is a retract of \dontshow{sjh}
\begin{equation}\label{sjh}
\eta\circ\fib:\fib(\cG)\to\fib^2(\cG).
\end{equation}
As both homotopy and categorical equivalences are preserved by retracts, it's enough to show that (\ref{sjh}) is a homotopy and categorical equivalence.
The former is the content of Corollary \ref{fibsq}, while the latter follows since $(\ref{sjh})$ is fully faithful, essentially surjective, and has a left inverse $\mu:\fib^2(\cG)\to\fib(\cG)$.
\qed
\vspace{.3cm}

The following result ensures that the fibrant replacement functor behaves well with respect to the operation of gluing cells.

\begin{lemma}\label{sip}\dontshow{sip}
Fibrant replacement preserves pushouts along closed saturated inclusions;
that is, for any closed saturated inclusion $\cH_Z\to\cH$ and any map $\cH_Z\to\cG$, the natural map \dontshow{gdy}
\begin{equation}\label{gdy}
\fib(\cH)\underset{\fib(\cH_Z)}{\sqcup}\fib(\cG)
\to
\fib\big(\cH\!\underset{\,\,\cH_Z}{\sqcup}\!\cG\big)
\end{equation}
is an isomorphism.
\end{lemma}

\proof
We first show that \dontshow{hdo}
\begin{equation}\label{hdo}
\cH_n\!\!\underset{\,\,\,(\cH_Z\!)_n}{\sqcup}\!\cG_n
\cong
\big(\cH\!\underset{\,\,\cH_Z}{\sqcup}\!\cG\big)_n.
\end{equation}
Clearly, we have
$\cH_0\sqcup_Z\cG_0
\cong
(\cH\sqcup_{\,\cH_Z}\!\cG)_0$, so the first non trivial case is \dontshow{ihs}
\begin{equation}\label{ihs}
\cH_1\!\!\underset{\,\,\,(\cH_Z\!)_1}{\sqcup}\!\cG_1
\cong
\big(\cH\!\underset{\,\,\cH_Z}{\sqcup}\!\cG\big)_1.
\end{equation}
The pushout $\cH\sqcup_{\,\cH_Z}\!\cG$ is the quotient of the free groupoid on $\cH_1\sqcup_{\,(\cH_Z\!)_1}\!\cG_1$
by relations specifying the multiplication on $\cH_2\sqcup_{\,(\cH_Z\!)_2}\!\cG_2$.
By Lemma \ref{pushpull}, the latter is isomorphic to
\[
%\cH_2\!\!\underset{\,\,\,(\cH_Z\!)_2}{\sqcup}\!\cG_2\simeq
\big(\cH_1\!\!\underset{\,\,\,(\cH_Z\!)_1}{\sqcup}\!\cG_1\big)\underset{\,\,\cH_0\sqcup_Z\cG_0}{\times}
\big(\cH_1\!\!\underset{\,\,\,(\cH_Z\!)_1}{\sqcup}\!\cG_1\big),
\]
so the relations fully determine the multiplication on $\cH_1\sqcup_{\,(\cH_Z\!)_1}\!\cG_1$ and (\ref{ihs}) follows.
The rest of (\ref{hdo}) is a straightforward application of Lemma \ref{pushpull}.

Since geometric realization commutes with pushouts of simplicial spaces, it follows that
\[
\|\cH\|\underset{\,\|\cH_Z\!\|}{\sqcup}\|\cG\|
\cong
\big\|\cH\underset{\,\,\cH_Z}{\sqcup}\cG\big\|
\qquad\text{and}\qquad
\|\cE\cH\|\underset{\,\|\cE\cH_Z\!\|}{\sqcup}\|\cE\cG\|
\cong
\big\|\cE\big(\cH\underset{\,\,\cH_Z}{\sqcup}\cG\big)\big\|.
\]
At the level of objects, we have already seen that (\ref{gdy}) is an isomorphism.
At the level of arrows, we can now compute
\[
\begin{split}
\hspace{2.2cm}\fib\big(\cH\!\underset{\,\,\cH_Z}{\sqcup}\!\cG\big)_1
&\cong
\Big[\big\|\cE\big(\cH\underset{\,\,\cH_Z}{\sqcup}\cG\big)\big\|
\underset{\,\,\cH_0\sqcup_Z\cG_0}{\times}
\big\|\cE\big(\cH\underset{\,\,\cH_Z}{\sqcup}\cG\big)\big\|\Big]\Big/\sim\\
&\cong
\Big[\big(\|\cE\cH\|\underset{\,\|\cE\cH_Z\!\|}{\sqcup}\|\cE\cG\|\big)
\underset{\,\,\cH_0\sqcup_Z\cG_0}{\times}
\big(\|\cE\cH\|\underset{\,\|\cE\cH_Z\!\|}{\sqcup}\|\cE\cG\|\big)\Big]\Big/\sim
\intertext{and}
\big(\fib(\cH)\underset{\fib(\cH_Z)}{\sqcup}\fib(\cG)\big)_1
&\cong
\fib(\cH)_1\underset{\fib(\cH_Z)_1}{\sqcup}\fib(\cG)_1\hspace{1.5cm}\text{\footnotesize (by the same argument as (\ref{ihs}))}\\
&\cong
\Big[\|\cE\cH\|\underset{\cH_0}{\times}\|\cE\cH\|\Big]/{\scriptstyle\sim}
\underset{\textstyle[ \begin{matrix}{\scriptstyle \|\cE\cH_Z\|\underset{Z}{\times}\|\cE\cH_Z\|}
\end{matrix} ]\scriptstyle /\scriptscriptstyle\sim}{\sqcup}
\Big[\|\cE\cG\|\underset{\cG_0}{\times}\|\cE\cG\|\Big]/{\scriptstyle\sim}\\
&\cong
\Big[\|\cE\cH\|\underset{\cH_0}{\times}\|\cE\cH\|
\underset{\|\cE\cH_Z\|\underset{Z}{\times}\|\cE\cH_Z\|}{\sqcup}
\|\cE\cG\|\underset{\cG_0}{\times}\|\cE\cG\|\Big]\Big/\sim.
\end{split}
\]
The result then follows from Lemma \ref{pushpull}.
\qed

\begin{corollary}\label{fcic}\dontshow{fcic}
The fibrant replacement of a cellular topological groupoid is cellular.\qed
\end{corollary}
We also have the following slight variation of Lemma \ref{pcannoy}.
\begin{lemma}
If $\cG$ is paracompact then $\fib(\cG)$ is paracompact. \qed
\end{lemma}

\section{Topological Stacks \rm\dontshow{sec:S}}\label{sec:S}

\subsection{The descent condition \rm\dontshow{s:ds}}\label{s:ds}

For us, a {\em Grothendieck topology} on a category $\C$ consists of a collection of morphisms of $\C$, called {\em covers}, subject to the usual conditions that isomorphisms are covers, pullbacks of covers are covers, and composites of covers are covers (compare with \cite[Chapter III, Definition 2]{MM94}).
%Let $\TT$ be a category equipped with a Grothendieck topology.
%We will mostly be interested in the case in which $\C$ is the category of topological spaces equipped with the \'etale topology (see Definition %\ref{d:cov}).
A {\em lax presheaf of groupoids} on $\C$ is a weak version of a contravariant functor from $\C$ to groupoids (see Section \ref{A:1} for a precise definition).
Given a lax presheaf of groupoids $\X$ on $\C$, an object $T\in\C$, and a cover $U\to T$, we may evaluate $\X$ on the iterated fibered products $U$, $U\times_T U$, and $U\times_T U \times_T U$.
The various projections then induce groupoid morphisms which can be assembled into the following diagram: \dontshow{2dg}
\begin{equation}\label{2dg}
\X(U)\rrarrow \X\big(U\times_T U\big)\rrrarrow \X\big(U\times_T U\times_T U\big).
\end{equation}
Let $\Gamma$ denote the subcategory of the simplicial indexing category $\Delta$ indexing the diagram (\ref{2dg}), so that (\ref{2dg}) describes a lax functor from $\Gamma$ to groupoids. %(see Appendix \ref{A:1} for the definition of a lax functor).

\begin{definition}
A lax presheaf of groupoids $\X$ on $\C$ is a {\em stack} if, for each object $T$ of $\C$ and each cover $U\to T$ of $T$, the natural map \dontshow{llm}
\begin{equation}\label{llm}
\X(T)\to\underset{\Gamma}{\holim}
\Big\{\X(U)\rrarrow \X\big(U\times_T U\big)\rrrarrow \X\big(U\times_T U\times_T U\big)
\Big\}
\end{equation}
is an equivalence of groupoids.
The $2$-category of {\em stacks on $\C$} is the full $2$-subcategory of the category of lax presheaves of groupoids on $\C$ consisting of the stacks.
\end{definition}

\begin{remark}
We do not require that the groupoids $\X(T)$ be small since that would rule out our main examples of stacks.
Instead, we only require that they be essentially small, which is to say (categorically) equivalent to small groupoids.
\end{remark}

From now on, we restrict to the case in which $\C$ is the category of topological
%compactly generated Hausdorff - (not needed!) -
spaces\footnote
{Technically, we should also pick a Grothendieck universe and restrict ourselves to those spaces which are elements of our universe.
}, equipped with the \'etale topology (see Definition \ref{d:cov}).
We refer to the resulting $2$-category of stacks on $\C$ as {\em topological stacks} or simply {\em stacks}.

We adopt the same conventions \caseone, \casetwo $\-$ as for topological groupoids, meaning that in case \casetwo $\-$ we restrict to the sub-$2$-category of faithful morphisms $\Y\to\X$.
In order to explain what this means we introduce a little bit of terminology.
Given a point $p\in\X(*)$ of a stack $\X$, its automorphism sheaf is the sheaf of groups given by
$$
T\mapsto \aut_{\X(T)}(\pr_T^*(p)),
$$
where $\pr_T:T\to *$ is the projection.
A map of sheaves $\Y\to\X$ is said to be a closed inclusion if for every space $X$ and map $X\to\X$
(by which we mean a map from the representable sheaf $\hom(-,X)$),
the pullback $X\times_\X \Y$ is represented by a space $Y$ and the corresponding map $Y\to X$ is a closed inclusion.
Finally, a map of stacks is said to be faithful if it induces closed inclusions on all automorphism sheaves of points.

\begin{remark}
On topological spaces, the \'etale topology and the usual ``open sets''  topology have equivalent theories of sheaves.
It follows that we could have done everything above with the usual topology and we would have obtained an equivalent 2-category of topological stacks.
\end{remark}

To be completely explicit, we give a concrete model for the lax limit (\ref{llm}) used in the definition of a stack.
As before, write $\Gamma$ for the subcategory of $\Delta$ indexing the diagram (\ref{2dg}), let $\X$ be a lax functor from $\Gamma$ to groupoids, and let $\X(0)$, $\X(1)$, $\X(2)$ be the images of the three objects of $\Gamma$.
Then the objects of the groupoid
\[
\underset{\Gamma}{\holim}
\Big\{\X(0)
\put(8.5,5.5){$\scriptstyle d^0$}
\put(8.5,-6){$\scriptstyle d^1$}
\hspace{.05cm}\mbox{\,\put(0,-2){$\longrightarrow$}\put(0,2){$\longrightarrow$}\hspace{.7cm}}
\X(1)
\put(8,8.5){$\scriptstyle d^0$}
\put(8,1){$\scriptstyle d^1$}
\put(8,-7.5){$\scriptstyle d^2$}
\hspace{.05cm}\mbox{\,\put(0,-3){$\longrightarrow$}\put(0,1){$\longrightarrow$}\put(0,5){$\longrightarrow$}\hspace{.7cm}}
\X(2)\Big\}
\]
are pairs $(x,\sigma)$, with $x$ an object of $\X(0)$ and $\alpha:d^1x\to d^0x$ an arrow of $\X(1)$ satisfying
\[
\X(d^0,d^1)d^0\alpha\X(d^0,d^1)^{-1}\X(d^2,d^0)d^2\alpha=
\X(d^1,d^0)d^1\alpha\X(d^1,d^1)^{-1}\X(d^2,d^1),
\]
(here $\X(d^i,d^j)$ denotes the composition 2-morphism $\X(d^j)\X(d^i)\Rightarrow\X(d^j d^i)$ of $\X$, as explained in (\ref{2F1})), and an arrow $(y,\beta)\to(x,\alpha)$ is an arrow $g:y\to x$ of $\X(0)$ with $\alpha(d^1g)=(d^0g)\beta$.

\subsection{Topological stacks from topological groupoids \rm\dontshow{s:tstg}}\label{s:tstg}

Associated to a topological groupoid $\cG$ is the {\em moduli stack} $\M_\cG$ of principal $\cG$-bundles.
That is, we have a functor
$$
\M:\{\text{Topological groupoids}\}\longrightarrow\{\text{Topological stacks}\}
$$
which sends the topological groupoid $\cG$ to the stack $\M_\cG:\{\text{Topological spaces}\}\to\{\text{Groupoids}\}$ defined by
\[
\M_\cG(T):=\{\text{Principal $\cG$-bundles on $T$}\}.
\]
Note that $\M_\cG(T)$ is essentially small since by (\ref{u2}), we have an equivalence
$$
\M_G(T)\simeq \underset{U\in\cov^2(\cH_0)}{\hocolim}\Hom(\pair_T U,\GG),
$$
and the colimit can be taken over any small cofinal collection of covers of $T$, such as those of the form $\coprod U_\alpha\to T$ with $\{U_\alpha\}$ an ordinary open cover of $T$.

\begin{remark}\label{R34}\dontshow{R34}
The functor $\M_\GG$ is the {\em stackification} of the presheaf of groupoids $\Hom(-,\GG)$.
Indeed, the functor $\Hom(-,\GG)$ classifies trivial principal $\GG$-bundles and their isomorphisms.
Provided there exists a space which carries a locally but not globally trivial principal $\GG$-bundle, this is not a stack.
However, the inclusion of globally trivial $\GG$-bundles into locally trivial $\GG$-bundles is initial among presheaf maps from $\Hom(-,\cG)$ with target a stack.
\end{remark}

\begin{example}
The functor $\M_{\cB G}$ assigns to a space $T$ the groupoid of principal $G$-bundles on $T$.
The functor $\M_{\cE G}$ assigns to a space $T$ the groupoid of principal $G$-bundles on $T$ equipped with a global section.
The forgetful functor $\M_{\cE G}\to\M_{\cB G}$ is then a principal $G$-bundle in the category of topological stacks.
\end{example}

Although a morphism of groupoids $\cH\to\cG$ induces a morphism of stacks $\M_\cH\to\M_\cG$,
it is typically not the case that an arbitrary stack map $\M_\cH\to\M_\cG$ comes from a morphism of groupoids.
Usually, $\cH$ must be replaced by $\cH_U$ (see Example \ref{txh}) for some cover $U\to\cH_0$ before such a map can be produced.

\begin{lemma}\label{cii}\dontshow{cii}
Let $\cG$ be a topological groupoid and $U\to \cG_0$ be a cover.
Then the induced map $\M_{\cG_U}\to\M_\cG$ is an equivalence of stacks.
\end{lemma}

\proof
The functor
\[
\qquad\begin{matrix}
\cG_1\acts P\\
\downarrow\downarrow\hspace{.1cm}\swarrow\hspace{.1cm}\downarrow\hspace{.1cm}\\
\cG_0\hspace{.4cm}T\hspace{0cm}
\end{matrix}
\qquad\mapsto\qquad
\begin{matrix}
U\!\times_{\cG_0}\!\cG_1\!\times_{\cG_0}\!U\acts U\!\times_{\cG_0}\!P\\
\hspace{.75cm}\downarrow\downarrow\hspace{.4cm}\swarrow\hspace{.5cm}\downarrow\\
\hspace{.8cm}U\hspace{1.2cm}T
\end{matrix}
\]
which assigns to a principal $\cG$-bundle over $T$ a corresponding principal $\cG_U$-bundle over $T$ has an inverse
\[
\begin{matrix}
U\!\times_{\cG_0}\!\cG_1\!\times_{\cG_0}\!U\acts Q\\
\hspace{1.4cm}\downarrow\downarrow\hspace{.35cm}\swarrow\hspace{.25cm}\downarrow\\
\hspace{1.5cm}U\hspace{.9cm}T
\end{matrix}
\qquad\mapsto\qquad
\begin{matrix}
\cG_1\acts Q/\!\sim\\
\downarrow\downarrow\hspace{.2cm}\swarrow\hspace{.2cm}\downarrow\hspace{.4cm}\\
\cG_0\hspace{.6cm}T\hspace{.3cm}
\end{matrix},\qquad
\]
where the equivalence relation is given by $(u,e,u')\cdot q\sim q$ for every identity arrow $e\in \cG_1$.
\qed
\vspace{.3cm}

We could also prove Lemma \ref{cii} using (\ref{u2}) by arguing that the map
\begin{equation*}
\underset{V\in\cov^2(T)}{\hocolim}\,\Hom(\pair_T V,\cG_U)\longrightarrow\underset{V\in\cov^2(T)}{\hocolim}\,\Hom(\pair_T V,\cG)
\end{equation*}
has an inverse.
It is given by sending a map $\pair_T V\to \cG$ to its pullback along the projection $\cG_U\to\cG$.
As a corollary of Lemma \ref{cii}, we get the following result.

\begin{proposition}\label{ijd}\dontshow{ijd}
The map
\[
\M_\cG\longrightarrow\M_{\fib(\cG)},
\]
induced by the inclusion $\eta:\cG\to\fib(\cG)$, is an equivalence of stacks.
\end{proposition}

\proof
Recall from Section \ref{s:gr} that the space $\|\cG\|=\fib(\cG)_0$ has a canonical cover $U_{\|\cG\|}\to\fib(\cG)_0$.
Thus we may form the pullback groupoid $\fib(\cG)_{U_{\|\cG\|}}$, with object space $(\fib(\cG)_{U_{\|\cG\|}})_0=U_{\|\cG\|}$
and arrow space $(\fib(\cG)_{U_{\|\cG\|}})_1=U_{\|\cG\|}\times_{\fib(\cG_0)}\fib(\cG)_1\times_{\fib(\cG)_0}U_{\|\cG\|}$.
We shall denote the arrows in $\fib(\cG)_{U_{\|\cG\|}}$ by
\[
(r_0\cdot x_0\stackrel{\alpha_1}{\leftarrow}\cdots\stackrel{\alpha_i}{\leftarrow}\underline{r_i\cdot x_i}\stackrel{\alpha_{i+1}}{\leftarrow}\cdots
\stackrel{\alpha_n}{\leftarrow}r_n\cdot x_n\stackrel{\beta}{\leftarrow}
s_0\cdot y_0\stackrel{\gamma_1}{\leftarrow}\cdots\stackrel{\gamma_j}{\leftarrow}\underline{s_j\cdot y_j}\stackrel{\gamma_{j+1}}{\leftarrow}
\cdots\stackrel{\gamma_m}{\leftarrow}s_m\cdot y_m).
\]
The map $\cG\to\fib(\cG)_{U_{\|\cG\|}}$, which sends $x$ to $(\underline{1\cdot x})$ and $x\stackrel{\alpha}{\leftarrow} y$ to
$(\underline{1\cdot x}\stackrel{\alpha}{\leftarrow}\underline{1\cdot y})$, has a categorical inverse $\fib(\cG)_{U_{\|\cG\|}}\to\cG$ given by
$$
(r_0\cdot x_0\stackrel{\alpha_1}{\leftarrow}\cdots\stackrel{\alpha_i}{\leftarrow}\underline{r_i\cdot x_i}\stackrel{\alpha_{i+1}}{\leftarrow}\cdots\stackrel{\alpha_n}{\leftarrow}r_n\cdot x_n)\mapsto x_i
$$
and
\[
\begin{matrix}
(r_0\cdot x_0\stackrel{\alpha_1}{\leftarrow}\cdots\stackrel{\alpha_i}{\leftarrow}\underline{r_i\cdot x_i}\stackrel{\alpha_{i+1}}{\leftarrow}\cdots\stackrel{\alpha_n}{\leftarrow}r_n\cdot x_n \stackrel{\beta}{\leftarrow}
s_0\cdot y_0\stackrel{\gamma_1}{\leftarrow}\cdots\stackrel{\gamma_j}{\leftarrow}\underline{s_j\cdot y_j}\stackrel{\gamma_{j+1}}{\leftarrow}
\cdots\stackrel{\gamma_m}{\leftarrow}s_m\cdot y_m)
\\
\raisebox{10pt}{\rotatebox{-90}{$\mapsto$}}
\\
x_i\put(2,6){$\scriptstyle\alpha_{i+1}\cdots\alpha_n\beta\gamma_1\cdots\gamma_j$}
\longleftarrow\!\!\!-\!\!\!-\!\!\!-\!\!\!-\!\!\!-\!\!\!-\!\!\!-\!\!\!-\!\!\!-\:y_j
\end{matrix}.
\]
Together with Lemma \ref{cii}, this gives equivalences $\M_\cG\simeq\M_{\fib(\cG)_{U_{\|\cG\|}}}\simeq\M_{\fib(\cG)}$.
%and the result follows by Lemma \ref{cii} since $\M_{\fib(\cG)_{U_{\|\cG\|}}}\simeq\M_{\fib(\cG)}$.
\qed
\vspace{.3cm}

Given two topological groupoids $\cH$ and $\cG$, we wish to have a concrete model for the groupoid of stack maps from $\M_\cH$ to $\M_\cG$.
Thus we define
\[
\M_\cG(\cH):=\begin{cases}
\big\{\text{Principal $\cG$-bundles on $\cH$}\big\}\qquad&\text{in case \caseone,}\\
\big\{\text{Faithful principal $\cG$-bundles on $\cH$}\big\}\qquad&\text{in case \casetwo,}
\end{cases}
\]
where the faithful bundles are those defined in Remark \ref{fb}.
Combining (\ref{u2}) and (\ref{u4}), this can then be rewritten as \dontshow{sf}
\begin{equation}\label{sf}
\M_\cG(\cH)\simeq\underset{U\in\cov^2(\cH_0)}{\hocolim}\,\Hom(\cH_U,\cG).
\end{equation}
%provided one interprets the colimit $2$-categorically.

Using Lemma \ref{cii} we can produce a map from $\M_\cG(\cH)$ to $\Hom(\M_\cH,\M_\cG)$ by composing the right hand side of (\ref{sf})
with the inverse of $\M_{\cH_U}\to\M_\cH$.
In terms of principal bundles, an element \dontshow{pnt}
\begin{equation}\label{pnt}
\begin{matrix}
\cG_1\acts P\acted \cH_1\\
\downarrow\downarrow\hspace{.1cm}\swarrow\hspace{.3cm}\searrow\hspace{.1cm}\downarrow\downarrow\hspace{.1cm}\\
\cG_0\hspace{1.17cm}\cH_0
\end{matrix}\quad\in\hspace{.1cm} \M_\cG(\cH)
\end{equation}
would then act by \dontshow{hmt}
\begin{equation}\label{hmt}
\begin{matrix}
\cH_1\acts Q\\
\downarrow\downarrow\hspace{.1cm}\swarrow\hspace{.15cm}\downarrow\hspace{.05cm}\\
\hspace{0cm}\cH_0\hspace{.4cm}T
\end{matrix}
\qquad
\mapsto
\qquad
\begin{matrix}
\cG_1\acts (P\times_{\cH_0}Q)/\cH_{\phantom{1}}\\
\downarrow\downarrow\hspace{.3cm}\swarrow\hspace{.8cm}\downarrow\hspace{1.1cm}\\
\cG_0\hspace{1.35cm}T\hspace{.95cm}.
\end{matrix}
\end{equation}
The map $\M_\cG(\cH)\to\Hom(\M_\cH,\M_\cG)$ is an equivalence of groupoids whose
inverse $\Hom(\M_\cH,\M_\cG) \to \M_\cG(\cH)$ is obtained as follows.
Given a natural transformation $p\in\Hom(\M_\cH,\M_\cG)$, we can evaluate it
on the principal $\cH$-bundle $\cH_1\to\cH_0$ to get a principal $\cG$-bundle $P\to\cH_0$.
The isomorphism $d_0^*(\cH_1\to\cH_0)\simeq d_1^*(\cH_1\to\cH_0)$ provides an isomorphism $d_0^*P\simeq d_1^*P$
which by naturality satisfies the required cocycle condition, resulting in a principal $\cG$-bundle on $\cH$.

\begin{proposition}\label{BGstr}\dontshow{BGstr}
Let $P$ be a principal $G$-bundle over a non-empty space $T$.
Then the stack associated to the gauge groupoid of $P$ is equivalent to the stack associated to $\cB G$.
Moreover, if $\M_\cG$ is equivalent to $\M_{\cB G}$ then $\cG\cong\gauge_G(P)$ for some principal $G$-bundle $P$.
\end{proposition}

\proof
A point $p\in P$ induces a morphism of principal $G$-bundles $G\to P$, and hence a morphism of groupoids $\cB G\cong\gauge_G(G)\to\gauge_G(P)$.
On the level of stacks, this map has an inverse, namely the map classifying the principal $G$-bundle $\pair(P)\to\gauge_G(P)$.

Conversely, suppose $\cG$ is a groupoid equipped with an equivalence $\varphi:\M_{\cG}\to\M_{\cB G}$.
Then $\varphi$ classifies a principal $G$-bundle $P$ on $\cG$, corresponding to a cover $U\to\cG_0$ and a groupoid map $\cG_U\to\cB G$.
Let $\cP:=(d_1^*P\rrarrow P)$ be the groupoid with source and target given by action (\ref{aktm}) and projection, respectively, so that $\cP\to\cG$ is a principal $G$-bundle in the category of topological groupoids.
The result is a (lax) cartesian square
$$
\xymatrix
{\M_\cP \ar[r]\ar[d] & \M_{\cE G}\ar[d]\\
\M_\cG\ar[r] & \M_{\cB G}.}
$$
in which the horizontal maps are equivalences.
But $\M_{\cE G}\cong\M_{\pair(G)}$ is equivalent to the trivial stack, and therefore so is $\M_\cP$.
At the level of groupoids, this means that the canonical map $\cP\to\pair(P)$ is an isomorphism.
Hence
\[
\cG\cong\cP/G\cong\pair(P)/G\cong\gauge_G(P)
\]
is the gauge groupoid of the principal $G$-bundle $P\to\cG_0$.
\qed

\begin{corollary}\label{a=g}\dontshow{a=g}
If $H$ is a closed subgroup of $G$ satisfying the assumptions of Lemma \ref{jspq},
then the stack associated to $G\ltimes (G/H)$ is equivalent to $\cB H$.
\end{corollary}

\proof
By Lemma \ref{jspq}, the action groupoid $G\ltimes (G/H)$ is equivalent to the gauge groupoid of the principal $H$-bundle $G\to G/H$.
\qed

\subsection{Fibrant topological groupoids are stacks \rm\dontshow{s:AS}}\label{s:AS}

As explained in Remark \ref{R34}, the functor $\Hom(-,\cG)$ is typically not a stack.
It turns out, though, that if $\cG$ is fibrant (and provided we restrict to the subcategory of paracompact spaces) then $\Hom(-,\cG)$ is a stack.

\begin{proposition}\label{ppda}\dontshow{ppda}
Let $\cH$ be a topological groupoid with paracompact object space.
Then for any topological groupoid $\cG$, there is a natural equivalence of groupoids
\[
\Hom(\cH,\fib\cG)\simeq \M_\cG(\cH).
\]
\end{proposition}

\proof
Given an object of $\M_\cG(\cH)$, classified by a cover $f:U\to\cH_0$ and a map $\rho:\cH_U\to\cG$, we construct a morphism $\theta:\cH\to\fib(\cG)$.
Since $\cH_0$ is paracompact, we may appropriately refine $U\to\HH_0$ and pick a partition of unity $(U,\varphi,<)$.
Given $y\in\cH_0$, let $\{y_i\}_{i=0}^n$ be its preimages under $f$, let $x_i:=\rho(y_i)\in\cG_0$, and let $\alpha_i\in\cG_1$ be the images of the morphisms
$(y_{i-1},1_y,y_i)\in U\times_{\cH_0}\cH_1\times_{\cH_0}U=(\cH_U)_1$ under $\rho$.
We then define $\theta(y)\in\fib(\cG)_0$ by
\[
\theta(y):=(\varphi(y_0)\cdot x_0\stackrel{\alpha_1}{\leftarrow}\cdots\stackrel{\alpha_n}{\leftarrow}\varphi(y_n)\cdot x_n).
\]
Given a point $z\in \cH_0$, with image $\theta(z)=(\varphi(z_0)\cdot v_0\stackrel{\gamma_1}{\leftarrow}\cdots\stackrel{\gamma_m}
{\leftarrow}\varphi(z_m)\cdot v_m)$,
and a morphism $g:y\leftarrow z$, we let $\beta:=\rho(y_n,g,z_0)$ and define
\[
\theta(g):=(\varphi(y_0)\cdot x_0\stackrel{\alpha_1}{\leftarrow}\cdots\stackrel{\alpha_n}{\leftarrow}\varphi(y_n)\cdot x_n\stackrel{\beta}{\leftarrow}\varphi(z_0)\cdot v_0\stackrel{\gamma_1}{\leftarrow}\cdots\stackrel{\gamma_m}{\leftarrow}\varphi(z_m)\cdot v_m).
\]
This construction identifies $\Hom(\cH,\fib\cG)$ with a groupoid whose objects are pairs consisting of a partition of unity $(U,\varphi,<)$ and a homomorphism $\rho:\HH_U\to\GG$.
This groupoid also models the $2$-categorical colimit $\hocolim\Hom(\cH_U,\cG)$, which in turn is equivalent to $\M_\GG(\HH)$ by (\ref{sf}).
\qed

\begin{corollary}\label{paf}\dontshow{paf}
If $\cG$ is fibrant and $\cH$ has paracompact object space, then
\[
\Hom(\cH,\cG)\simeq\Hom(\M_\cH,\M_\cG).
\]
\end{corollary}
\proof
By Lemma \ref{hece}, the groupoids
$\cG$ and $\fib\cG$ are categorically equivalent, and thus $\Hom(\cH,\cG)\simeq\Hom(\cH,\fib\cG)$.
The result follows since $\Hom(\M_\cH,\M_\cG)\simeq \M_\cG(\cH)$.
\qed
\vspace{.3cm}

In the special case where $\cH$ is equivalent to the unit groupoid of a topological space $T$ (see Example \ref{extp}), the above corollary
specializes an equivalence
\[
\Hom(T,\fib\cG)\simeq\M_\cG(T),
\]
and we obtain the following result.

\begin{corollary}
Let $\cG$ be a fibrant topological groupoid.
Then the functor $T\mapsto \Hom(T,\cG)$ is a stack on the category of paracompact topological spaces.\qed
\end{corollary}

\subsection{The homotopy type of a topological stack \rm\dontshow{sec:hts}}\label{sec:hts}

As with spaces, the {\em weak homotopy type} of a stack is its best possible approximation by a CW-complex.

\begin{definition}
The {\em category of CW-complexes over a stack} $\X$ is the 2-category whose objects are pairs $(X,\varphi)$, with $X$ a CW-complex and $\varphi:\M_X\to\X$ a stack map, and whose arrows from $\psi:\M_Y\to\X$ to $\varphi:\M_X\to\X$ are pairs $(f,\alpha)$, with $f:Y\to X$ a map of CW-complexes and $\alpha$ a $2$-morphism $\varphi\circ f\Rightarrow\psi$.
\end{definition}

Even though the category of CW-complexes over $\X$ is not topologically enriched,
there's still a good notion of homotopy between morphisms $(f,\alpha),(g,\beta):(\psi:\M_Y\to\X)\to(\varphi:\M_X\to\X)$.
Namely, it's a morphism $(h,\gamma)$ from $\psi\circ\pr_1:\M_{Y\times[0,1]}\to\M_Y\to\X$ to $\varphi:\M_X\to\X$,
such that $h$ is a homotopy from $f$ to $g$ and $\gamma:\varphi\circ h\Rightarrow \psi\circ\pr_1$ is a
2-morphism whose restriction to the endpoints of $[0,1]$ are $\alpha$ and $\beta$, respectively.
We then let $\mathrm{Ho}(\mathrm{CW}/\X)$ be the homotopy category of CW-complexes over $\X$.
Its objects are CW-complexes over $\X$, and its morphisms are homotopy classes in the above sense.

\begin{definition}
The weak homotopy type of a stack $\X$ is the homotopy type of a terminal object of the homotopy category $\mathrm{Ho}(\mathrm{CW}/\X)$ of CW-complexes over $\X$.
\end{definition}

\begin{remark}
All stacks possess weak homotopy types, even though we shall only prove it for those of the form $\M_\cG$.
\end{remark}

\begin{remark}
A {\em homotopy terminal} object of a topologically enriched category is an object $T$ such that, for any other object $U$, the space $\map(U,T)$ is contractible.
Clearly a homotopy terminal object $T$ is terminal in the homotopy category; for the converse, it is enough to also require that our topologically enriched category be tensored over topological spaces.
Indeed, if $T$ is terminal in the homotopy category, then for any object $U$ we have that
\[
[S^n,\map(U,T)]\cong[S^n\otimes U,T]\cong*
\]
and thus that $\map(U,T)$ is contractible.
\end{remark}

\begin{remark}
The category of CW-complexes over a stack $\X$ is more naturally enriched over the category of sheaves of sets on topological spaces.
However, it still makes sense to talk about weakly contractible sheaves, and if $T$ is a terminal object of $\mathrm{Ho}(\mathrm{CW}/\X)$ then the sheaf $\map_{\X}(-,\M_T)$ is indeed weakly contractible.
\end{remark}

\begin{remark}
It is also possible to enrich topological stacks over the homotopy category of topological spaces.
For this, we require the notion of the mapping stack:
given stacks $\X$ and $\Y$, let $\stackMap(\Y,\X)$ denote the stack
$$
\stackMap(\Y,\X)(T):=\Hom(\M_T\times\Y,\X).
$$
Then we simply define $\map(\Y,\X)$ to be the weak homotopy type of the mapping stack $\stackMap(\Y,\X)$.
\end{remark}

\begin{proposition}
The weak homotopy type of $\M_\cG$ is the weak homotopy type of $\|\cG\|$.
In particular, the weak homotopy type of $\M_\cG$ exists.
\end{proposition}

\proof
There is a canonical map $\M_{\|\cG\|}\to\M_\cG$ which classifies the universal principal $\cG$-bundle $\|\cE\cG\|\to\|\cG\|$.
Assuming for the moment that $\|\cG\|$ is a CW-complex, this map is an object of the category of CW-complexes over $\M_\cG$, which we claim is homotopy terminal.
Indeed, if $T$ is a CW-complex equipped with a map $\M_T\to\M_\cG$, then, since $\|\cG\|$ (or more precisely its image in the homotopy category)
classifies isomorphism classes of principal $\cG$-bundles over paracompact spaces,
we see that there is precisely one homotopy class of map $T\to\|\cG\|$ over $\M_\cG$.

If $\|\cG\|$ is not a CW-complex, then a CW-replacement $\cof(\|\cG\|)$ of $\|\cG\|$ has the property that any map $T\to\|\cG\|$ with $T$ a CW-complex
factors uniquely up to homotopy through the canonical map $\cof(\|\cG\|)\to\|\cG\|$.
It follows that $\cof(\|\cG\|)$ is terminal in $\mathrm{Ho}(\mathrm{CW}/\M_\cG)$.
\qed
\vspace{.3cm}

There is a similar notion of a CW-complex over a topological groupoid.
More precisely, for any topological groupoid $\cG$, define the homotopy category $\mathrm{Ho}(\mathrm{CW}/\cG)$ of CW-complexes over $\cG$ as the quotient category of $\mathrm{Ho}(\mathrm{CW}/\cG_0)$ obtained by identifying two morphisms which differ by a natural transformation.

\begin{theorem}\label{wty}\dontshow{wty}
Let $\cG$ be a fibrant groupoid and $\cH$ a groupoid with paracompact object space.
Then the spaces $\map(\cH,\cG)$ and $\map(\M_\cH,\M_\cG)$ are naturally weakly homotopy equivalent.
\end{theorem}

\proof
For any CW-complex $T$ we have groupoid equivalences
\[
\Hom\big(\M_T,\stackMap(\M_\cH,\M_\cG)\big)\simeq\Hom(\M_T\times\M_\cH,\M_\cG)\simeq\Hom(T\times\cH,\cG)\simeq\Hom\big(T,\Map(\cH,\cG)\big),
\]
where the second equivalence follows from Corollary \ref{paf}.
%Proposition \ref{ppda}.
Passing to equivalence classes of objects and then modding out by homotopies, it follows that the homotopy categories of CW-complexes over
$\stackMap(\M_\cH,\M_\cG)$ and over $\Map(\cH,\cG)$ are equivalent.

To finish the proof, we need to identify the homotopy category of CW-complexes over the topological groupoid $\Map(\cH,\cG)$ with
the homotopy category of CW-complexes over the topological space $\map(\cH,\cG)$.
This amounts to showing that two naturally isomorphic maps $f,g:T\to\Map(\cH,\cG)$ are necessarily also homotopic.
We show instead that the adjoint maps $\tilde f,\tilde g:T\times\cH\to \cG$ are homotopic.
Indeed, the homotopy
\[
T\times\cH\times[0,1]\to \cG
\]
between $\tilde f$ and $\tilde g$ is obtained by composing the homotopy $T\times\cH\times[0,1]\to \fib(\cG)$ constructed in (\ref{aqy})
with the retraction $\fib(\cG)\to \cG$ given by Proposition \ref{ccFR}.
\qed

\section{$\Orb$-spaces \rm\dontshow {sec:OSp}}\label{sec:OSp}

\subsection{The category $\Orb$ \rm\dontshow{catorb}}\label{catorb}

Given a topological group $G$, recall that $\cU G:=\fib(\cB G)\cong\gauge_G(EG)$ is the universal gauge groupoid of $G$.
Let $\torb$ be the essentially small topological category whose class of objects is our fixed family $\ccF$ of allowed isotropy groups and whose morphisms are given by
\dontshow{opu}
\begin{equation}\label{opu}
\torb(H,G):=\map(\cU H,\cU G).
\end{equation}
Note that the right hand side of (\ref{opu}) is either the space of all maps $\cU H\to\cU G$ or the subspace of faithful maps $\cU H\to\cU G$, depending upon whether we are in case \caseone \- or \casetwo.
%Also recall that, depending if we are in case \caseone \- or \casetwo, the right hand side of (\ref{opu}) means either all maps $\cU H\to \cU G$
%or faithful maps $\cU H\to\cU G$.

In applications it may be desirable to replace $\torb$ by an equivalent but better behaved topological category.
We therefore fix, once and for all, an essentially small topological category $\Orb$ and a weak equivalence $\Orb\to\torb$
(see Definition \ref{da4}).
Abusing terminology, we refer to the objects of $\Orb$ as {\em allowed isotropy groups} and designate them by the same letters as their images in $\torb$.

\begin{definition}
An {\em $\Orb$-space} is a contravariant functor from $\Orb$ to spaces.
\end{definition}

\begin{remark}
The main theorems of this paper assert that $\Orb$-spaces comprise a good homotopy theoretic model of orbispaces.
Thus the precise choice of indexing category is not important, since by Lemma \ref{Aee} two different choices of $\Orb$ induce Quillen equivalent model categories of $\Orb$-spaces.
\end{remark}

%We have avoided choosing a particular indexing category $\Orb$ in order to have the upmost flexibility in applications.
%One possibility which always works is to let $\Orb=\torb$.
\begin{remark}
The choice $\Orb=\torb$ is conceptually satisfying but annoying to work with in practice.
A smaller model for $\Orb$ is given by \dontshow{vd}
\begin{equation}\label{vd}
\Orb(H,G):=\|\Map(\cB H,\cB G)\|.
\end{equation}
Here, composition is defined as the retraction from Lemma \ref{Kri}
\[
\|\Map(\cB K,\cB H)\|\times\|\Map(\cB H,\cB G)\|\to\|\Map(\cB K,\cB H)\times\Map(\cB H,\cB G)\|
\]
followed by the realization
\(
\|\Map(\cB K,\cB H)\times\Map(\cB H,\cB G)\|\to\|\Map(\cB K,\cB G)\|
\)
of the composition in topological groupoids,
and the functor $\Orb\to\torb$ is given by restricting the map \dontshow{vp}
\begin{equation}\label{vp}
\fib\big(\Map(\cB H,\cB G)\big)\times \fib(\cB H)\stackrel{r}{\longrightarrow}\fib\big(\Map(\cB H,\cB G)\times \cB H\big)\to\fib(\cB G).
\end{equation}
to the objects of $\fib(\Map(\cB H,\cB G))$.
\end{remark}

\begin{lemma}\label{l:vd}\dontshow{l:vd}
The map
\[
\|\Map(\cB H,\cB G)\|\to\map(\cU H,\cU G)
\]
defined by (\ref{vp}) is a homotopy equivalence.
\end{lemma}

\proof
The map $\|\Map(\cB H,\cB G)\|\to\map(\cB H,\cU G)$
given by sending $(r_0\cdot\varphi_0\stackrel{g_1}{\leftarrow}\cdots\stackrel{g_n}{\leftarrow}r_n\cdot\varphi_n)$ to the map
\[
\begin{split}
\cB H\to\cU G:\:
*&\mapsto (r_0\cdot *\stackrel{g_1}{\leftarrow}\cdots\stackrel{g_n}{\leftarrow}r_n\cdot *)\\
h&\mapsto
(r_0\cdot *\stackrel{g_1}{\leftarrow}\cdots\stackrel{g_n}{\leftarrow}r_n\cdot *
\put(2,6){$\scriptstyle
(g_1\ldots g_n)^{-1}\varphi_0(h)
$}
\longleftarrow\!\!\!-\!\!\!-\!\!\!-\!\!\!-\!\!\!-\!\!\!-\!\!\!-\!\!\!-\:
r_0\cdot *\stackrel{g_1}{\leftarrow}\cdots\stackrel{g_n}{\leftarrow}r_n\cdot *)
\end{split}
\]
is an isomorphism.
The result follows since we have a commutative diagram
\[
\xymatrix@R=.5cm{
\|\Map(\cB H,\cB G)\|\ar[rr]^{\simeq}\ar[dr]&&\map(\cB H,\cU G)\\
&\map(\cU H,\cU G)\ar[ur]
}
\]
in which the restriction $\map(\cU H,\cU G)\to \map(\cB H,\cU G)$ is an equivalence by Proposition \ref{fibfactor}.
\qed
\vspace{.3cm}

\begin{remark}
We have isomorphisms of topological groupoids $\Map(\cB H,\cB G)\cong\map(H,G)\rtimes G$, where $\map(H,G)$ denotes the space of group
homomorphisms (monomorphisms in case \casetwo) from $H$ to $G$, and $G$ acts on $\map(H,G)$ by conjugation.
Using equation (\ref{x:hq}), we may rewrite (\ref{vd})\dontshow{hqv}
\begin{equation}\label{hqv}
\Orb(H,G):=\|\Map(\cB H,\cB G)\|\cong\map(H,G)\times_G EG
\end{equation}
as the Borel construction for the conjugation action of $G$ on $\map(H,G)$.
\end{remark}

A continuous contravariant from $\Orb$ to spaces is the same thing as a space with a continuous left action by $\Orb$, regarded as a category object in the category of spaces\footnote{For this to make sense, $\Orb$ should be taken to be genuinely small instead of just essentially small.}.
Adopting our notational conventions for topological groupoids, write $\Orb=(\Orb_1\rrarrow\Orb_0)$, where $\Orb_0$ is the (discrete) space of objects of $\Orb$ and
$$
\Orb_1:=\coprod_{H,G\in\Orb_0}\Orb(H,G)
$$
is the space of arrows of $\Orb$.

Viewing $\Orb_0$ as a category with only identity arrows, the inclusion of $\Orb_0$ into $\Orb$ induces a forgetful functor from $\Orb$-spaces to $\Orb_0$-spaces, the latter being isomorphic to the category of spaces over $\Orb_0$.
Left adjoint to the forgetful functor is the {\em free} $\Orb$-space functor $\FF$, defined by
\[
\FF(T):=\Orb_1\times_{\Orb_0}T,
\]
with left $\Orb$-action given by
$$
\Orb_1\times_{\Orb_0}\FF(T)\simeq\Orb_1\times_{\Orb_0}\Orb_1\times_{\Orb_0}T\stackrel{m\times 1}{\longrightarrow} \Orb_1\times_{\Orb_0}T\simeq\FF(T).
$$
In terms of functors, the value of $\FF(T)$ on the orbit $H$ is
$$
\FF(T)(H):=\coprod_{G\in\Orb_0}\Orb_1(H,G)\times T(G),
$$
where $T(G)$ is value of the $\Orb_0$-space $T$ at $G$, which is to say the fiber of $T\to\Orb_0$ over $G\in\Orb_0$.

For the reader's convenience, we briefly recall a topological version of the Bousfield-Kan,
or projective, model structure on the category of continuous functors from a small topologically-enriched category to a cofibrantly generated topological model category (see \cite{Hir03} for the basics on cofibrantly generated model categories).
The weak equivalences and the fibrations are defined by evaluation on objects of the indexing category,
and the generating (trivial) cofibrations are obtained by applying free functors (left adjoint to evaluation on objects) to the usual generating (trivial) cofibrations for spaces.

More precisely, let $\delta_G$ denote the $\Orb_0$-space defined by $\delta_G(H)=\emptyset$ if $H\not = G$ and $\delta_G(H)=*$ if $H=G$.
Then the generating cofibrations and generating trivial cofibrations are given by \dontshow{gtw}
\begin{equation}\label{gtw}
\FF(S^{n-1}\times\delta_G)\rightarrow \FF(D^n\times\delta_G)\qquad
\text{and}\qquad \FF(D^n\times\delta_G) \rightarrow \FF(D^n\times[0,1]\times\delta_G)
\end{equation}
respectively.

One consequence of this formalism is the existence of a cofibrant replacement functor
\[
\cof:\{\text{$\Orb$-spaces}\}\longrightarrow\{\text{$\Orb$-spaces}\}.
\]
It comes equipped with a natural transformation $\epsilon:\cof\to 1$ such that for each $\Orb$-space $X$ the map $\epsilon(X):\cof(X)\to X$
is a weak equivalence of $\Orb$-spaces, and the resulting two maps $\cof^2(X)\to\cof(X)$ are homotopy equivalences.
Moreover, $\cof(X)$ is always a cellular object in the sense of cofibrantly generated model categories,
meaning that it's a colimit of pushouts along maps of the form $\FF(S^{n-1}\times\delta_G)\rightarrow \FF(D^n\times\delta_G)$.

\subsection{Comparison with the category of topological groupoids {\rm\dontshow{sec:CwG}}}\label{sec:CwG}

In this section we construct a pair of adjoint functors
\[
L:\:\big\{\Orb\text{-spaces}\big\}\put(7,3){$\longrightarrow$}\put(7,-2){$\longleftarrow$}\hspace{1cm}\big\{\text{Topological groupoids}\big\}\::R
\]
which, suitably derived, will implement the equivalence of homotopy theories between $\Orb$-spaces
and topological groupoids.

The right adjoint is easy: a topological groupoid $\cG$ gives rise to an $\Orb$-space
\begin{equation*}
R(\cG):G\mapsto \map(\cU G,\cG)
\end{equation*}
represented by $\cG$, and the action $\Orb(H,G)\times \map(\cU G,\cG)\to \map(\cU H,\cG)$ is given by the reference functor $\Orb\to\Orb'$.
Since $R$ is only well-behaved on the subcategory of fibrant groupoids, we write
\[
\RR(\cG):= R (\fib(\cG))
\]
for the $\Orb$-space represented by the fibrant replacement of $\cG$.
Even though we are not really using the formalism of model categories in the context of topological groupoids, we think of the composite $\RR = R\circ\fib$ as the {\em right derived} functor of $R$.

\begin{lemma}\label{freeres}\dontshow{freeres}
An $\Orb$-space $X$ is canonically the coequalizer
%of the arrows
%$
%\FF^2(X)\rrarrow\FF(X),
%$
of the pair of maps
$$
\FF^2(X)=\Orb_1\times_{\Orb_0}\Orb_1\times_{\Orb_0} X\rrarrow\Orb_1\times_{\Orb_0} X=\FF(X)
$$
induced by the action $\Orb_1\times_{\Orb_0} X\to X$ and the
multiplication $\Orb_1\times_{\Orb_0}\Orb_1\to\Orb_1$.
\end{lemma}

\proof
Since $\FF(X)\to X$ coequalizes $\FF^2(X)\rrarrow\FF(X)$, we get a map
from the coequalizer to $X$.
That map is split surjective since there's a section $X\to \FF(X):x\mapsto (1,x)$.
To see that it's injective, consider another preimage $(f,y)\in\FF(X)$ of $x$.
Then $(1,f,y)\in\FF^2(X)$ maps to $(f,y)$ and $(1,x)$, thus proving that they're equal in the coequalizer.
\qed
\vspace{.3cm}

We now describe the left adjoint $L$ of the functor $R$.
By the above lemma, $L$ is determined by its restriction to the full subcategory of free $\Orb$-spaces by the formula
\[
L(X)=\coeq\big(L\FF^2(X)\rrarrow L\FF(X)\big).
\]
Moreover, as any $\Orb_0$-space is a coproduct (indexed by a space, not just a set) of the $\Orb$-spaces $\delta_G$, it's enough to compute $L\FF(\delta_G)$.
But the groupoid $L\FF(\delta_G)$ must satisfy
\[
\begin{split}
\map\big(L\FF(\delta_G),\cG\big)&=\map_{\Orb}\big(\FF(\delta_G),R\cG\big)\\
&=\map_{\Orb_0}(\delta_G,R\cG)
=R\cG(G)=\map(\cU G,\cG),
\end{split}
\]
so we must have $L\FF(\delta_G)=\cU G$.
Letting $\cU:=\coprod_{G\in\Orb_0}\cU G$ be the disjoint union of the universal gauge groupoids,
it follows that
\[
L\FF(T)=\coprod_{G\in\Orb_0}\cU G\times T(G)=\cU\times_{\Orb_0} T,
\]
and that for a general $\Orb$-space $X$,
\begin{equation}
LX=\coeq\Big(\cU\times_{\Orb_0}\Orb_1\times_{\Orb_0}X\rrarrow\cU\times_{\Orb_0}X\Big).
\end{equation}
We write $\LL$ for the composite functor $L\circ\cof$, and regard it as the {\em left derived} functor of $L$.

\begin{lemma}\label{Lcm}\dontshow{Lcm}
The functor $L$ commutes with attaching cells; that is, the natural map
\[
L(X)\underset{S^{n-1}\times\cU G}{\sqcup}D^n\times\cU G\,\,\to\,\, L\Big(X\underset{\FF(S^{n-1}\times\delta_G)}{\sqcup}\FF(D^n\times\delta_G)\Big)
\]
is an isomorphism.
\end{lemma}

\proof
Since $L$ is a left adjoint, it commutes with pushouts.
The result then follows from the isomorphisms
$L\FF(S^{n-1}\times\delta_G)\cong S^{n-1}\times \cU G$ and $L\FF(D^n\times\delta_G)\cong D^n\times \cU G$.
\qed

\begin{corollary}\label{cnf}\dontshow{cnf}
For any $\Orb$-space $X$, the groupoid $\LL X$ is cellular. \qed
\end{corollary}

\begin{lemma}\label{heLX}\dontshow{heLX}
For any $\Orb$-space $X$, the canonical map $\eta:\LL X\to\fib\LL X$ is a homotopy equivalence of groupoids.
\end{lemma}

\proof
Since $\cof(X)$ is a cellular $\Orb$-space, we can write it as
$\colim_\alpha\cof_\alpha(X)$,
where $\alpha$ runs through the set of ordinals smaller than a given ordinal $\alpha_0$.
We then have
\[
\cof_{\alpha+1}(X)=\cof_\alpha(X)\underset{\FF(S^{n-1}\times\delta_G)}{\sqcup}\FF(D^n\times\delta_G)
\qquad\text{and}\qquad
\cof_\beta(X)=\underset{\alpha<\beta}{\colim}\,\cof_\alpha(X)
\]
for ordinals $\alpha<\alpha_0$ and limit ordinals $\beta<\alpha_0$.
We show by induction that for $\alpha\le \alpha_0$ the map \dontshow{bls}
\begin{equation}\label{bls}
L\cof_\alpha(X)\longrightarrow\fib L \cof_\alpha(X)
\end{equation}
is a homotopy equivalence, and that the homotopy inverse and the two homotopies between the composites and the identities can be taken
compatibly with those for smaller $\alpha$.

So let's assume that (\ref{bls}) is a homotopy equivalence.
We must show that
\[
L\cof_{\alpha+1}(X)\longrightarrow\fib L \cof_{\alpha+1}(X)
\]
is also a homotopy equivalence, and that the homotopy inverse and two homotopies can be taken compatibly with those of (\ref{bls}).
By Lemmas \ref{Lcm} and \ref{sip}, we can write
\[
\begin{split}
L\cof_{\alpha+1}(X)&\cong
L\cof_{\alpha}(X)\underset{S^{n-1}\times\cU G}{\sqcup}D^n\times\cU G,\\
\fib L\cof_{\alpha+1}(X)&\cong
\fib L\cof_{\alpha}(X)\underset{S^{n-1}\times\fib \cU G}{\sqcup}D^n\times\fib \cU G.
\end{split}
\]
By corollary \ref{fibsq}, the map $\cU G = \fib \cB G \to \fib^2 \cB G = \fib \cU G$ is a homotopy equivalence.
The result then follows from Lemma \ref{ytb}.

We still need to show that, assuming (\ref{bls}) is a homotopy equivalence for all $\alpha$ smaller than a given limit ordinal $\beta$,
then it's also a homotopy equivalence for $\beta$.
But this is trivial given the compatibility conditions that we assumed for the homotopy inverses and the two homotopies.
\qed
\vspace{.3cm}

Recall that a map of $\Orb$-spaces $Y\to X$ is a weak equivalence if for every group $G\in\Orb_0$, the evaluation $Y(G)\to X(G)$ is
a weak equivalence of spaces.

\begin{lemma}\label{isco}\dontshow{isco}
For any $\Orb$-space $X$, the canonical map $\cof (X)\to RL\cof (X)$ is a weak equivalence.
\end{lemma}

\proof
Write $\cof(X)$ as $\colim_\alpha\cof_\alpha(X)$, with
\[
\cof_{\alpha+1}(X)=\cof_\alpha(X)\underset{\FF(S^{n-1}\times\delta_G)}{\sqcup}\FF(D^n\times\delta_G).
\]
Assuming by induction that $\cof_\alpha(X)\to RL\cof_\alpha(X)$ is a weak equivalence,
we must show that the same holds for $\cof_{\alpha+1}(X)\to RL\cof_{\alpha+1}(X)$.
In other words, we must show that for any $H\in\Orb_0$,
the evaluation $\cof_{\alpha+1}(X)(H)\to RL\cof_{\alpha+1}(X)(H)$
is a weak equivalence of spaces.
Evaluating $\cof_{\alpha+1}(X)$ at $H$, we get
\begin{alignat}{1}
\cof_{\alpha+1}(X)(H)&\cong\cof_\alpha(X)(H)\underset{S^{n-1}\times\Orb(H,G)}\sqcup D^n\times\Orb(H,G).\nonumber\\
\intertext{Evaluating $RL\cof_{\alpha+1}(X)$ at $H$, we get}
RL\cof_{\alpha+1}(X)(H)&\cong\map\big(\cU H,L\cof_{\alpha+1}(X)\big)\nonumber\\
&\cong\map\Big(\cU H,L\cof_\alpha(X)\underset{S^{n-1}\times\cU G}{\sqcup}D^n\times\cU G\Big)\label{azu}\\
&\approx\map\big(\cU H,L\cof_\alpha(X)\big)\underset{S^{n-1}\times\map(\cU H,\cU G)}{\sqcup}D^n\times\map(\cU H,\cU G)\nonumber\\
&\cong RL\cof_\alpha(X)(H)\underset{S^{n-1}\times\map(\cU H,\cU G)}{\sqcup}D^n\times\map(\cU H,\cU G),\nonumber
\end{alignat}
where the second equality holds by Lemmas \ref{Lcm} and the third equality holds by Lemma \ref{orpush}.
The result then follows since $\Orb(H,G)\to \map(\cU H,\cU G)$ is a weak equivalence, and weak equivalences are preserved by pushouts along Hurewicz cofibrations.

We must also show that, assuming $\cof_\alpha(X)\to RL\cof_\alpha(X)$
is a weak equivalence for all $\alpha$ smaller than a given limit ordinal $\beta$,
then $\cof_\beta(X)\to RL\cof_\beta(X)$ is also a weak equivalence.
By Lemma \ref{orpush} and the computation (\ref{azu}), the inclusion $RL\cof_\alpha(X)(H) \hookrightarrow RL\cof_{\alpha+1}(X)(H)$ is a Hurewicz cofibration.
Clearly, $\cof_\alpha(X)(H) \hookrightarrow \cof_{\alpha+1}(X)(H)$ is also a Hurewicz cofibration.
Weak equivalences being preserved by (transfinite) colimits along Hurewicz cofibrations, it follows that the map
\[
\cof_\beta(X)(H) =
\underrightarrow{\lim}\,
\cof_\alpha(X)(H)
\,\,\rightarrow\,\,
\underrightarrow{\lim}\,
RL\cof_\alpha(X)(H)=
\underrightarrow{\lim}\,
\map(\cU H,L\cof_\alpha(X))
\]
is a weak equivalence.
By Lemma \ref{orfilt}, we can then rewrite the right hand side as
\[
\map\big(\cU H,
\underrightarrow{\lim}\,
L\cof_\alpha(X)\big)
\cong \map\big(\cU H,L\,
\underrightarrow{\lim}\,
\cof_\alpha(X)\big)\cong RL \cof_\beta(X)(H),
\]
yielding the desired weak equivalence.
\qed

\begin{proposition}\label{du}\dontshow{du}
The canonical map $\cof X \to \RR\LL X$ is a weak equivalence of $\Orb$-spaces.
\end{proposition}

\proof
The map $\cof X \to \RR\LL X$ can be factored as
\[
\cof X\to RL\cof X\to R\fib L \cof X=\RR\LL X.
\]
The first map is a weak equivalence by Lemma \ref{isco}, while the second map is a homotopy equivalence by Lemma \ref{heLX}.
\qed

\subsection{Homotopy groups {\rm\dontshow{sec:HG}}}\label{sec:HG}

Given a groupoid $\cG$ and a point $q\in\tmap(\cB G,\cG)$, we can define the homotopy groups $\pi_n^G(\cG,q)$ of $\cG$, based at the point $q$,
as the set of path components in the space of maps $S^n\times\fib\cB G\to\fib\cG$
whose restriction to $\fib\cB G$ is equal to $\fib (q)$.
That is, $\pi_n^G(\cG,q)$ is the set of path components of the fiber over $\fib (q)$ of the map
$$
\map(S^n\times\fib\cB G,\fib\cG)\longrightarrow\map(\fib\cB G,\fib\cG).
$$
For $n=0$, it is preferable to use the unbased version
$\pi_0^G(\cG):=\pi_0\map(\fib\cB G,\fib\cG)$.

\begin{definition}
A map of topological groupoids $\cH\to\cG$ is a {\em weak equivalence} if, for all allowed isotropy groups $G\in\ccF$,
the map $\pi_0^G(\cH)\to \pi_0^G(\cG)$ is an isomorphism,
and, for all $q\in\tmap(\cB G,\cH)$ and positive integers $n$, the map
$\pi_n^G(\cH,q)\to\pi_n^G(\cG,f\circ q)$ is an isomorphism.
\end{definition}

\begin{remark}
Unraveling the definition, a map of groupoids $\cH\to\cG$ is a weak equivalence iff the induced map of $\Orb$-spaces $\RR\cH\to\RR\cG$ is a weak equivalence.
\end{remark}

We have defined the homotopy groups of $\cG$ to be those of the spaces $\RR\GG(G)=\map(\cU G,\fib\cG)$.
However, there's nothing special about the choice of $\cU G$.
If $P\to X$ is another universal $G$-bundle with paracompact base, then $\cO:=\gauge(P)$ and $\cU G=\gauge(\|\cE G\|)$ are homotopy equivalent groupoids.
The homotopy groups of $\cG$ therefore also agree with those of the space $\map(\cO,\fib\cG)$.
If $\cG$ is in addition fibrant (and paracompact), then the groupoids $\cG$ and $\fib(\cG)$ are homotopy equivalent by Lemma \ref{hece}.
Omitting the base points from the notation, we can then write \dontshow{lqg}
\begin{equation}\label{lqg}
\pi_*^G(\cG)\,\cong\,\pi_*\map(\cO,\cG).
\end{equation}

\begin{proposition}[Whitehead theorem]
\label{whT}\dontshow{whT}
Let $\cH$, $\cG$ be topological groupoids which are both cellular and fibrant, and let $f:\cH\to\cG$ be a weak equivalence.
Then $f$ is a homotopy equivalence.
\end{proposition}

\proof
Write $\cG=\colim_\alpha\cG_\alpha$, where \dontshow{PKz}
\begin{equation}\label{PKz}
\cG_{\alpha+1}=\cG_\alpha\underset{S^{n-1}\times\cO}{\sqcup}D^n\times\cO
\qquad\text{and}\qquad
\cG_\beta=\underset{\alpha<\beta}{\colim}\,\cG_\alpha
\end{equation}
for limit ordinals $\beta$.
Since $\cG$ is paracompact and fibrant,
the orbits $\cO$ used in (\ref{PKz}) are necessarily also paracompact and fibrant,
and so we can use them to compute $\pi_*^G$ as in (\ref{lqg}).

By induction on $\alpha$, we construct maps $g_\alpha:\cG_\alpha\to\cH$,
and homotopies $h_\alpha:\cG_\alpha\times[0,1]\to\cG$ between $f\circ g_\alpha$ and $i_\alpha:\cG_\alpha\hookrightarrow\cG$.
We also make sure
that the $g_\alpha$ and $h_\alpha$ are compatible with those for smaller $\alpha$.

Suppose we have $g_\alpha$ and $h_\alpha$ for a given $\alpha$,
and let's assume for simplicity that the $n$ in (\ref{PKz}) is at least $2$ (the cases $n=0,1$ are left to the reader).
The attaching map $S^{n-1}\times\cO\to\cG_\alpha$ represents zero in $\pi_{n-1}^G(\cG)$.
Its image under $g_\alpha$ is then also zero in $\pi_{n-1}^G(\cH)$, so we can extend $g_\alpha$ to a map $\tilde g_{\alpha+1}:\cG_{\alpha+1}\to\cH$.
The union
\[
(f\circ\tilde g_{\alpha+1})\cup h_\alpha\cup i_{\alpha+1}\,:\,
(D^n\times \cO)\cup(S^{n-1}\times[0,1]\times\cO)\cup (D^n\times \cO)\,\to\,\cG
\]
now defines an element $x\in\pi_n^G(\cG)$, which is exactly the obstruction to extending $h_\alpha$ to a homotopy $h_{\alpha+1}$
between $f\circ \tilde g_{\alpha+1}$ and $i_{\alpha+1}$.
So let $g_{\alpha+1}:\cG_{\alpha+1}\to\cH$ be another extension of $g_\alpha$, obtained by subtracting $(f_*)^{-1}(x)$ from $\tilde g_{\alpha}$.
The obstruction then vanishes for $g_{\alpha+1}$, and we can extend $h_\alpha$ to a homotopy $h_{\alpha+1}$
between $f\circ g_{\alpha+1}$ and $i_{\alpha+1}$.

Now suppose that we have constructed $g_\alpha$ and $h_\alpha$ for all ordinals $\alpha$ smaller than a given limit ordinal $\beta$,
and that they are all compatible with each other.
Then we jut set $g_\beta:=\cup\,g_\alpha$, $h_\beta:=\cup\,h_\alpha$.

The map $g$ constructed above is a right homotopy inverse of $f$.
Applying the same argument to $g$, we obtain a right homotopy inverse of $g$.
Hence $g$ is a two sided homotopy inverse of $f$.
\qed
\vspace{.3cm}

Unfortunately, there is no Quillen model structure on the category of topological groupoids
such that the weak equivalences are the ones described above, the cofibrant objects are the retracts of cellular groupoids,
and the fibrant objects are those introduced in Section \ref{se:Fr}.
Nevertheless, with this notion of weak equivalence, together with the classes of cellular and fibrant objects,
the category of topological groupoids behaves much like a topological model category.

\begin{proposition}\label{miwe}\dontshow{miwe}
Let $\cG\to\cG'$ be a weak equivalence between fibrant topological groupoids.
Then for any cellular groupoid $\cH$, the induced map \dontshow{ohz}
\begin{equation}\label{ohz}
\map(\cH,\cG)\to\map(\cH,\cG')
\end{equation}
is a weak equivalence of spaces.
\end{proposition}

\proof
Let $\ccC$ be the class of topological groupoids $\cH$ such that for every weak equivalence $f:\cG\to\cG'$ between fibrant groupoids,
the induced map (\ref{ohz}) is a weak equivalence.
We shall prove the result in four steps, denoted $a)$, $b)$, $c)$, and $d)$.

In step $a)$,
we show that for any $G\in\ccF$ and any principal $G$-bundle $P$ over a paracompact base, we have
$\gauge(P)\in\ccC$.
We then show in $b)$ that if $\cH\in\ccC$ then $X\times \cH\in\ccC$, for any CW-complex $X$.
This proves in particular that $S^{n-1}\times \gauge(P)$ and $D^n\times \gauge(P)$ are in $\ccC$.
In $c)$, we assume that
$\cH$ is in $\ccC$ and show that $\cH\sqcup_{S^{n-1}\times\gauge(P)}D^n\times\gauge(P)$ is also in $\ccC$.
Finally, for our last step $d)$, we show that
if $\cH=\colim_{\alpha<\beta}\cH_\alpha$ is a cellular groupoid and $\cH_\alpha\in\ccC$, for all $\alpha<\beta$, then $\cH_\beta\in\ccC$.

$a)$
For any fibrant groupoid $\cG$ we have weak equivalences
\[
\map\big(\gauge(P),\cG\big)
\approx\map\big(\M_{\gauge(P)},\M_\cG\big)
\approx\map\big(\M_{\cU G},\M_{\fib\cG}\big)
\approx\map\big(\cU G,\fib\cG\big),
\]
where the first and third equivalences hold by Theorem \ref{wty}, and the second equivalence holds by Propositions \ref{ijd} and \ref{BGstr}.
But the functor $\cG\mapsto \map\big(\cU G,\fib\cG\big)$ preserves weak equivalences, so $\gauge(P)\in\ccC$.

$b)$
The functor $\map(X\times\cH,-)\cong\map(X,\map(\cH,-))$ preserves weak equivalences (when restricted to fibrant groupoids)
because the functors $\map(\cH,-)$ and $\map(X,-)$ both do.

$c)$
The mapping space $\map(\cH\sqcup_{S^{n-1}\times\gauge(P)}D^n\times\gauge(P),\cG)$
can be rewritten as a pullback
\begin{equation*}
\map(\cH,\cG)\underset{\map(S^{n-1}\times\gauge(P),\cG)}{\times}\map\big(D^n\times\gauge(P),\cG\big).
\end{equation*}
The restriction $\map(D^n\times\gauge(P),\cG)\to\map(S^{n-1}\times\gauge(P),\cG)$ is a Hurewicz fibration.
The functors $\map(\cH,-)$, $\map(S^{n-1}\times\gauge(P),-)$, and $\map(D^n\times\gauge(P),-)$
send weak equivalences between fibrant groupoids to weak equivalences of spaces.
And since pullbacks along fibrations preserve weak equivalences, so does the functor
\[
\map\Big(\cH\underset{S^{n-1}\times\gauge(P)}{\sqcup}D^n\times\gauge(P),-\Big).
\]

$d)$
Each restriction $\map(\cH_{\alpha+1},\cG)\to\map(\cH_\alpha,\cG)$ in the limit diagram
\[
\map\big(\underset{\alpha}{\colim}\,\cH_\alpha,\cG\big)\cong \underset{\alpha}{\lim}\map(\cH_\alpha,\cG)
\]
is a Hurewicz fibration.
The result then follows since each functor $\map(\cH_\alpha,-)$ sends weak equivalences between fibrant groupoids to weak equivalences of spaces,
and limits along fibrations preserve weak equivalences.
\qed

\begin{proposition}\label{mtf}\dontshow{mtf}
Let $\cH\to\cH'$ be a weak equivalence of cellular topological groupoids.
Then for any fibrant groupoid $\cG$, the induced map \dontshow{roz}
\begin{equation}\label{roz}
\map(\cH',\cG)\to\map(\cH,\cG)
\end{equation}
is a homotopy equivalence.
\end{proposition}

\proof
Let $\cG$ be a fibrant groupoid.
By Propositions \ref{fibfactor} and \ref{hece}, the fibrant replacement $\eta$ induces homotopy equivalences
\[
\begin{split}
\map(\cH,\cG)&\stackrel{\scriptscriptstyle\sim}{\rightarrow}
\map(\cH,\fib\cG)\stackrel{\scriptscriptstyle\sim}{\leftarrow}
\map(\fib\cH,\fib\cG),\\
\map(\cH',\cG)&\stackrel{\scriptscriptstyle\sim}{\rightarrow}
\map(\cH',\fib\cG)\stackrel{\scriptscriptstyle\sim}{\leftarrow}
\map(\fib\cH',\fib\cG).
\end{split}
\]
The result follows since by Proposition \ref{whT} and Corollary \ref{fcic}, the map $\fib(\cH)\to\fib(\cH')$ is a homotopy equivalence.
\phantom{.}\qed

\begin{proposition}\label{dc}
For each topological groupoid $\cG$, the map $\LL\RR\cG\to\fib\cG$ is a weak equivalence.
\end{proposition}

\proof
The map $\LL\RR\cG\to\fib\cG$ is a weak equivalence of groupoids if and only if $\RR\LL\RR \cG\to \RR\fib \cG$ is a weak
equivalence of $\Orb$-spaces.
The result then follows from the following commutative diagram.
\[
\xymatrix{
&R\fib\ar@{=}[rr]\ar[dr]&&R\fib\ar[dr]^{\simeq}_{\text{(Corollary \ref{fibsq})\,}}\\
\cof R\fib\ar[ur]^{\simeq}\ar[dr]^{\simeq}_{\text{(Lemma \ref{isco})}}&&RLR\fib\ar[ur]\ar[dr]&&R\fib\fib\\
&RL\cof R\fib\ar[ur]\ar[dr]^{\simeq}_{\text{(Lemma \ref{heLX})}}&&R\fib L R\fib\ar[ur]\\
&&R\fib L\cof R\fib\ar[ur]
}
\]
\qed

\subsection{The main theorems \rm\dontshow{sec:mt}}\label{sec:mt}

We are finally in a position to prove our main theorems.
The idea is that $\Orb$-spaces model the homotopy theory of fibrant cellular topological groupoids, which in turn model the homotopy theory of cellular topological stacks, or {\em orbispaces}.
Moreover, $\Orb$-spaces enjoy all of the excellent formal properties of a diagram category,
making them the most useful model of the homotopy theory of orbispaces for many applications.
Finally, the theory of $\Orb$-spaces is quite flexible, for as input we specify an arbitrary family of orbit types
(the members of which must satisfy a mild paracompactness condition and whose isomorphism classes must comprise a proper set)
and choose whether to allow arbitrary maps or restrict to the representable ones.

\begin{theorem}\label{mt1}\dontshow{mt1}
For any $\Orb$-space $X$ the derived unit map \dontshow{unit}
\begin{equation}\label{unit}
\cof X\to\RR\LL X
\end{equation}
is a weak equivalence of $\Orb$-spaces.
For any topological groupoid $\cG$ the derived counit map \dontshow{counit}
\begin{equation}\label{counit}
\LL\RR\cG\to\fib\cG
\end{equation}
is a weak equivalence of topological groupoids.
\end{theorem}

\proof
This is just Propositions \ref{du} and \ref{dc}, respectively.
\qed

\begin{theorem}\label{mt2}\dontshow{mt2}
For any pair of topological groupoids $\cG$ and $\cH$, with $\cH$ cellular, the derived counit (\ref{counit}) induces a weak equivalence of mapping spaces
\begin{equation*}
\map\big(\cH,\fib\cG\big)\,\approx\,\map\big(\!\cof\RR\cH,\RR\cG\big).
\end{equation*}
For any pair of $\Orb$-spaces $X$ and $Y$ the derived unit (\ref{unit}) induces a weak equivalence of mapping spaces
\begin{equation*}
\map\big(\!\cof Y,X\big)\,\approx\,\map\big(\LL Y,\fib\LL X\big).
\end{equation*}
\end{theorem}

\proof
By Theorem \ref{mt1}, the derived
unit $\cof X\to\RR\LL X$
and
counit $\LL\RR\cH\to\fib\cH$
are weak equivalences.
The maps \dontshow{hah}
\begin{equation}\label{hah}
\map(\cH,\fib\cG)\leftarrow\map(\fib\cH,\fib\cG)\to\map(\LL\RR\cH,\fib\cG)\simeq\map(\cof\RR\cH,\RR\cG),
\end{equation}
are weak equivalences by Propositions \ref{fibfactor} and \ref{mtf} respectively.
The maps
\[
\map(\cof Y,X)\leftarrow\map(\cof Y,\cof X)\to\map(\cof Y,\RR\LL X)\simeq\map(\LL Y,\fib\LL X).
\]
are weak equivalences because $\Orb$-spaces form a topological model category.
\qed
\vspace{.3cm}

In fact the maps in (\ref{hah}) are homotopy equivalences; we stated the result as we did for esthetic reasons.
Also, as all $\Orb$-spaces are fibrant, we didn't need to ask that $X$ be fibrant.

\begin{theorem}\label{mt3}\dontshow{mt3}
Let $\cH$ be a cellular groupoid, and $\cG$ be a fibrant groupoid.
Then we have a weak equivalence
\[
\map(\M_\cH,\M_\cG)\approx\map(\cH,\cG).
\]
\end{theorem}

\proof
This follows from Theorem \ref{wty} and Lemma \ref{pcrap}.
\qed
\vspace{.3cm}

The above theorems show, in particular, that the homotopy categories of cellular topological stacks,
cellular fibrant topological groupoids and cofibrant $\Orb$-spaces are equivalent.

\subsection{Applications to equivariant homotopy theory \rm\dontshow{sec:gq}}\label{sec:gq}

In this section we restrict to case \casetwo.
This means, among other things, that the attaching maps used for constructing the cellular stacks are assumed to be faithful.

Let $G$ be a topological group and let $\ccF$ be a family of closed subgroups of $G$.
If $\ccF$ is closed under taking subgroups, then Elmendorf's Theorem \ref{et's} produces a space $E_\ccF:=\Psi($terminal $\Orb_G$-space$)$
with the property that  $(E_\ccF)^H$ is contractible for all $H\in\ccF$ and empty otherwise (see \cite{Lue05} for more applications of $E_\ccF G$).
Let $B_\ccF G:=[E_\ccF G/G]$ be the quotient stack.
The homotopy theory of cellular stacks over $B_\ccF G$
is then equivalent to the homotopy theory of cellular $G$-spaces with stabilizers in $\ccF$,
via functors
%\[
%\begin{split}
%\big\{\text{Stacks over $B_\ccF G$}\big\}\put(7,3){$\longrightarrow$}\put(7,-2){$\longleftarrow$}\hspace{1cm}&\big\{\text{$G$-spaces with stabilizers in $\ccF$}\big\}\\
%\big(\X\to B_\ccF G\big)\hspace{.4cm}\hspace{.3cm}\mapsto\hspace{.7cm}&\big(\X\times_{B_\ccF G}E_\ccF G\acted G\big)\\
%\big([M\times E_\ccF G/G]\to[E_\ccF G/G]\big)
%\hspace{.5cm}\raisebox{1.8mm}{\rotatebox{180}{$\mapsto$}}\hspace{.6cm}
%&\big(M\acted G\big)
%\end{split}
%\]
\begin{eqnarray*}
\big\{\text{Stacks over $B_\ccF G$}\big\}&\put(-8,3){$\longrightarrow$}\put(-8,-2){$\longleftarrow$}&
%\put(7,3){$\longrightarrow$}\put(7,-2){$\longleftarrow$}\hspace{1cm}&
\big\{\text{$G$-spaces with stabilizers in $\ccF$}\big\}\\
\big(\X\to B_\ccF G\big)&\mapsto&\big(\X\times_{B_\ccF G}E_\ccF G\acted G\big)\\
\big([M\times E_\ccF G/G]\to[E_\ccF G/G]\big)
%\hspace{.5cm}
&\raisebox{1.8mm}{\rotatebox{180}{$\mapsto$}}&
%\hspace{.6cm}
\big(M\acted G\big)
\end{eqnarray*}
So given a model for the homotopy theory of orbispaces, we get a model for the homotopy theory of $G$-spaces by taking a slice category.
The more surprising fact is that one can also go the other way.

%Given topological groups $H$, $K$, recall that we are in case \casetwo, and that $\map(H,K)$ denotes the space of group homomorphisms which are closed
%inclusions.
\begin{definition}\label{Dvc}\dontshow{Dvc}
Let $G$ be a topological group and $\ccF$ a family of closed subgroups of $G$ (not necessarily closed under taking subgroups) such that the quotient projections $G\to G/H$ are locally trivial bundles for all $H\in\ccF$.
Then the group $G$ is {\em $\ccF$-contractible}
if for every pair of subgroups $H, K \in \ccF$, the map \dontshow{trA}
\begin{equation}\label{trA}
\big\{g\in G\,\big|\,H^g<K\big\}\rightarrow \map(H,K)
\end{equation}
is a trivial Serre fibration. Here $H^g$ denotes the conjugate subgroup $g^{-1}Hg$.
%To avoid pathologies, we also require that the maps $G\to G/K$, $K\in\ccF$ be locally trivial bundles.
\end{definition}

Our terminology is justified by noting that if $\{1\}\in\ccF$, then $\ccF$-contractible implies weakly contractible.
%If the elements of $\ccF$ are compact Lie groups, then Definition \ref{Dvc} simplifies to the following:

\begin{lemma}
If the subgroups comprising $\ccF$ are compact Lie groups, then $G$ is $\ccF$-contractible if and only if for every monomorphism $f:H\to K$, the space \dontshow{pfh}
\begin{equation}\label{pfh}
\big\{g\in G\,\big|\,\mathrm{Ad}(g)|_H=f\big\}
\end{equation}
is weakly contractible.
\end{lemma}

\proof
Clearly, if $G$ is $\ccF$-contractible, then the spaces (\ref{pfh}) are weakly contractible.
The converse holds because the map (\ref{trA}) is $K$-equivariant and the space $\map(H,K)$ is a disjoint union of $K$-orbits.

\qed

\begin{example}\label{xuh}\dontshow{xuh}
Let $\mathcal H$ be an infinite dimensional Hilbert space and let $U(\mathcal H)$ be its group of unitary automorphisms under the norm topology.
Let $\ccF$ be the collection of subgroups $H<U(\mathcal H)$ which are compact Lie
and such that each irreducible representation of $H$ appears infinitely many times in $\mathcal H$.
Then $U(\mathcal H)$ is $\ccF$-contractible.

Indeed, given a monomorphism $H\to K$ between elements of $\ccF$, the inclusion $H\hookrightarrow U(\mathcal H)$ and the composite
$H\to K\hookrightarrow U(\mathcal H)$ are equivalent representations of $H$.
It follows that (\ref{pfh}) is non-empty.
The space (\ref{pfh}) has a simply transitive action of the centralizer of $H$ in $U(\mathcal H)$, so it's enough to show that this centralizer is contractible.
Since each irreducible representation of $H$ appears infinitely many times in $\mathcal H$, we may write $\mathcal H$ as
\begin{equation*}
\mathcal H\simeq\bigoplus_{\rho}\mathcal K_\rho\otimes \rho
\end{equation*}
where $\mathcal K_\rho$ are infinite dimensional Hilbert spaces, and the direct sum is indexed over all the irreducible representations $\rho$ of $H$.
It follows from Schur's lemma that the centralizer of $H$ in $U(\mathcal H)$ is given
by \dontshow{zUf}
\begin{equation}\label{zUf}
Z_{U(\mathcal H)}(H)\cong \prod_\rho U(\mathcal K_\rho).
\end{equation}
One might be tempted to argue that by Kuiper's theorem \cite{Kui65}, each factor $U(\mathcal K_\rho)$ is contractible, and thus so is (\ref{zUf}).
However, the topology on (\ref{zUf}) is not the product topology,
rather it's the topology induced from the ``max'' metric on the product of the metric spaces $U(\mathcal K_\rho)$.
The space (\ref{zUf}) is nevertheless contractible,
as can be seen by applying $H$-fixed points to the $H$-equivariant version of Kuiper's theorem \cite[Proposition A3.1]{AS04}.
\end{example}

We now go back to the setup of Section \ref{se:ctg}, so that $\ccF$ is a class of arbitrary topological groups and not just a set of subgroups of a fixed topological group.
The following is the main result of this section.

\begin{theorem}\label{Tvc}\dontshow{Tvc}
Let $\ccF$ be as in Section \ref{se:ctg},
and let $G$ be a topological group with a family $\ccF'$ of closed subgroups satisfying $\ccF'\subset\ccF$.
Suppose that for each group $H$ in $\ccF$ there exists a group $H'$ in $\ccF'$ which is isomorphic to $H$.
Let $\Orb_G$ be the category of $G$-orbits with isotropy in $\ccF'$, and let $\Orb$ be the topological category introduced in Section \ref{catorb}.
Then if $G$ is $\ccF'$-contractible, the categories $\Orb_G$ and $\Orb$ are weakly equivalent (in the sense of Definition \ref{da4}).
\end{theorem}

\proof
Let us fix a particular model for $\Orb$, as in (\ref{hqv}); that is,
\[
\Orb(H,K)=\map(H,K)\times_K EK.
\]
Replacing $\Orb_G$ by an equivalent category, we may identify its objects with $\ccF'$.
The morphism spaces are then given by
\[
\Orb_G(H,K)=\big\{g\in G\,\big|\,H^g\subset K\big\}/K.
\]
Let us also introduce an auxiliary category $\mathcal O$ on the same objects $\ccF'$, with morphism spaces
\[
\mathcal O(H,K):=\big\{g\in G\,\big|\,H^g\subset K\big\}\times_K EK
\]
and composition as in $\Orb$.
We then have two functors $\Orb_G\leftarrow \mathcal O\rightarrow \Orb$.
The first one is a weak equivalence because the projection
$\{g\in G\,\big|\,H^g\subset K\}\times_K EK\to \{g\in G\,\big|\,H^g\subset K\}/ EK$
is a locally trivial bundle with contractible fibers.
The second one is a weak equivalence because
$\{g\in G\,\big|\,H^g\subset K\}\times_K EK\to \map(H,K)\times_K EK$
is a trivial Serre fibration.
\qed
\vspace{.3cm}

Combining Theorem \ref{Tvc} with our main theorems, Elmendorf's theorem, and Lemma \ref{Aee}, we obtain the following result.

\begin{corollary}
Let $\ccF$, $G$, and $\ccF'$ be as in Theorem \ref{Tvc}.
Then the category of orbispaces with isotropy groups in $\ccF$ has the same homotopy theory as the category of $G$-spaces with isotropy groups in $\ccF'$.
\end{corollary}

For the above results to be of any use, we need to show that there are enough contractible groups.
We accomplish this by a variation on the small object argument.

\begin{proposition} \label{rb:} \dontshow{rb:}
Let $G_0$ be a topological group and $\ccF_0$ a family of closed subgroups.
Then there exists a group $G$, containing $G_0$ as a closed subgroup, which is $\{H^g| g\in G, H\in \ccF_0\}$-contractible.
Moreover, if $G_0$ and the quotients $G_0/H$, $H\in\ccF_0$, are cellular as spaces, then $G$ may be taken to be cellular as well.
\end{proposition}

\proof
We let $G$ be the colimit of the groups $G_i$, $i\in\mathbb N$, inductively defined as follows.
Let $\mathcal{S}$ be the set of all diagrams \dontshow{srg}
\begin{equation}\label{srg}
\begin{matrix}\xymatrix{
S^{n-1}\ar[r]^(.3)\alpha\ar[d]&\big\{g\in G_i\,\big|\,H^g<K\big\}\ar[d]\\
D^n\ar[r]^(.4)\beta&\map(H,K),
}\end{matrix}
\end{equation}
where $H,K<G_i$ are in $\ccF_i:=\{H^g| g\in G_i, H\in \ccF_0\}$.
For each $(\alpha,\beta)\in \mathcal{S}$, let $G_i(\alpha,\beta)$ be the free topological group on $G_i\cup_\alpha D^n$ subject to the relations
\[
\begin{split}
g\cdot g'&=g g'\hspace{1.7cm}\forall g,g'\in G_i\,,\\
x^{-1}\cdot h \cdot x &= \beta(x)(h)\qquad \forall x\in D^n, h\in H.
\end{split}
\]
The group $G_i(\alpha,\beta)$ is an HNN-extension, from which it follows that $G_i\to G_i(\alpha,\beta)$ is a closed inclusion.
More precisely, an element $y\in G_i(\alpha,\beta)$ can be written as \dontshow{hnc}
\begin{equation}\label{hnc}
\qquad \qquad y=g_0x_{11}^{\pm 1}x_{12}^{\pm 1}\cdots x_{1n_1}^{\pm 1}g_1x_{21}^{\pm 1}\ldots \ldots g_{m-1}x_{m1}^{\pm 1}\cdots x_{mn_m}^{\pm 1}g_m,\qquad \quad {\scriptstyle (n_i>0)},
\end{equation}
with $g_i\in G_i$ and $x_{ij}\in D^n$.
Assuming that $g_1,\ldots,g_{m-1}\not\in H$, $x_{ij}\not\in S^{n-1}$, and that (\ref{hnc}) is a reduced word,
this expression is then unique up to an action of $H^m$.
It follows that $G_i(\alpha,\beta)$ is built from $G_i$ by successively attaching spaces of the form
$(G_i^{m+1}\times (D^n)^k)/H^m$ along unions of subspaces of the form
$(G_i^m\times H\times (D^n)^k)/H^m$, $(G_i^{m+1}\times (D^n)^{k-1}\times S^{n-1})/H^m$, and $(G_i^{m+1}\times (D^n)^{k-1})/H^m$.
A further analysis shows that if $G_i$ and $G_i/K$, $K\in\ccF$, are cellular, then the same holds for $G_i(\alpha,\beta)$ and $G_i(\alpha,\beta)/K$.

Let $G_{i+1}$ be the amalgamated product over $G_i$ of all the $G_i(\alpha,\beta)$.
Amalgamated products also have normal forms, so the map $G_i\to G_{i+1}$ is a closed inclusion.
Moreover, if all the $G_i(\alpha,\beta)$ and $G_i(\alpha,\beta)/H$ are cellular, then the same holds for $G_{i+1}$ and $G_{i+1}/H$.

To see that $G:=\colim G_i$ is $(\cup_i\ccF_i)$-contractible, consider groups $H,K \in \cup_i\ccF_i$
and a lifting problem
\begin{equation}\label{SRG}
\begin{matrix}\xymatrix{
S^{n-1}\ar[r]^(.32)\alpha\ar[d]&\big\{g\in G\,\big|\,H^g<K\big\}\ar[d]\\
D^n\ar[r]^(.4)\beta\ar@{.>}[ur]^(.45){\exists ?}&\map(H,K).
}\end{matrix}
\end{equation}
Since $S^{n-1}$ and $D^n$ are compact, the diagram (\ref{SRG}) factors at some finite stage of the colimit and we are therefore reduced to a diagram like (\ref{srg}).
We then have a lift $D^n\to \{g\in G_i(\alpha,\beta)\,\big|\,H^g<K\}$ by definition of $G_i(\alpha,\beta)$.
\qed

\begin{corollary}\label{xtg}\dontshow{xtg}
For any family of topological groups $\ccF$,
there exists a topological group $G$ and a family of closed subgroups $\ccF'$ satisfying the hypotheses of Theorem \ref{Tvc}.
Moreover, if all the groups in $\ccF$ are cellular (as spaces), we may choose $G$ to be cellular as well.
\end{corollary}

\proof
Let $S$ be a set of representatives of isomorphism classes of elements of $\ccF$, and let
\[
G_0:=\underset{S'\subset S}{\colim}\,\big(\prod_{H\in S'}H\big),
\]
where the colimit ranges over all finite subsets of $S$.
Then apply Proposition \ref{rb:}.
\qed

\begin{remark}
In Corollary \ref{xtg}, it is easy to arrange for $\ccF'$ to be equal to the family of those subgroups of $G$ belonging to $\ccF$
(as opposed to being merely contained in it).
We have phrased Theorem \ref{Tvc} as we did for greater flexibility, and so that it applies to Example \ref{xuh}.
\end{remark}

Given an action of a group $G$ on a space $X$, let $[X/G]:=\M_{G\ltimes X}$ be the stack-theoretic quotient of $X$ by $G$.
Stacks of the form $[X/G]$ are called {\em global quotients}.
Recall that a cellular stack is a stack of the form $\M_\cG$, where $\cG$ is a cellular groupoid.
As a corollary of Theorem \ref{Tvc}, we see that every cellular orbispace is a global quotient.

\begin{theorem}
Let $\ccF$, $G$, and $\ccF'$ be as in Theorem \ref{Tvc}, and suppose that $G$ is $\ccF'$-contractible.
Then given any cellular (with respect to $\ccF$) stack $\Y$, there exists a $G$-space $X$ such that $\Y\simeq [X/G]$.
\end{theorem}

\proof
Combining our main theorems, Lemma \ref{Aee}, Theorem \ref{Tvc}, and Elmendorf's theorem,
we may find a $G$-space $X'$ such that $[X'/G]$ is homotopy equivalent to $\Y$.
In particular, we get a map $f:\Y\to [X'/G]$.
Let $\X$ be the pullback stack
\[
\xymatrix{
\X\ar[r]\ar[d]&\M_{X'}\ar[d]\\
\Y\ar[r]^(.4)f&[X'/G].
}
\]
The action of $G$ on $X'$ induces an action on $\X$ making it into a principal $G$-bundle over $\Y$.
Moreover, since we are in case \casetwo, the map $f$ is necessarily injective on stabilizer groups,
so the stack $\X$ is of the form $\M_X$ for some space $X$.
It follows that $\Y\simeq [X/G]$.
\qed

\appendix
\section{Appendix}

\subsection{Lax presheaves of groupoids {\rm\dontshow{A:1}}}\label{A:1}

Let $\C$ be a small category.
A lax presheaf of (non-topologically-enriched) groupoids $\X$ on $\C$ assigns to each object $T$ of $\C$ a groupoid $\X(T)$,
to each arrow $f:U\to T$ of $\C$ a $1$-morphism $\X(f):\X(T)\to\X(U)$,
and to each composable pair $g:V\to U,f:U\to T$ of arrows of $\C$ a $2$-morphism \dontshow{2F1}
\begin{equation}\label{2F1}
\begin{matrix}
\xymatrix{
&\X(U)\ar[rd]^{\X(g)}\ar@2[d]|{\X(f,g)}&\\
\X(T)\ar[ur]^{\X(f)}\ar[rr]_{\X(fg)} & &\X(V)}
\end{matrix}
\end{equation}
such that for each triple of composable arrows $h:W\to V,g:V\to U,f:U\to T$ of $\C$ the pentagon
\begin{equation*}
\begin{matrix}
\xymatrix@C=6pt@R=.7cm{
&\big(\X(h)\X(g)\big)\X(f)\ar@{=}[rr]\ar@2[dl]_(.6){\X(g,h)\X(f)\,\,\,\,}&&\X(h)\big(\X(g)\X(f)\big)\ar@2[dr]^(.6){\,\,\,\X(h)\X(f,g)}\\
\X(gh)\X(f)\ar@2[drr]_{\X(f,gh)}&&&&\X(h)\X(fg)\ar@2[lld]^{\X(fg,h)}\\
&&\X(fgh)}
\end{matrix}
\end{equation*}
commutes.

A $1$-morphism $\varphi:\Y\to\X$ of lax presheaves of groupoids on $\C$ assigns
to each object $T$ of $\C$ a $1$-morphism $\varphi(T):\Y(T)\to\X(T)$ and to each arrow $f:U\to T$ of $\C$ a $2$-morphism
$$
\xymatrix@C=1.5cm{
\Y(T)\ar[r]^{\varphi(T)}\ar[d]_{\Y(f)}&\X(T)\ar[d]^{\X(f)}\\
\Y(U)\ar[r]_{\varphi(U)}\ar@2{->}[ur]|{\varphi(f)}&\X(U)}
$$
such that for each pair of composable arrows $g:V\to U,f:U\to T$ of $\C$ the pentagon
$$
\xymatrix@C=6pt@R=.7cm{
&\X(g)\varphi(U)\Y(f)\ar@2[rr]^{\X(g)\varphi(f)}&&\X(g)\X(f)\varphi(T)\ar@2[dr]^(.6){\,\,\X(f,g)\varphi(T)}\\
\varphi(V)\Y(g)\Y(f)\ar@2[drr]_{\varphi(V)\Y(f,g)}\ar@2[ur]^(.4){\varphi(g)\Y(f)\,\,}&&&&\X(fg)\varphi(T)\\
&&\varphi(V)\Y(fg)\ar@2[urr]_{\varphi(fg)}
}
$$
commutes.
The composite of a composable pair of $1$-morphisms $\psi:\Z\to\Y$ and $\varphi:\Y\to\X$ is the $1$-morphism $\varphi\psi:\Z\to\X$ defined by $(\varphi\psi)(T):=\varphi(T)\psi(T)$ and $(\varphi\psi)(f):=\varphi(f)\psi(f)$.
The fact that groupoids form a $2$-category implies that the requisite diagram commutes.

A $2$-morphism $\alpha:\psi\Rightarrow\varphi$ between a pair of $1$-morphisms $\varphi,\psi:\Y\to\X$ assigns to each object $T$ of $\C$ a $2$-morphism $\alpha(T):\psi(T)\Rightarrow\varphi(T)$ such that for each arrow $f:U\to T$ of $\C$ the square \dontshow{sQe}
\begin{equation}\label{sQe}
\begin{matrix}
\xymatrix@C=1.5cm{
\psi(U)\Y(f)\ar@2[r]^{\alpha(U)\Y(f)}\ar@2[d]_{\psi(f)}&\varphi(U)\Y(f)\ar@2[d]^{\varphi(f)}\\
\X(f)\psi(T)\ar@2[r]^{\X(f)\alpha(T)}&\X(f)\varphi(T)}
\end{matrix}
\end{equation}
commutes.
Horizontal and vertical composition of $2$-morphisms are defined in the obvious way and they satisfy the axioms of a 2-category.

\begin{remark}
Our definition of lax presheaf on $\C$ is formally the same as a contravariant weak $2$-functor from $\C$ to an arbitrary $2$-category $\D$,
provided the $2$-morphisms in $\D$ are invertible.
If there are non-invertible $2$-morphisms in $\D$ then we would require an extra condition to ensure that identity morphisms in $\D$ are taken to invertible $1$-morphisms.
\end{remark}

\subsection{$2$-categorical colimits \rm\dontshow{A:2}}\label{A:2}
Given 2-categories $\C$, $\D$ and lax contravariant functors $\X,\Y:\C\to\D$,
a {\em modification} $\alpha:\psi\Rrightarrow\varphi$ between a pair of natural transformations
$\psi,\varphi:\Y\Rightarrow \X$ consists of a family
$$
\big\{\alpha(T):\psi(T)\Rightarrow\varphi(T)\big\}_{T\in\C}
$$
of $2$-morphisms of $\D$ such that for each arrow $f:U\to T$ of $\C$, the square (\ref{sQe}) commutes\footnote{
This differs with the definition of modification given in \cite[section 7.3]{Bor94}
since we consider lax as opposed to strict functors.}.
A {\em cone} on a lax contravariant functor $\F:\C\to \D$ with vertex $V\in\D$
is a lax natural transformation $\F\Rightarrow\Delta_V$ from $\F$ to the constant functor $\Delta_V:\C\to\D$ at $V$.
The cones on $\F$ with vertex $V$ form the objects of a category $\cone(\F,V)$
whose arrows are the modifications.

\begin{definition}
A {\em lax colimit} of $\F:\C\to\D$ is a pair $(V,\pi)$ with $V$ an object of $\D$ and $\pi:\F\Rightarrow\Delta_V$ a cone on $\F$ with vertex $V$,
which is universal in the sense that for all $W\in \D$,
the functor \dontshow{upfc}
\begin{equation}\label{upfc}
\D(V,W)\to\cone(\F,W),
\end{equation}
induced by precomposition with $\pi$, is an equivalence of categories.
We denote a lax colimit of $\F$ by $\underset{T\in\C}{\hocolim}\,\F(T)$, or just $\hocolim \F$ when the context is clear.
\end{definition}

If the target 2-category $\D$ is the category of groupoids, we can give an explicit model for the lax colimit.
The objects of $\colim\F$ are the disjoint union of the objects of the groupoids $\F(T)$, and the morphisms are defined by generators and relations.
%There are two sets of generators and four sets of relations.
Writing $p^*(\sigma)$ in place of $\F(p)(\sigma)$, the generators consist of arrows
\vspace{-.3cm}

\begin{itemize}
\setlength{\itemindent}{6pt}
\setlength{\itemsep}{-3pt}
{\vspace{.4cm}\item[$({\scriptstyle\bullet\,\to\,\bullet})$]
\parbox{14.05cm}{%\vspace{.4cm}
$\alpha:\sigma\to\tau$ for each object $T\in\C$ and arrow $\alpha:\sigma\to \tau$ of $\F(T)$,}}
%we have a generator $\alpha$ between the corresponding objects $\sigma,\tau$ of $\hocolim\F$.}
{\vspace{.4cm}\item[$({\scriptstyle\bullet\,\mapsto\,\bullet})$]
\parbox{14.05cm}{%\vspace{.4cm}
$p_\sigma:\sigma\to p^*(\sigma)$ for each arrow $p:U\to T$ of $\C$ and object $\sigma\in\F(T)$,}}
%we have a generator $p_\sigma:\sigma\to p^*(\sigma)$, where $p^*(\sigma)$ is another notation for $\F(p)(\sigma)$.}
\end{itemize}
and are subject to the relations\vspace{-.1cm}

\begin{itemize}
\setlength{\itemindent}{10pt}
\setlength{\itemsep}{-1pt}
\item[$\Big(
\put(3,-2.5){\rotatebox{57}{$\scriptstyle\to$}}
\put(10,4){\rotatebox{-55}{$\scriptstyle\to$}}
\begin{matrix}
\scriptstyle\bullet\\
\scriptstyle\bullet\,\to\,\bullet
\end{matrix}\Big)$]
\parbox{13.95cm}{%\vspace{.4cm}
$\alpha\beta=\gamma$ for each pair of composable arrows $\alpha,\beta\in\F(T)$ with composite $\gamma\in\F(T)$,}
%the relation $\alpha\beta=\gamma$ should also hold in $\colim\F$.}
\item[$\Big(
\put(1,5){\rotatebox{-90}{$\scriptscriptstyle\mapsto$}}%\put(8,5){\rotatebox{-90}{$\scriptscriptstyle\mapsto$}}
\put(15.5,5){\rotatebox{-90}{$\scriptscriptstyle\mapsto$}}
\begin{matrix}
\scriptstyle\bullet\,\to\,\bullet\\
\scriptstyle\bullet\,\to\,\bullet
\end{matrix}\Big)$]
\parbox{13.95cm}{%\vspace{.4cm}
$p_\tau\alpha=p^*(\alpha)p_\sigma$ for each arrow $p:U\to T$ in $\C$, and arrow $\alpha:\sigma\to \tau$ in $\F(T)$,}
%one has the relation $p_\tau\alpha=p^*(\alpha)p_\sigma$.}
\item[$\Big(
\put(2.5,-2.2){\rotatebox{61.5}{$\scriptstyle\mapsto$}}
\put(9.3,4.7){\rotatebox{-44}{$\scriptstyle\mapsto$}}
\put(15.5,-2.5){$\scriptstyle\bullet$}
\begin{matrix}
\scriptstyle\bullet\hspace{.07cm}\\
\scriptstyle\bullet\,\mapsto\,\bullet
\end{matrix}
\Big)$]
\parbox{13.95cm}{%\vspace{.4cm}
$(pq)_\sigma\!=\!\F(p,q)(\sigma)\,q_{p^*\!(\sigma)}p_\sigma$ for each pair of composable arrows $p,q\in\C$ and object $\sigma\in\F(T)$,}
%we have $(pq)_\sigma=\F(p,q)(\sigma)\,q_{p^*\!(\sigma)}p_\sigma$,
%where $\F(p,q)$ is the 2-morphism (\ref{2F1}).}
\item[$\Big(
\put(1.5,6.5){\rotatebox{237}{$\scriptstyle\mapsto$}}
\put(8,4){\rotatebox{-90}{$\scriptscriptstyle\Mapsto$}}
\put(10,4.5){\rotatebox{-57}{$\scriptstyle\mapsto$}}
\begin{matrix}
\scriptstyle\bullet\\
\scriptstyle\bullet\,\to\,\bullet
\end{matrix}\Big)$]
\parbox{13.95cm}{%\vspace{.4cm}
$q_\sigma=\F(\zeta)(\sigma)\,p_\sigma$ for each 2-morphism $\zeta:p\Rightarrow q:U\to T$ in $\C$ and object $\sigma\in\F(T)$.}
%we have $q_\sigma=\F(\zeta)(\sigma)\,p_\sigma$.}
\end{itemize}
Here $\F(p,q)$ is the 2-morphism of (\ref{2F1}).
\noindent Note that if $\F$ is strict, as will be the case in our example,
the third relation simplifies to

\begin{itemize}
\setlength{\itemindent}{10pt}
\setlength{\itemsep}{-1pt}
\item[$\Big(
\put(3,-2.5){\rotatebox{57}{$\scriptstyle\mapsto$}}
\put(10,4){\rotatebox{-55}{$\scriptstyle\mapsto$}}
\begin{matrix}
\scriptstyle\bullet\\
\scriptstyle\bullet\,\mapsto\,\bullet
\end{matrix}
\Big)$]
\parbox{13.95cm}{%\vspace{.4cm}
%Given composable arrows $p,q$ of $\C$ and an object $\sigma$ of $\F(T)$, we have
$(pq)_\sigma=q_{p^*\!(\sigma)}p_\sigma$ for each pair of composable arrows $p,q\in\C$ and object $\sigma\in\F(T)$.
}
\end{itemize}
This particular model for $\hocolim\F$ comes with a preferred cone $(V,\pi)\in\cone(\F,\hocolim\F)$ satisfying a stricter universal property than (\ref{upfc}):
given a groupoid $\cG$, the induced map
\[
\Hom(\hocolim\F,\cG)\to\cone(\F,\cG)
\]
is not merely an equivalence but actually an isomorphism of groupoids.

We now come to the main purpose of this section, namely to verify that the groupoid $\M_\cG(\cH)$ of principal $\cG$-bundles on $\cH$ is equivalent
to the lax colimit over $\cov^2(\cH_0)$ of the groupoids $\Hom(\cH_U,\cG)$.
Given two groupoids $\cG$, $\cH$, and a cover $U\to\cH_0$ of $\cH_0$, recall from Section \ref{secPGB} (second half of the proof of Proposition \ref{zuo})
that a homomorphism $\sigma:\cH_U\to\cG$ induces a principal $\cG$-bundle $P(\sigma)$ on $\cH$.

\begin{proposition}\label{pbasm} \dontshow{pbasm}
Let $\cG$ and $\cH$ be topological groupoids and let $\cov^2(\cH_0)$ be the 2-category defined in (\ref{d2cv}).
Then the functor
$$
\Hom\big(\cH_{(-)},\cG\big):\cov(\cH_0)\to\{\text{\rm Groupoids}\}
$$
%a cover $U\to\cH_0$ to the groupoid $\Hom(\cH_U,G)$,
extends over the inclusion $\cov(\cH_0)\hookrightarrow\cov^2(\cH_0)$ and induces an equivalence
%such that map sending a groupoid homomorphism $\sigma:\cH_U\to\cG$ to the principal $\cG$-bundle $P(\sigma)$ then induces an equivalence \dontshow{keep}
\begin{equation}\label{keep}
\Phi\,:\,\,\underset{U\in\cov^2(\cH_0)}{\hocolim}\,
\Hom\big(\cH_U,\cG\big)
\,\stackrel{\scriptscriptstyle\sim}{\longrightarrow}\,\M_\cG(\cH)
\end{equation}
%associating to a groupoid homomorphism $\sigma:\cH_U\to\cG$ its principal $\cG$-bundle $P(\sigma)$
%is an equivalence of groupoids.
from its lax colimit to the groupoid of principal $\cG$-bundles on $\cH$.
\end{proposition}

\proof
Given morphisms $q_1, q_2:V\to U$ in $\cov(\cH_0)$, there's a unique 2-morphism between them in $\cov^2(\cH_0)$.
The extension of $\Hom(\cH_{(-)},\cG)$ to $\cov^2(\cH_0)$ assigns to that 2-morphism the natural transformation
$q_1^*\Rightarrow q_2^*:\Hom(\cH_U,\cG)\to\Hom(\cH_V,\cG)$ sending $\sigma\in\Hom(\cH_U,\cG)$ to the morphism
$v\mapsto \sigma(q_1(v),1,q_2(v))$ of $\Hom(\cH_V,\cG)$.

In the first half of the proof of Proposition \ref{zuo},
we have constructed for any given principal $\cG$-bundle $P$ on $\cH$,
an object of $\hocolim\Hom(\cH_{(-)},\cG)$ whose image under $\Phi$ is isomorphic to $P$.
It follows that $\Phi$ is essentially surjective.

Let $\M'_\cG(\cH)$ be the groupoid obtained by restricting $\M_\cG(\cH)$ along
$$
\Phi_0:(\hocolim\Hom(\cH_{(-)},\cG))_0\to (\M_\cG(\cH))_0.
$$
The objects of $\M'_\cG(\cH)$ are by definition those of $\hocolim\Hom(\cH_{(-)},\cG)$,
and the morphisms are given by
\[
\Hom_{\M'_\cG(\cH)}(\text{\,-\,},\text{\,-\,}):=\Hom_{\M_\cG(\cH)}\big(\Phi(\text{\,-\,}),\Phi(\text{\,-\,})\big).
\]
We then have a functor
\[
\Phi':\underset{U}{\hocolim}\,\Hom\big(\cH_U,\cG\big)\to\M'_\cG(\cH)
\]
given by the identity on objects and by $\Phi$ on morphisms.
Clearly, $\Phi$ is an equivalence if and only if $\Phi'$ is an isomorphism.
To show that $\Phi'$ is an isomorphism, we introduce a functor
\[
\Psi:\M'_\cG(\cH)\to\underset{U}{\hocolim}\,\Hom\big(\cH_U,\cG\big)
\]
and show that it's a strict inverse of $\Phi'$.

Given two objects $\sigma:\cH_U\to\cG$, $\tau:\cH_V\to\cG$ of $\M'_\cG(\cH)$, a morphism from $\sigma$ to $\tau$ is by definition a map
$f:P(\sigma)\to P(\tau)$ between the associated principal $\cG$-bundles.
Letting $W:=U\times_{\cH_0}V$, and $p:W\to U$, $q:W\to V$ be the projections,
the map $f$ induces a natural transformation $\alpha_f:p^*(\sigma)\Rightarrow q^*(\tau):\cH_W\to\cG$.
The functor $\Psi$ is then given by $\Psi(f):=q_\tau^{-1}\alpha_f\,p_\sigma$,
where $p_\sigma$, $q_\tau$ are the generators $({\scriptstyle\bullet\,\mapsto\,\bullet})$ corresponding to $p$ and $q$.
Graphically, this can be described as:
\[
\Psi:
\Big[f:\sigma\to\tau\Big]
\mapsto
\Big[\sigma\put(8,7){$\scriptstyle p_\sigma$}\longmapsto p^*(\sigma)\stackrel{\alpha_f}{\longrightarrow}
q^*(\tau)\stackrel{q_\tau}{\raisebox{4.3pt}{\rotatebox{180}{$\longmapsto$}}}\tau
\Big].
\]
The natural transformation $\alpha_f$ is made so that $\Phi(\alpha_f)=f$
(modulo the obvious identifications $P(\sigma)\simeq P(p^*\sigma)$, $P(\tau)\simeq P(q^*\tau)$).
It follows that $\Phi'\Psi=1$.

To see that $\Psi$ respects composition, consider three covers $U,V,W$ of $\cH_0$ and let
\[
\begin{matrix}
p_1:U\times_{\cH_0}V\hspace{.2mm}\to\hspace{.5mm}U
&q_1:V\times_{\cH_0}W\to V
&r_1:U\times_{\cH_0}W\to W\\
p_2:U\times_{\cH_0}W\to U
&q_2:U\times_{\cH_0}V\hspace{.5mm}\to\hspace{.5mm} V
&r_2:V\times_{\cH_0}W\to W\\
p_3:U\times_{\cH_0}V\times_{\cH_0}W\to U
&q_3:U\times_{\cH_0}V\times_{\cH_0}W\to V
&r_3:U\times_{\cH_0}V\times_{\cH_0}W\to W
\end{matrix}
\]
be the projections.
Let $\sigma:\cH_U\to\cG$, $\tau:\cH_V\to\cG$, $\nu:\cH_W\to\cG$ be objects,
and let $g:P(\sigma)\to P(\tau)$, $f:P(\tau)\to P(\nu)$ be morphisms in $\M'_\cG(\cH)$.
Then $\Psi(fg)=\Psi(f)\Psi(g)$ holds because the obvious diagram
{
\psfrag{1}{$\tau$}
\psfrag{2}{$q_2^*\tau$}
\psfrag{3}{$q_3^*\tau$}
\psfrag{4}{$q_1^*\tau$}
\psfrag{5}{$p_1^*\sigma$}
\psfrag{6}{$p_3^*\sigma$}
\psfrag{7}{$r_3^*\nu$}
\psfrag{8}{$r_2^*\nu$}
\psfrag{9}{$\sigma$}
\psfrag{10}{$p_2^*\sigma$}
\psfrag{11}{$r_1^*\nu$}
\psfrag{12}{$\nu$}
\psfrag{13}{\rotatebox{60}{$\scriptstyle\alpha_g$}}
\psfrag{14}{\rotatebox{-60}{$\scriptstyle\alpha_f$}}
\psfrag{15}{$\scriptstyle\alpha_{fg}$}
\[
\begin{matrix}
\epsfig{file=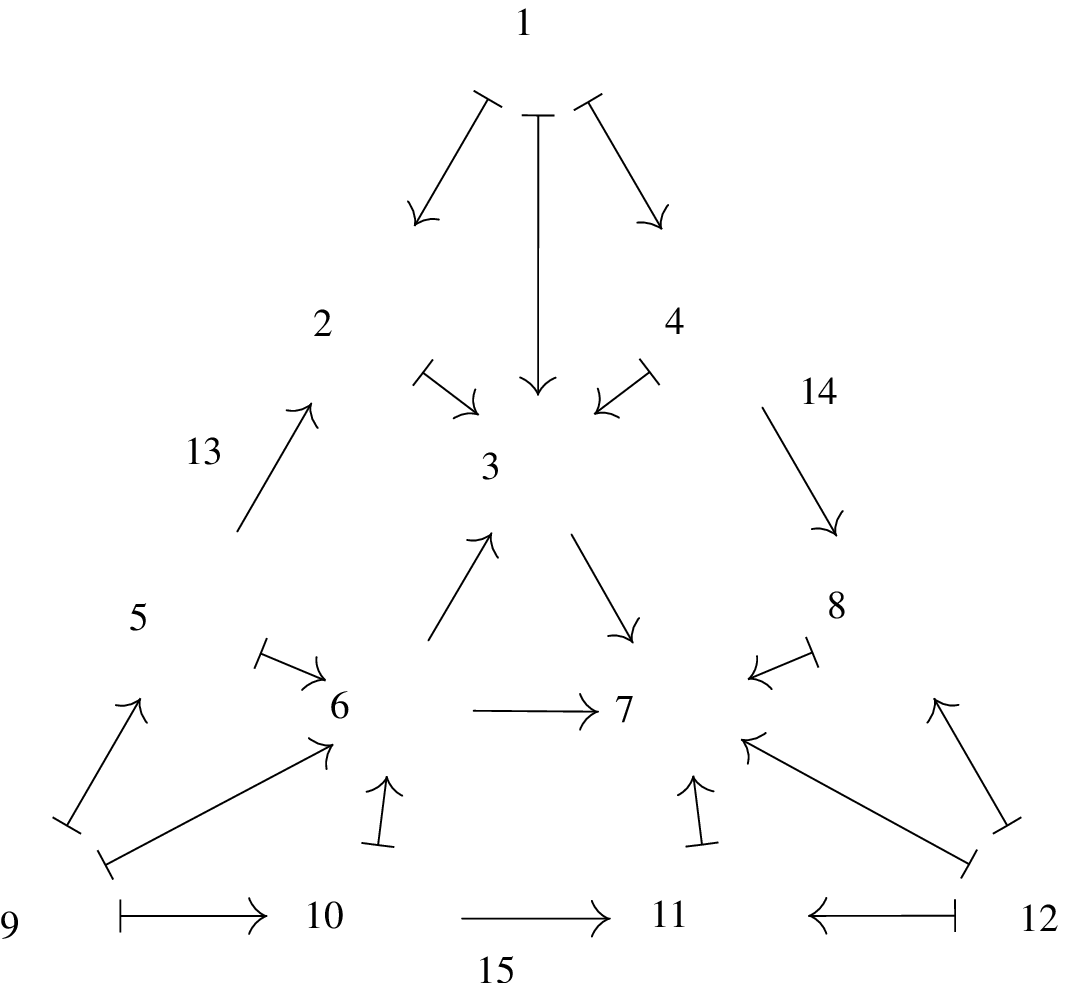,height=4cm}
\end{matrix}
\]
}
is commutative in $\hocolim\Hom(\cH_{(-)},\cG)$.

So far, we have only used the relations
$\Big(
\put(3,-2.5){\rotatebox{57}{$\scriptstyle\to$}}
\put(10,4){\rotatebox{-55}{$\scriptstyle\to$}}
\begin{matrix}
\scriptstyle\bullet\\
\scriptstyle\bullet\,\to\,\bullet
\end{matrix}\Big)$,
$\Big(
\put(1,5){\rotatebox{-90}{$\scriptscriptstyle\mapsto$}}
\put(15.5,5){\rotatebox{-90}{$\scriptscriptstyle\mapsto$}}
\begin{matrix}
\scriptstyle\bullet\,\to\,\bullet\\
\scriptstyle\bullet\,\to\,\bullet
\end{matrix}\Big)$, and
$\Big(
\put(3,-2.5){\rotatebox{57}{$\scriptstyle\mapsto$}}
\put(10,4){\rotatebox{-55}{$\scriptstyle\mapsto$}}
\begin{matrix}
\scriptstyle\bullet\\
\scriptstyle\bullet\,\mapsto\,\bullet
\end{matrix}
\Big)$.
The fourth one
$\Big(
\put(1.5,6.5){\rotatebox{237}{$\scriptstyle\mapsto$}}
\put(8,4){\rotatebox{-90}{$\scriptscriptstyle\Mapsto$}}
\put(10,4.5){\rotatebox{-57}{$\scriptstyle\mapsto$}}
\begin{matrix}
\scriptstyle\bullet\\
\scriptstyle\bullet\,\to\,\bullet
\end{matrix}\Big)$, involving the 2-morphisms of $\cov^2(\cH_0)$, is used to check that $\Psi\Phi'=1$.

We first check that $\Psi\Phi'$ is the identity when applied to the generators $({\scriptstyle\bullet\,\mapsto\,\bullet})$.
Let $V,U\to\cH_0$ be covers, $p:V\to U$ a map between them, and $\sigma$ an object of $\Hom(\cH_U,\cG)$.
We then have a corresponding generator $p_\sigma:\sigma\to p^*(\sigma)$ of $\hocolim\Hom(\cH_{(-)},\cG)$.
Applying $\Phi'$ to it, we get the morphism $\sigma\to p^*(\sigma)$ of $\M'_\cG(\cH)$ corresponding to the obvious isomorphism $\iota: P(\sigma)\to P(p^*\sigma)$.
Applying $\Psi$ to $\iota$, we then get \dontshow{ske}
\begin{equation}\label{ske}
\Psi\Phi':\Big[\sigma\put(8,7){$\scriptstyle p_\sigma$}\longmapsto p^*(\sigma)\Big]
\mapsto\Big[\sigma\put(8,7){$\scriptstyle q_\sigma$}\longmapsto q^*(\sigma)
\stackrel{\alpha_\iota}{\longrightarrow}r^*(p^*(\sigma))\stackrel{r_\sigma}{\raisebox{4.3pt}{\rotatebox{180}{$\longmapsto$}}}\tau\Big],
\end{equation}
where $q: U\times_{\cH_0}V\to U$, $r: U\times_{\cH_0}V\to V$ are the projections,
and $\alpha_\iota$ is given by $\alpha_\iota(u,v)=\sigma(u,1,p(v)):\sigma(u)\to\sigma(p(v))$.
By
$\Big(
\put(3,-2.5){\rotatebox{57}{$\scriptstyle\mapsto$}}
\put(10,4){\rotatebox{-55}{$\scriptstyle\mapsto$}}
\begin{matrix}
\scriptstyle\bullet\\
\scriptstyle\bullet\,\mapsto\,\bullet
\end{matrix}
\Big)$ and
$\Big(
\put(1.5,6.5){\rotatebox{237}{$\scriptstyle\mapsto$}}
\put(8,4){\rotatebox{-90}{$\scriptscriptstyle\Mapsto$}}
\put(10,4.5){\rotatebox{-57}{$\scriptstyle\mapsto$}}
\begin{matrix}
\scriptstyle\bullet\\
\scriptstyle\bullet\,\to\,\bullet
\end{matrix}\Big)$,
the square
{
\psfrag{1}{$\sigma$}
\psfrag{2}{$q^*\sigma$}
\psfrag{3}{$p^*\sigma$}
\psfrag{4}{$r^*p^*\sigma$}
\psfrag{5}{$\scriptstyle q_\sigma$}
\psfrag{6}{$\scriptstyle p_\sigma$}
\psfrag{7}{\rotatebox{-38}{$\scriptstyle (pr)_\sigma$}}
\psfrag{8}{$\scriptstyle \alpha_\iota$}
\psfrag{9}{$\scriptstyle r_{p^*\sigma}$}
\[
\begin{matrix}
\epsfig{file=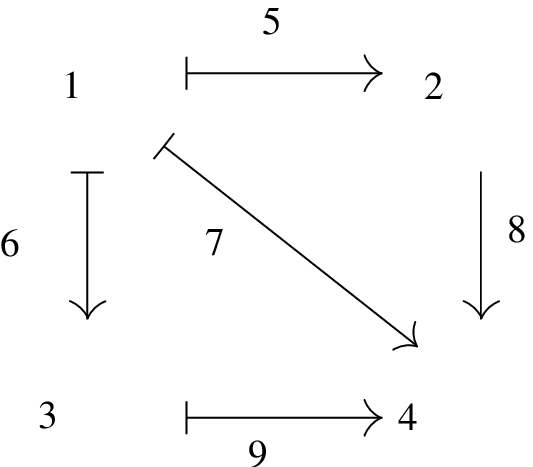,height=2cm}
\end{matrix}
\]
}
is commutative in $\colim\Hom(\cH_{(-)},\cG)$.
It follows that (\ref{ske}) is the identity mapping.

Let now $U\to\cH_0$ be a cover and $\beta:\sigma\Rightarrow\tau:\cH_U\to\cG$ an arrow in $\Hom(\cH_U,\cG)$.
We verify that $\Psi\Phi'=1$ on the corresponding generator $({\scriptstyle\bullet\,\to\,\bullet})$ of $\hocolim\Hom(\cH_{(-)},\cG)$.
Letting $f:P(\sigma)\to P(\tau)$ be the principal $\cG$-bundle map corresponding to $\beta$, we get \dontshow{sqo}
\begin{equation}\label{sqo}
\Psi\Phi':\Big[\sigma\stackrel{\beta}{\longrightarrow}\tau
\Big]
\mapsto\Big[\sigma\put(8,7){$\scriptstyle p_\sigma$}\longmapsto p^*(\sigma)
\stackrel{\alpha_f}{\longrightarrow}q^*(\tau)\stackrel{q_\tau}{\raisebox{4.3pt}{\rotatebox{180}{$\longmapsto$}}}\tau\Big],
\end{equation}
where $p,q:U\times_{\cH_0}U\to U$ are the two projections.
The transformation $\alpha_f:p^*(\sigma)\Rightarrow q^*(\tau)$ evaluated on $(u,v)\in U\times_{\cH_0}U$ is then given by the diagonal of the commutative square
\[
\xymatrix{
\sigma(u)\ar[r]^{\beta(u)}\ar[d]_{\sigma(u,1,v)}\ar[dr]|{\alpha_f(u,v)}&\tau(u)\ar[d]^{\tau(u,1,v)}\\
\sigma(v)\ar[r]_{\beta(v)}&\tau(v).
}
\]
Let $\gamma:p^*(\tau)\Rightarrow q^*(\tau)$, $\gamma(u,v):=\tau(u,1,v):\tau(u)\to\tau(v)$
be the natural transformation corresponding to the right vertical arrow in the above square,
and let $t:U\times_{\cH_0}U \to U\times_{\cH_0}U$ be the flip.
Applying
$\Big(
\put(1.5,6.5){\rotatebox{237}{$\scriptstyle\mapsto$}}
\put(8,4){\rotatebox{-90}{$\scriptscriptstyle\Mapsto$}}
\put(10,4.5){\rotatebox{-57}{$\scriptstyle\mapsto$}}
\begin{matrix}
\scriptstyle\bullet\\
\scriptstyle\bullet\,\to\,\bullet
\end{matrix}\Big)$
to the (unique) 2-morphism $1\Rightarrow t$ in $\cov^2(\cH_0)$,
and to the object $p^*(\tau)$ of $\Hom(\cH_{U\times_{\cH_0}U},\cG)$, we see that $t_{p^*\tau} =\gamma$.
By
$\Big(
\put(3,-2.5){\rotatebox{57}{$\scriptstyle\to$}}
\put(10,4){\rotatebox{-55}{$\scriptstyle\to$}}
\begin{matrix}
\scriptstyle\bullet\\
\scriptstyle\bullet\,\to\,\bullet
\end{matrix}\Big)$,
$\Big(
\put(1,5){\rotatebox{-90}{$\scriptscriptstyle\mapsto$}}
\put(15.5,5){\rotatebox{-90}{$\scriptscriptstyle\mapsto$}}
\begin{matrix}
\scriptstyle\bullet\,\to\,\bullet\\
\scriptstyle\bullet\,\to\,\bullet
\end{matrix}\Big)$, and
$\Big(
\put(3,-2.5){\rotatebox{57}{$\scriptstyle\mapsto$}}
\put(10,4){\rotatebox{-55}{$\scriptstyle\mapsto$}}
\begin{matrix}
\scriptstyle\bullet\\
\scriptstyle\bullet\,\mapsto\,\bullet
\end{matrix}
\Big)$,
the square
{
\psfrag{1}{$\sigma$}
\psfrag{2}{$p^*\sigma$}
\psfrag{3}{$\tau$}
\psfrag{4}{$q^*\tau$}
\psfrag{5}{$\scriptstyle p_\sigma$}
\psfrag{6}{$\scriptstyle \beta$}
\psfrag{7}{$p^*\tau$}
\psfrag{8}{$\scriptstyle \alpha_f$}
\psfrag{9}{$\scriptstyle q_\tau$}
\psfrag{10}{$\scriptstyle p_\tau$}
\psfrag{11}{$\scriptstyle \gamma$}
\psfrag{12}{$\scriptstyle p^*\!\beta$}
\[
\begin{matrix}
\epsfig{file=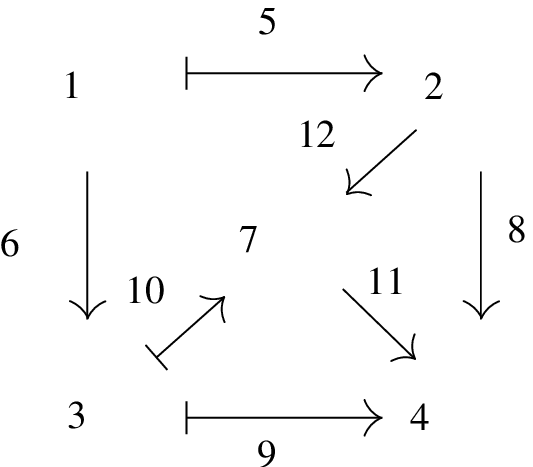,height=2cm}
\end{matrix}
\]
}
is commutative in $\colim\Hom(\cH_{(-)},\cG)$. It follows that (\ref{sqo}) is the identity.
\qed
\vspace{.3cm}

\begin{remark}
We believe that the canonical map \dontshow{nca}
\begin{equation}\label{nca}
\underset{U\in\cov(\cH_0)}{\hocolim}\,\Hom\big(\cH_U,\cG\big)
\,\longrightarrow\,\M_\cG(\cH)
\end{equation}
is always an equivalence; nevertheless, we claim that the left hand side of (\ref{nca})
is the wrong object to consider.
To illustrate our opinion, we present the following counterexample:
let $\cH:=\pair(\RR/\ZZ)$, $\cG:=*$, and let $\C\subset\cov(\RR/\ZZ)$ be the full subcategory of $\cov(\RR/\ZZ)$ whose objects are the covers
\[
U_n:=\coprod_{i=1}^{2^n}\big({\textstyle\frac{i-1}{2^n},\frac{i+1}{2^n}}\big)\to \RR/\ZZ,\qquad n\ge 0.
\]
Then the map
\[
\underset{U\in\C}{\hocolim}\,\Hom\big(\cH_U,\cG\big)
\,\longrightarrow\,\M_\cG(\cH)
\]
is {\em not} an equivalence.
\end{remark}

\proof
Clearly $\M_\cG(\cH)\simeq *$, so it suffices to show that $\hocolim_\C \Hom(\cH_{(-)},\cG)$ is not contractible.
Since $\hocolim_\C \Hom(\cH_{(-)},\cG)\simeq \hocolim_\C *$ is a connected groupoid whose automorphism groups of objects
are isomorphic to $\pi_1|\C|$, it's enough to show that $|\C|$ has non-trivial fundamental group.

Let $\D$ be the category with two objects $A,B$ and two non-identity arrows $\alpha,\bar\alpha:B\to A$.
Clearly, $\pi_1|\D|=\ZZ$.
Let $u_0,\bar u_0\in U_0$ be the two preimages of $\frac{1}{2}\in\RR/\ZZ$ and let $u_n\in U_n$, $n\ge 1$, be the unique preimage of $\frac{1}{2}\in\RR/\ZZ$.
Define a functor $F:\C\to \D$ by
\[
\begin{split}
&F(U_0)=A,\\
&F(U_n)=B\qquad\quad\text{for}\quad n\ge 1,\\
&F\big(f:U_n\to U_1\big)=\begin{cases}
\alpha &\quad\text{if}\quad f(u_n)=u_0\\
\bar\alpha &\quad\text{if}\quad f(u_n)=\bar u_0.
\end{cases}
\end{split}
\]
The map $F_*:\pi_1|\C|\to \pi_1|\D|$ is clearly surjective, from which it follows that $\pi_1|\C|\not =0$.
\qed

\subsection{Equivalent topological categories have equivalent diagram categories \rm\dontshow{A:3}}\label{A:3}

Given a topologically-enriched category $\C$, its homotopy category $\mathrm{Ho}(\C)$ is the category with objects those of $\C$ and
with arrows the path components
\[
\hom_{\mathrm{Ho}(\C)}(C,C'):=\pi_0\map_\C(C,C')
\]
of the morphism spaces of $\C$.
An arrow of $\C$ is said to be an {\em homotopy equivalence} if it becomes an isomorphism in $\mathrm{Ho}(\C)$.

If $\C$ is small, we may form the category of $\C$-spaces, with objects the continuous contravariant functors from $\C$ to spaces
and arrows their natural transformations.
We also have a category of $\C_0$-spaces, whose objects are collections of spaces indexed by the set of objects of $\C$.
Henceforth we shall write $\map_\C$ and $\map_{\C_0}$ for mapping spaces in the categories of $\C$-spaces and $\C_0$-spaces, respectively.
The forgetful functor $\{\C\text{-spaces}\}\to\{\C_0\text{-spaces}\}$ has a continuous left adjoint denoted $\FF_\C$,
or just $\FF$ when the category is clear from the context.
Following the argument in Lemma \ref{freeres}, we see that any $\C$-space $X$ is canonically the coequalizer
of the two natural maps $\FF^2 X\rrarrow\FF X$.

The Bousfield-Kan, or projective, model structure on $\C$-spaces has objectwise weak equivalences and objectwise fibrations.
Letting $\delta_C$ be the $\C_0$-space assigning the one point space to $C$ and the empty set to all other objects of $\C$,
the generating cofibrations and trivial cofibrations are then given by \dontshow{aul}
\begin{equation}\label{aul}
\FF(S^{n-1}\times\delta_C)\rightarrow \FF(D^n\times\delta_C)\qquad
\text{and}\qquad \FF(D^n\times\delta_C) \rightarrow \FF(D^n\times[0,1]\times\delta_C),
\end{equation}
respectively.
A general (trivial) cofibration is then a retract of a transfinite colimit of pushouts along generating (trivial) cofibrations
\cite[Corollary 10.5.23]{Hir03}.

\begin{definition}\label{da4}\dontshow{da4}
A continuous functor $f:\C\to\D$ between topologically enriched categories
is a {\em weak equivalence} if it is homotopically fully faithful and homotopically essentially surjective;
that is, $f$ must induce weak equivalences $\map(C,C')\to\map(fC,fC')$ for all pairs of objects $C$ and $C'$ of $\C$,
and any object $D$ of $\D$ must be homotopy equivalent to an object of the form $fC$ for some object $C$ of $\C$.
\end{definition}

\begin{lemma}\label{Aee}\dontshow{Aee}
Let $\C$ and $\D$ be small topologically enriched categories and let
$f:\C\to\D$ be a continuous functor.
Then the restriction functor $f^*$ from $\D$-spaces to $\C$-spaces admits a continuous left adjoint $f_!$ such that $(f_!,f^*)$ is a Quillen pair.
Moreover, if $f$ is a weak equivalence, then $(f_!,f^*)$ is a Quillen equivalence. % {\bf of topological model categories?}
\end{lemma}

\proof
Let $f_+:\{\C_0\text{-spaces}\}\to\{\D_0\text{-spaces}\}$ be the functor defined by the formula
\[
f_+(X)(D):=\coprod_{C\in f^{-1}(D)}X(C).
\]
That is, $f_+$ regards a $\C_0$-space as a $\D_0$-space via $f_0:\C_0\to\D_0$, and is left adjoint to the pullback functor $f^*$.
Given a free $\C$-space $\FF_\C X$, the adjunctions
$$
\map_\C(\FF_\C X,f^*Y)\cong
\map_{\C_0}(X,f^* Y)\cong
\map_{\D_0}(f_+ X,Y)\cong
\map_\D(\FF_\D f_+ X,Y)
$$
show that one must have $f_!\FF_\C X\simeq \FF_\D f_+ X$.
More generally, since a $\C$-space $X$ is always the coequalizer of $\FF_\C^2 X\rrarrow \FF_\C X$, we see that
\[
f_!X
=\coeq\big(\FF_\D f_+\FF_\C X\rrarrow \FF_\D f_+ C\big),
\]
where the two maps are induced by the left action of $\FF_\C$ on $X$ and the right action of $\FF_\C$ on $\FF_\D f_+$.

Now that we have our adjoint pair, we must show that $f^*$ preserves fibrations and trivial fibrations,
and that $f_!$ preserves cofibrations and trivial cofibrations.
Clearly, if $Y\to Y'$ is a (trivial) fibration, then, for any $C\in \C$, the map $f^*Y(C)\cong Y(fC)\to Y'(fC)\cong f^*Y'(C)$ is a (trivial) fibration of spaces,
so $f^*Y\to f^*Y'$ is a (trivial) fibration.

To check that $f_!$ preserves cofibrations and trivial cofibrations it's enough, as $f_!$ preserves colimits and retracts,
to show that $f_!$ takes the generating (trivial) cofibration to (trivial) cofibrations.
Indeed, letting $i$ denote the generating (trivial) cofibration $i: S^{n-1}\to D^n$ or $i:D^n\to D^n\times [0,1]$ for the usual model structure on topological spaces,
we see that
\[
f_!\FF_\C(i\times \delta_C)\cong \FF_\D f_+ (i\times \delta_C)\cong \FF_\D(i\times \delta_{fC}).
\]
It follows that $(f_!,f^*)$ is a Quillen pair.

Assuming now that $f$ is a weak equivalence, we show that $(f_!,f^*)$ is a Quillen equivalence.
Namely, if $X$ is a cofibrant $\C$-space and $Y$ is a fibrant $\D$-space (that is, an arbitrary $\D$-space),
we show that a $\C$-space map $X\to f^*Y$ is an equivalence if and only if its adjoint $f_!X\to Y$ is an equivalence.
Given a $\C$-space $X$ and an object $C\in\C$, let us first examine under what conditions the map \dontshow{crs}
\begin{equation}\label{crs}
X(C)\cong \map_\C(\FF_\C\delta_C,X)\stackrel{\textstyle f_!}{\longrightarrow}\map_\D(f_!\FF_\C\delta_C,f_!X)\cong \map_{\D_0}(\delta_{fC},f_!X)\cong f_!X(fC)
\end{equation}
is a weak equivalence.
If $X=\FF_\C\delta_{C'}$, then the above map can be written as
\[
\FF_\C\delta_{C'}(C)\cong\map(C,C')\stackrel{\textstyle f}{\longrightarrow}
\map(fC,fC')\cong \FF_\D\delta_{fC'}(fC)\cong f_!\FF_\C\delta_{C'}(fC),
\]
and is a weak equivalence since $f$ is homotopically fully faithful.
A cofibration $X\to X'$ of $\C$-spaces induces Hurewicz cofibrations on both sides of (\ref{crs}).
As the class of $\C$-spaces $X$ for which (\ref{crs}) is an equivalence is closed under
the operations of taking products with disks, pushing out along generating cofibrations, transfinite filtered colimits of these pushout maps, and retracts,
it follows that (\ref{crs}) is a weak equivalence for every cofibrant $\C$-space $X$.

Now suppose $X$ is a cofibrant $\C$-space and $Y$ is a fibrant $\D$-space.
Given a map $X\to f^*Y$, with adjoint $f_!X\to Y$, and an object $C$ in $\C$, we get a commutative triangle
\[
\xymatrix{
&f_!X(fC)\ar[dr]&\\
X(C)\ar[ur]^{\sim}\ar[rr]&&Y(fC)
}
\]
in which the left-most arrow is always an equivalence.

By the two out of three property, $X(C)\to Y(fC)$ is a weak equivalence for all $C\in\C$
if and only if $f_!(fC)\to Y(fC)$ is a weak equivalence for all $C\in\C$.
Since $f$ is homotopically essentially surjective, the latter property is equivalent to
$f_!(D)\to Y(D)$ being a weak equivalence for all $D\in\D$.
Hence $f_! Y\to X$ is a weak equivalence.
\qed

\subsection{Some point-set topology}

\begin{lemma}\label{ytb}\dontshow{ytb}
Let $A$, $A'$, $B$, and $B'$ be objects of a topologically enriched category, and suppose that we are given homotopy equivalences
$f:A\to A'$, $g:B\to B'$ and maps $\varepsilon:S^{n-1}\times B\to A$, $\varepsilon':S^{n-1}\times B'\to A'$ satisfying $f\circ \varepsilon=\varepsilon'\circ(1\times g)$.
Then the resulting map
\[
f\cup(1\times g):A\cup_\varepsilon(D^n\times B)\longrightarrow A'\cup_{\varepsilon'}(D^n\times B')
\]
is a homotopy equivalence.
Moreover, the homotopy inverse of $f\cup(1\times g)$ and the two homotopies between the composites and the identities can be taken
compatibly with those for $f$.
\end{lemma}

\proof
We just construct the homotopy inverse.
Let $f':A'\to A$, $g':B'\to B$ be homotopy inverses of $f$ and $g$, and let $h:[0,1]\times S^{n-1}\times B'\to A$
be the homotopy between $\varepsilon\circ(1\times g')$ and $f'\circ\varepsilon'$ obtained by composing the homotopies
\[
\varepsilon\circ(1\times g') \sim f'\circ f\circ\varepsilon\circ(1\times g')=
f'\circ \varepsilon'\circ(1\times g)\circ(1\times g') \sim f'\circ\varepsilon'.
\]
Then
\[
A'\cup_{\varepsilon'}(D^n\times B')\cong
A'\cup_{\varepsilon'}([0,1]\times S^{n-1}\times B')\cup (D^n\times B')
\put(4,6){$\scriptstyle f'\cup h\cup (1\times g')$}
\:{-\!\!\!-\!\!\!-\!\!\!-\!\!\!-\!\!\!-\!\!\!-\!\!\!
\longrightarrow}\:A\cup_\varepsilon(D^n\times B).
\]
is the desired homotopy inverse of $f\cup (1\times g)$.
\qed
\vspace{.3cm}

For this last lemma we write $+$ in place of $\sqcup$ for better readability.

\begin{lemma}\label{pushpull}
Consider the following commutative diagram of topological spaces
$$
\xymatrix{
U_1 \ar[d]\ar[r] & Y_1 \ar[d] & Z_1 \ar[l]\ar[d]\ar[r] & X_1\ar[d]\\
U_3       \ar[r] & Y_3        & Z_3 \ar[l]      \ar[r] & X_3      \\
U_2 \ar[u]\ar[r] & Y_2 \ar[u] & Z_2 \ar[l]\ar[u]\ar[r] & X_2\ar[u]}
$$
in which the maps $Z_i\hookrightarrow Y_i$ are closed inclusions with complementary open inclusions $U_i\hookrightarrow Y_i$.
Write
$$
P_U := U_1\times_{U_3} U_2,\qquad P_Y := Y_1\times_{Y_3} Y_2,\qquad P_Z := Z_1\times_{Z_3} Z_2,\qquad P_X:= X_1\times_{X_3} X_2
$$
for the pullbacks of the $U_i, Y_i, Z_i, X_i$, respectively, and
$$
Q_i := Y_i +_{Z_i} X_i
$$
for the pushout of the closed inclusion $Z_i\hookrightarrow Y_i$ along $Z_i\to X_i$.
Then the natural map
$$
P_Y +_{P_Z} P_X\to Q_1\times_{Q_3} Q_2,
$$
from the pushout of the column-wise pullbacks to the pullback of the row-wise pushouts, is an isomorphism.
\end{lemma}

\proof
The continuous bijections $U_i + Z_i\to Y_i$ induce continuous bijections $U_i + X_i\to Q_i$; similarly, the continuous bijection $P_U + P_Z\to P_Y$ induces a continuous bijection $P_U + P_X\to P_Y +_{P_Z} P_X$.
It follows that we have a continuous bijection $P_U + P_X\to Q_1\times_{Q_3} Q_2$, and therefore that the map $P_Y +_{P_Z} P_X\to Q_1\times_{Q_3} Q_2$ is a continuous bijection.

Consider the commutative diagram
$$
\xymatrix{
P_Y + P_X \ar[d]\ar[r] & P_Y + (Z_1\times_{X_3} X_2) + (X_1\times_{X_3} Z_2) + P_X\ar[d]\ar[r] & (Y_1 + X_1)\times (Y_2 + X_2)\ar[d]\\
P_Y +_{P_Z} P_X \ar[r] & Q_1\times_{Q_3} Q_2                                        \ar[r] & Q_1\times Q_2}
$$
in which (as one can easily check using the various continuous bijections mentioned above) the right-hand square is a pullback, the vertical maps are quotient maps, and the horizontal maps --- save for possibly the map $P_Y +_{P_Z} P_X\to Q_1\times_{Q_3} Q_2$ in question --- are closed inclusions.
Since the map
$$
(Z_1\times_{X_3} X_2) + (X_1\times_{X_3} Z_2)\to X_1\times_{X_3} X_2 = P_X
$$
induces a retraction
$$
P_Y + (Z_1\times_{X_3} X_2) + (X_1\times_{X_3} Z_2) + P_X\to P_Y + P_X
$$
of the upper-left-most closed inclusion, we see that $P_Y +_{P_Z} P_X$ is a quotient of $P_Y + Z_1\times_{X_3} X_2 + X_1\times_{X_3} Z_2 + P_X$.
Since $Q_1\times_{Q_3} Q_2$ is a quotient of the same space, it follows that the continuous bijection $P_Y +_{P_Z} P_X\to Q_1\times_{Q_3} Q_2$ is a homeomorphism.
\qed

\bibliography{../main}
\bibliographystyle{plain}
\end{document}